\documentclass
[10pt,a4paper]{article}

\usepackage{amsmath,amsfonts,amssymb,amsthm,cite,graphics,graphicx}
\usepackage{txfonts}
\usepackage{mathrsfs}
\usepackage[colorlinks=true]{hyperref}
\usepackage{tabularx}
\hypersetup{urlcolor=blue, citecolor=red}
\usepackage{hyperref}
\usepackage{geometry}
\newtheorem{cor}{Corollary}[section]
\newtheorem{thm}{Theorem}[section]
\newtheorem{lem}{Lemma}[section]
\newtheorem{rem}{Remark}[section]
\newtheorem{prop}{Proposition}[section]
\newtheorem{defi}{Definition}[section]
\theoremstyle{definition}

\begin{document}
\setlength{\topmargin}{-0.5in} \setlength{\baselineskip}{0.22in}

\title{$H^{s,p}$ estimates for weak solutions of fractional nonlocal equations in divergence form }

\maketitle
\maketitle

\centerline{\scshape Zhenjie Li }
\medskip
{\footnotesize
 \centerline{357836376@qq.com}
  \centerline{School of Mathematical Sciences,
   Guangxi Minzu University, Nanning, 530006,  China}

}
\centerline{\scshape Lihe Wang }
\medskip
{\footnotesize
  \centerline{lihe-wang@uiowa.edu}
  \centerline{Department of Mathematics, University of Iowa, Iowa City, IA, 52242, USA;}

}
\centerline{\scshape Chunqin Zhou }
\medskip
{\footnotesize
 \centerline{cqzhou@sjtu.edu.cn}
 \centerline{School of Mathematical Sciences,
 Shanghai Jiao Tong University, Shanghai, 200240, China}

}

\begin{center}
\parbox{0.8\hsize}{\small
\noindent{\bf Abstract}
 This paper establishes interior  $H^{s,p}$ regularity for weak solutions to fractional nonlocal equations in divergence form with general kernels. Under the $(\delta, R)$-vanishing condition on the coefficients, we prove that for $p>2$, if the source term $f\in L^p_{loc}$, then the solution $u\in H^{s,p}_{loc}$. In contrast to Byun and Kim \cite{BK2025}, our result holds in the borderline case for $f$ without requiring differentiability assumptions. Our approach simplifies previous methods by using the stability properties of the fractional gradient to derive the good-inequality directly, thereby avoiding reliance on maximal functions. The proof integrates Vitali's covering lemma, Hardy-Littlewood maximal theory, self-improving estimates, and compactness arguments. This work extends classical $W^{1,p}$ regularity to fractional nonlocal settings under minimal kernel assumptions.

\
\

}

\end{center}

\pagestyle{myheadings} \thispagestyle{plain} \markboth{ }{}

\section{Introduction}

In this paper, we study the interior $H^{s,p}$ regularity for weak solutions to the following fractional nonlocal equation:
\begin{equation}\label{1.1}
 L^{2s}_Ku(x)=L^{s}_bf(x),\quad \hbox{ in }\Omega\subset\mathbb{R}^n,
\end{equation}
where $0<s<1$, $n>2s$ and $\Omega$ is a bounded domain. It is well known that for any  functions $u, f\in\mathcal{S}(\mathbb{R}^n)$, the nonlocal operators are defined as follows:
\begin{eqnarray}
L^{2s}_Ku(x)&=&C(n,s)P.V.\int_{\mathbb{R}^n}K(x,y)[u(x)-u(y)]dy\nonumber\\
&=&C\lim_{\epsilon\rightarrow0}\int_{\mathbb{R}^n\setminus B_\epsilon(x)}K(x,y)[u(x)-u(y)]dy,\nonumber
\end{eqnarray}
and
\begin{eqnarray}
L^{s}_bf(x)&=&C\left(\frac{n}{2},s\right)P.V.\int_{\mathbb{R}^n}b(x,y)[f(x)-f(y)]dy\nonumber\\
&=&C\lim_{\epsilon\rightarrow0}\int_{\mathbb{R}^n\setminus B_\epsilon(x)}b(x,y)[f(x)-f(y)]dy,\nonumber
\end{eqnarray}
where $\mathcal{S}(\mathbb{R}^n)$ is the Schwartz space of rapidly decaying $C^\infty$ functions in $\mathbb{R}^n$, the symbol P.V. indicates integration in the sense of Cauchy principal value, and
$$
C(n,s)=\left(\int_{\mathbb{R}^n}\frac{1-\cos(\zeta_1)}{|\zeta|^{n+2s}}d\zeta\right)^{-1},
$$
$K(x,y)$ and $b(x,y)$ are integral kernels. In what follows, we always assume that the kernels satisfy
\begin{equation}\label{1.2}
K(x,y)=\frac{A(x,y)}{|x-y|^{n+2s}}\quad\hbox{ with }\quad\lambda\leq A(x,y)\leq\Lambda,  \\
\end{equation}
\begin{equation}\label{1.3}
|b(x,y)|\leq\frac{\beta}{|x-y|^{n+s}},\\
\end{equation}
and
\begin{equation}\label{1.3.1}
\|L^s_bv\|_{L^2(\mathbb{R}^n)}\leq C(n,s,\beta)[v]_{H^s(\mathbb{R}^n)},\quad\forall v\in H^s(\mathbb{R}^n),
\end{equation}
for some positive constants $\lambda<\Lambda$ and $\beta$.
Note that (\ref{1.3.1}) holds for any translation-invariant kernel $b(x,y)=b(x-y)$ satisfying (\ref{1.3}). Throughout this paper, we  always assume  kernels $K$
and $b$  satisfy  assumptions  (\ref{1.2})-(\ref{1.3.1}).

We also can define the weak solution to (\ref{1.1}). We say $u\in H^{s}(\mathbb{R}^n)$ is a weak solution of (\ref{1.1}) if
\begin{equation}
\int_{\mathbb{R}^n}\int_{\mathbb{R}^n}K(x,y)[u(x)-u(y)][\phi(x)-\phi(y)]dydx=\int_{\mathbb{R}^n}f(x)L^s_b\phi(x)dx, \quad\forall\phi\in \mathcal{H}^s_0(\Omega),\nonumber\\
\end{equation}
where  $f\in L^2(\mathbb{R}^n)$.

$W^{1,p}$ estimates play a crucial role in the regularity theory of partial differential equations. The fractional counterpart,  $W^{s,p}$ estimates, has been investigated  by many mathematicians. Nowak\cite{N2020} studied weak solutions of  a nonlocal elliptic equation
\begin{equation}
 L^{2s}_Au(x)+bu=L^{2s}_Dg(x)+f(x),\quad \hbox{ in }\Omega,\nonumber\\
\end{equation}
and established the local $W^{s,p}$ regularity under the assumption that the kernel $A$ is elliptic, symmetric and translation invariant. The approach, initially introduced in \cite{CP1998} and further extended in \cite{W2003,BW2004},   relies on finding an appropriate approximation equation whose  weak solution admits  local $C^{0,1}$ estimates.  Consequently,  the desired local $W^{1,p}$ estimates can be obtained via the compactness method, the Hardy-Littlewood maximal function and  Vitali's covering Lemma. Moreover, Dong and Kim in \cite{DK2012} obtained a similar regularity result for the equation
$$
L^{2s}_Au-\lambda u=f,
$$
when $A$ is bounded and translation invariant. For further results in the nonlocal setting, we can  refer to \cite{N2021-1}, where the local $W^{s,p}$ estimate was extended to more general nonlinear nonlocal equations using the characterization of Bessel potential spaces and a perturbation approach. Additionally,  Biccari, Warma and Zuazua \cite{BWZ2017}  established the local $W^{2s,p}$ regularity for weak solutions to the homogeneous Dirichlet problem for the fractional Laplacian equation
\begin{equation}
\left\{
\begin{array}{lll}
(-\Delta)^su(x)=f(x),\quad &\hbox{ in }\Omega,\nonumber\\
 u=0,\quad &\hbox{ in }\mathbb{R}^n\setminus\Omega,\nonumber
\end{array}\right.
\end{equation}
on any open bounded domain $\Omega$, through a cut-off argument.
Cozzi\cite{C2017} derived  local $H^{2s-\epsilon,p}$ estimates for weak solutions of the nonlocal equation
\begin{equation}
L^{2s}_Ku(x)=f(x),\quad \hbox{ in }\Omega,\nonumber
\end{equation}
assuming  the kernel $K(x,y)$ is symmetric, globally comparable to the kernel of the fractional Laplacian and locally uniformly H\"{o}lder continuous, employing Nirenberg's translation method. By a similar process, higher differentiability and integrability in the nonlinear setting were proved in \cite{KMS2015}.
Byun and Kim investigated the fractional p-Laplacian-type equations
with discontinuous kernel coefficients
in divergence form in \cite{BK2025},
\begin{equation}
(-\Delta)^s_Au=(-\Delta)^{\frac{s}{p}}f\quad\text{ in }\Omega,\quad
\end{equation}
where $2\leq p<\infty$, and  established $W^{s+\sigma,q}$ estimates for any choice of pairs $(\sigma, q)$ with
$q\in(p,\infty)$ and $\sigma\in (0, \min\{\frac{s}{p-1}, 1-s\})$ under the assumption that $f\in W^{\sigma,q}$  and the associated
kernel coefficients have small BMO seminorms near the diagonal.

In this paper, we establish the  interior  $H^{s,p}$ estimates for weak solutions to the fractional nonlocal equation (\ref{1.1}).  Unlike Byun and Kim \cite{BK2025}, our result holds in the borderline case for $f$ without imposing differentiability assumptions. To the best of our knowledge, this problem is novel, as there are few existing results regarding $H^{s,p}$ regularity of equations of this type. Indeed, equation (\ref{1.1})  can be viewed as  a  general nonlocal version of  second-order elliptic equations in  divergence form, such as
\begin{equation}\label{1.4}
 \operatorname{div}(A\nabla u)= \operatorname{div}(\mathbf{f}).
\end{equation}
For this equation, Calder\'{o}n-Zygmund \cite{CZ1952}  used singular integral techniques to establish the following a priori estimate:
$$
\|\nabla u\|_{L^p(\Omega)}\leq C\|\mathbf{f}\|_{L^p(\Omega)} \quad\text{for any }p\in(1,\infty),
$$
 when $A=I_n$, where $C$ is independent of $u$ and $\mathbf{f}$. This estimate was later extended to the case where $A$ is uniformly continuous in \cite{M1966,S1972,C1980}.  Meyers in \cite{M1963}  demonstrated that this estimate fails for  certain discontinuous coefficients by providing a counterexample. Later, based on representation formulas involving singular integral operators and commutators, Di Fazio \cite{F1996} deduced this a priori estimate  for the corresponding Dirichlet problem  of equation (\ref{1.4}) under the assumption that $A\in VMO$ and $\partial\Omega\in C^{1,1}$. This result was subsequently weakened  to  include real and complex non-symmetric operators and $C^1$ boundary via an appropriate representation for the Green function in \cite{AQ2002}. Using  an approximation  method, Caffarelli and Peral \cite{CP1998} provided an alternative proof of $W^{1,p}$ regularity for linear elliptic equations  under the assumption that $\|A-I_n\|_{\infty}\leq\delta$. This was further extended under weaker assumptions on the coefficient $A$ and the regularity of the domain $\Omega$ in \cite{BW2004,W2003}. In particular, Byun and Wang \cite{BW2004} derived the optimal regularity requirements on the coefficients and domains for  the  $W^{1,p}$ regularity for the Dirichlet problem of equation (\ref{1.4}) using an  approach involving Vitali's covering Lemma, the compactness method and the Hardy-Littlewood maximal function.

In this paper, we remove the symmetry and translation-invariance assumptions in \cite{N2020,DK2012}. Our proof is simpler in that we exploit the integral properties of $\nabla^s_{\mathbb{R}^n}u(x)$ directly, rather than rely on $\mathcal{M}(\nabla^s_{B}u)(x)$ as in \cite{N2020}.

\begin{defi}
We say the function $A(x,y)$ is $(\delta,R)$-vanishing if
$$
\Big(\fint_{B_r(x_0)}\fint_{B_r(x_0)}[A(x,y)-\bar{A}_{B_r(x_0)\times B_r(x_0)}]^2dydx\Big)^{\frac{1}{2}}\leq\delta\quad\hbox{ for any }0<r<R,\hbox{ and  any }x_0\in\mathbb{R}^n,
$$
where $\bar{A}_{B_r(x_0)\times B_r(x_0)}=\fint_{B_r(x_0)}\fint_{B_r(x_0)}A(x,y)dydx$.
\end{defi}

The main result is
\begin{thm}Let $s\in(0,1)$. Assume that the kernels  $K$ and $b$ satisfying (\ref{1.2})-(\ref{1.3.1}), $f\in L^2(\mathbb{R}^n)$, and $u\in {H}^s(\mathbb{R}^n)$ is a weak solution of (\ref{1.1}) in $\Omega$.  Then there exists  a small $\delta=\delta(n,s,\lambda,\Lambda,\beta,p)$ such that if $B_8\subset\Omega$, $A(x,y)$ is $(\delta, 8)$-vanishing and  $f\in L^p(B_8)$ with $p\in (2, +\infty)$, then $u\in H^{s,p}(B_1)$ and
\begin{equation}
\|\nabla^s_{\mathbb{R}^n}u\|_{L^p(B_1)}\leq C\Big[\|u\|_{L^2(B_8)}+\Big(\int_{\mathbb{R}^n\setminus B_8}\frac{u^2(x)}{|x|^{n+2s}}dx\Big)^{\frac{1}{2}}+\|f\|_{L^p(B_8)}+\Big(\int_{\mathbb{R}^n\setminus B_8}\frac{|f|^2(x)}{|x|^{n+s}}dx\Big)^{\frac{1}{2}}\Big],\nonumber\\
\end{equation}
 where $C$ is a positive constant independent of $u$ and $f$.
\end{thm}

\begin{rem}
Our approach is  heavily influenced by \cite{BW2004} and  the deep, delicate  analysis in \cite{W2003}, where  the Hardy-Littlewood maximal function was used to overcome the instability of $|\nabla^2u(x)|\leq1$ in the setting of $W^{2,p}$ theory when the authors attempted to obtain the so-called good-$\lambda$ inequality. Unlike the case of $\nabla^2 u$, we observe that the condition $\nabla^s_{\mathbb{R}^n}u(x)\leq1$  possesses  certain stability properties, which enable us to derive the good-$\lambda$ condition more directly through $\nabla^s_{\mathbb{R}^n}u(x)$ instead of $\mathcal{M}(\nabla^s_{\mathbb{R}^n}u)(x)$. The main tools employed in this paper are the Vitali covering lemma and the Hardy-Littlewood maximal function. Our basic idea is to seek  an appropriate approximation equation for which the corresponding weak solution has an interior $C^\alpha$ ($\alpha>s$) estimate,  thereby allowing us to obtain the desired $H^{s,p}$ estimate via a modified fractional version  of the argument in  \cite{BW2004,W2003}. It turns out that the fractional version of the argument is more delicate and complicated than  the differential one.
\end{rem}

The remaining sections of this  paper are organized as follows. In Section 2, we introduce the  basic notation for fractional Sobolev spaces and some  preliminary tools. In Section 3, we derive the interior H\"{o}lder estimate,  which will be used in the following section. We establish a higher differentiability  and integrability result in Section 4. Finally, we establish the interior   $H^{s,p}$ estimate in Section 5.

\section{Notations and Preliminary Tools}\setcounter{equation}{0}

\subsection{Notations for Function Space}\setcounter{equation}{0}
Let $\Omega$ be a domain in $\mathbb{R}^n$. The fractional Sobolev space $W^{s,p}(\Omega)$,  also called Sobolev-Slobodeckij space, is defined as follows
 $$
 {W}^{s,p}(\Omega)=\Big\{u\in L^p(\Omega)\Big|\hbox{ }\frac{|u(x)-u(y)|}{|x-y|^{\frac{n}{p}+s}}\in L^p(\Omega\times\Omega)\Big\},\quad 0<s<1,\quad 1\leq p<\infty,
 $$
which is an  intermediary Banach space  between $L^p(\Omega)$ and $W^{1,p}(\Omega)$, endowed  with the norm
$$
\|u\|_{W^{s,p}(\Omega)}=\Big(\|u\|^p_{L^p(\Omega)}+[u]^p_{W^{s,p}(\Omega)}\Big)^{\frac{1}{p}},
$$
where
$$
[u]_{W^{s,p}(\Omega)}=\Big(\int_{\Omega}\int_{\Omega}\frac{|u(x)-u(y)|^p}{|x-y|^{n+ps}}dxdy\Big)^{\frac{1}{p}}
$$
is the Gagliardo seminorm of $u$. Let  $W^{s,p}_0(\Omega)$ be the closure of $C^\infty_0(\Omega)$ in the norm $\|\cdot\|_{W^{s,p}(\Omega)}$, then $W^{s,p}_0(\mathbb{R}^n)=W^{s,p}(\mathbb{R}^n)$.
For every open and bounded domain $\Omega\subset\mathbb{R}^n$, define $\mathcal{W}^{s,p}_0(\Omega)$ as the closure of $C^\infty_0(\Omega)$ with respect to the norm
$$
\|u\|_{\mathcal{W}^{s,p}_0(\Omega)}=\Big(\int_{\mathbb{R}^n}\int_{\mathbb{R}^n}\frac{|u(x)-u(y)|^p}{|x-y|^{n+ps}}dxdy\Big)^{\frac{1}{p}}.
$$
In fact, the space $\mathcal{W}^{s,p}_0(\Omega)$ can be equivalently defined as the closure of $C^\infty_0(\Omega)$ with respect to the norm
$$
\Big(\int_{\Omega}|u|^pdx\Big)^{\frac{1}{p}}+\Big(\int_{\mathbb{R}^n}\int_{\mathbb{R}^n}\frac{|u(x)-u(y)|^p}{|x-y|^{n+ps}}dxdy\Big)^{\frac{1}{p}}.
$$
If  $sp\neq1$  and $\Omega$ is an open bounded Lipschitz set in $\mathbb{R}^n$, the space $\mathcal{W}^{s,p}_0(\Omega)$ coincides with $W^{s,p}_0(\Omega)$, see \cite{BLP2014}.
For the case $sp=1$, the strict inclusion holds, i.e., $\mathcal{W}^{s,p}_0(\Omega)\subsetneqq W^{s,p}_0(\Omega)$, see\cite{BP2016}.

For any $0<s<1$ and $1\leq p<\infty$, the Bessel potential space $H^{s,p}(\mathbb{R}^n)$ is defined as
$$
H^{s,p}(\mathbb{R}^n)=\Big\{u\in L^p(\mathbb{R}^n)|\hbox{ }\mathscr{F}^{-1}[(1+|\xi|^{2})^{\frac{s}{2}}\mathscr{F}u]\in L^p(\mathbb{R}^n)\Big\},
$$
with the norm
$$
\|u\|_{H^{s,p}(\mathbb{R}^n)}=\Big\|\mathscr{F}^{-1}[(1+|\xi|^{2})^{\frac{s}{2}}\mathscr{F}u]\Big\|_{L^p(\mathbb{R}^n)},
$$
where $\mathscr{F}$ and $\mathscr{F}^{-1}$ denote the Fourier transform and the inverse Fourier transform, respectively. For any open set $\Omega\subset\mathbb{R}^n$, $H^{s,p}(\Omega)=\{u|_\Omega\hbox{ }|\hbox{ }u\in H^{s,p}(\mathbb{R}^n)\}$ with the norm $\|u\|_{H^{s,p}(\Omega)}=\inf\{\|v\|_{H^{s,p}(\mathbb{R}^n)}\hbox{ }|\hbox{ }v|_\Omega=u\}$. When $p=2$, $W^{s,2}(\mathbb{R}^n)=H^{s,2}(\mathbb{R}^n)\triangleq H^{s}(\mathbb{R}^n)$ \cite{NPV2012}.
 For any  open bounded domain $\Omega\subset\mathbb{R}^n$, we denote $\mathcal{W}^{s,2}_0(\Omega)=\mathcal{H}^s_0(\Omega)$, which is a Hilbert space with the scalar product
$$
<u,v>_{\mathcal{H}^s_0(\Omega)}=\int_{\mathbb{R}^n}\int_{\mathbb{R}^n}\frac{[u(x)-u(y)][v(x)-v(y)]}{|x-y|^{n+2s}}dxdy,
$$
and corresponding norm $$\|u\|_{\mathcal{H}^s_0(\Omega)}=[u]_{H^s(\mathbb{R}^n)}
=\Big(\int_{\mathbb{R}^n}\int_{\mathbb{R}^n}\frac{[u(x)-u(y)]^2}{|x-y|^{n+2s}}dxdy\Big)^{\frac{1}{2}}.$$
We recall that
$
\mathcal{H}^s_0(\Omega)=\Big\{u\in H^{s,2}(\mathbb{R}^n)|\hbox{ }  u=0 \hbox{ in }\mathbb{R}^n\setminus\Omega\Big\}
$ for any bounded domain with continuous boundary, see \cite{G1985}.\\

Let $\Omega\subset\mathbb{R}^n$ be an open set. For any measurable function $u$: $\mathbb{R}^n\rightarrow\mathbb{R}$, define
$$
\nabla^s_\Omega u(x)=\Big(\int_{\Omega}\frac{[u(x)-u(y)]^2}{|x-y|^{n+2s}}dy\Big)^\frac{1}{2}.
$$

\begin{rem}
Let $\Omega$ be a bounded $C^{0,1}$ domain or $\Omega=\mathbb{R}^n$. For any $s\in(0,1)$ and $p\in(\frac{2n}{n+2s},\infty)$, $u\in H^{s,p}(\Omega)$ if and only if $u$, $\nabla^s_\Omega u\in L^p(\Omega)$,
 see\cite{N2020,S1961,PS2017,P2015}.

\end{rem}

\begin{lem}(\cite{CC1995})
Suppose $f\geq0$ is a measurable function in a bounded  domain $\Omega$. Let $\theta>0$ and $m>1$ be constants. Then for any $p\in(0,\infty)$,
$$
f\in L^p(\Omega)\quad\hbox{ iff } \quad S=\sum_{k\geq1}m^{kp}|\{x\in\Omega\hbox{ }| \hbox{ }f(x)>\theta m^k\}|<\infty,
$$
and
$$
\frac{1}{C}S\leq\|f\|^p_{L^p(\Omega)}\leq C(|\Omega|+S),
$$
 where $C>0$ is a constant depending only on $\theta$, $m$, and $p$.
\end{lem}
\subsection{Preliminary Tools}\setcounter{equation}{0}

(1) The Hardy-Littlewood maximal function.

 Let $f$ be a locally integrable function. Then we say
$$
(\mathcal{M}f)(x)=\sup_{R>0}\frac{1}{|B_R(x)|}\int_{B_R(x)}|f(y)|dy.
$$
is the Hardy-Littlewood maximal function of $f$.

For any domain $\Omega\subset\mathbb{R}^n$, let $f\in L^1(\Omega)$ and define $\mathcal{M}_{\Omega}f=\mathcal{M}(f\chi_\Omega)$.
\begin{thm}(\cite{N2020})

(i) (weak 1-1 estimate) If $f\in L^1(\Omega)$ and $t>0$, then
$$
|\{x\in\Omega|\quad\mathcal{M}_\Omega(f)(x)>t\}|\leq\frac{C}{t}\int_{\Omega} |f|dx,
$$
where $C=C(n)$.

(ii) (strong $p$-$p$ estimate) If $f\in L^p(\Omega)$ with $p>1$, then $\mathcal{M}_\Omega f\in L^p(\Omega)$. Moreover,
$$
\|f\|_{L^p({\Omega})}\leq\|\mathcal{M}_\Omega f\|_{L^p({\Omega})}\leq C\|f\|_{L^p({\Omega})},
$$
where $C=C(n,p)$.
\end{thm}

\noindent (2) Vitali's covering Lemma.

\begin{lem}(\cite{S1993})
Let $E$ be a measurable set. Suppose that  a class of balls $B_\alpha$ covers $E$, i.e.,
$E\subset\bigcup_\alpha B_\alpha$.
Suppose the radius of $B_\alpha$ is bounded from above. Then there exists disjoint $\{B_{\alpha_i}\}^\infty_{i=1}\subset\{B_\alpha\}_\alpha$ such that
$E\subset\bigcup_i 5B_{\alpha_i}$,
where $5B_{\alpha_i}$ is the ball with five times the radius of $B_{\alpha_i}$. Consequently, we have
$|E|\leq 5^n\sum_i|B_{\alpha_i}|$.
\end{lem}

In order to get the interior $H^{s,p}$ estimate, we will use the following version of the Vitali covering lemma.

\begin{lem}(\cite{W2003})
Assume that $C$ and $D$ are measurable sets, $C\subset D\subset B_1$, and  that
 there exists $\epsilon>0$ such that
$|C|< \epsilon|B_1|$, and for all $x\in B_1$ and for all $r\in(0,1]$ with $|C\cap B_r(x)|\geq\epsilon|B_r(x)|$,
$B_r(x)\cap B_1\subset D$.
Then $|C|\leq 10^n\epsilon|D|$.
\end{lem}

\noindent (3) Fractional Sobolev-Poincar\'{e} inequalities
\begin{lem}[fractional Poincar\'{e} inequality\cite{KMS2015,EE2022}]
Let $s\in(0,1)$, $1\leq p<\infty$, $B\subset\mathbb{R}^n$ be a ball of radius $r$ and $w\in W^{s,p}(B)$. Then the following inequality holds:
\begin{equation}
\int_B|w-\bar{w}_B|^p\leq c(n,s,p)r^{ps}\int_B\int_B\frac{|w(x)-w(y)|^p}{|x-y|^{n+ps}}dydx.\nonumber
\end{equation}
\end{lem}

\begin{lem}[fractional Sobolev-Poincar\'{e} inequality\cite{KMS2015}]
(1)Let $s\in(0,1)$, $1\leq p<\infty$, $w\in W^{s,p}(B)$ where $B\subset\mathbb{R}^n$ is a ball of radius $r$. If $n>ps$, then the following inequality holds:
\begin{equation}
\Big(\int_B|w-\bar{w}_B|^{p^*}\Big)^{\frac{1}{p^*}}\leq C(n,s,p)\Big(\int_B\int_B\frac{|w(x)-w(y)|^p}{|x-y|^{n+ps}}dydx\Big)^{\frac{1}{p}},\nonumber\\
\end{equation}
 where $p^*=np/(n-ps)$.
\end{lem}

\section{Interior H\"{o}lder regularity}\setcounter{equation}{0}

Let $L_{2s}(\mathbb{R}^n)=\{u\in L^1_{loc}(\mathbb{R}^n)|\quad\int_{\mathbb{R}^n}\frac{|u(x)|}{1+|x|^{n+2s}}dx<\infty\}$. Then for any function $u\in L_{2s}(\mathbb{R}^n)$,
$$
\int_{\mathbb{R}^n\setminus B_R(x_0)}\frac{|u(x)|}{|x-x_0|^{n+2s}}dx \hbox{ is finite for all }R>0,x_0\in\mathbb{R}^n.
$$
\begin{thm}[Interior H\"{o}lder regularity]
Let $s\in(0,1)$.
Let $v\in H^{s}(\mathbb{R}^n)$ be a weak solution of
\begin{equation}\label{3.1}
L^{2s}_{K}v(x)=g\quad\hbox{ in } B_3,\\
\end{equation}
where the kernel $K$ satisfies assumptions (\ref{1.2}) and
 $$
\lim_{h\rightarrow0}\sup_{x,y\in D,|x-y|\leq\epsilon_0}|A(x+h,y+h)-A(x,y)|=0\hbox{ for any compact set }D\subset B_3,
$$
 for some $\epsilon_0>0$, and $g\in L^\infty(B_3)$. Then $v\in C^\alpha(B_1)$ for every $\alpha\in(0,\min\{2s,1\})$. Moreover,  there exists a constant $C=C(n,s,\lambda,\Lambda,\alpha)$ such that
$$
[v]_{C^\alpha(B_1)}\leq C(\|\nabla^s_{\mathbb{R}^n}v\|_{L^2(B_3)}+\|g\|_{L^\infty(B_3)}).
$$
\end{thm}
 We derive Theorem 3.1 from  the following results.
\begin{prop}(\cite{N2021})
Let $s\in(0,1)$ and let $K$  satisfy the assumptions in Theorem 3.1. Let $v\in W^{s,2}_{loc}(B_3)\cap L_{2s}(\mathbb{R}^n)$ be a weak solution of (\ref{3.1}) in $B_3$, where   $g\in L^\infty(B_2)$. Then  $v\in C^\alpha(B_1)$ for any $\alpha\in(0,\min\{2s,1\})$, and there exists a constant $C=C(n,s,\lambda,\Lambda,\alpha)$ such that
$$
[v]_{C^\alpha(B_1)}\leq C\Big(\|v\|_{L^2(B_2)}+\int_{\mathbb{R}^n\setminus{B_2}}\frac{|v(y)|}{|y|^{n+2s}}dy+\|g\|_{L^\infty(B_2)}\Big).
$$

\end{prop}
\begin{prop}(\cite{N2020})
Let $s\in(0,1)$. For any $v\in H^{s}(\mathbb{R}^n)$, we have
$$
\int_{\mathbb{R}^n\setminus{B_2}}\frac{v(y)^2}{|y|^{n+2s}}dy\leq C(\|\nabla^s_{\mathbb{R}^n}v\|^2_{L^2(B_3)}+\|v\|^2_{L^2(B_3)}),
$$
where $C$ is a constant depending on $n$ and $s$.
\end{prop}
\begin{proof}[Proof of Theorem 3.1]
Since $\tilde{v}=v-\overline{v}_{B_3}\in H^{s}(\mathbb{R}^n)$ is a weak solution of
$$
L^{2s}_{K}\tilde{v}=g\quad \hbox{ in }B_3.
$$
Obviously,  it is also a weak solution of this equation in $B_2$.
By Proposition 3.1, Proposition 3.2, H\"{o}lder's inequality and   fractional Poincar\'{e} inequality, for any $\alpha\in (0,\min\{2s,1\})$,
\begin{eqnarray}
[v]_{C^\alpha(B_1)}&=&[\tilde{v}]_{C^\alpha(B_1)}\nonumber\\
&\leq& C\Big(\|\tilde{v}\|_{L^2(B_2)}+\int_{\mathbb{R}^n\setminus{B_2}}\frac{|\tilde{v}(y)|}{|y|^{n+2s}}dy+\|g\|_{L^\infty(B_2)}\Big)\nonumber\\
&\leq&C\Big[\|\tilde{v}\|_{L^2(B_2)}+\Big(\int_{\mathbb{R}^n\setminus{B_2}}\frac{\tilde{v}(y)^2}{|y|^{n+2s}}dy\Big)^{\frac{1}{2}}+\|g\|_{L^\infty(B_2)}\Big]\nonumber\\
&\leq& C\Big(\|\nabla^s_{\mathbb{R}^n}\tilde{v}\|_{L^2(B_3)}+\|\tilde{v}\|_{L^2(B_3)}+\|g\|_{L^\infty(B_2)}\Big)\nonumber\\
&\leq& C\Big(\|\nabla^s_{\mathbb{R}^n}v\|_{L^2(B_3)}+\|g\|_{L^\infty(B_3)}\Big),\nonumber
\end{eqnarray}
where the constant $C$ depends on $n$, $s$, $\lambda$, $\Lambda$ and $\alpha$.
\end{proof}

\section{ Higher  differentiability and integrability }\setcounter{equation}{0}
 In this section, we mainly establish a higher  differentiability and integrability  result, Theorem 4.1, for the following equation:
\begin{equation}\label{4.+1}
 L^{2s}_{K}u(x)=L^{2s}_Dg(x)+L^{s}_bf(x),\quad \hbox{ in }\Omega,\\
\end{equation}
where $K$ and $b$ satisfy the assumptions (\ref{1.2})-(\ref{1.3.1}), and $ D$ satisfies
\begin{equation}\label{4.+2}
|D(x,y)|\leq\frac{\tilde{\beta}}{|x-y|^{n+2s}}\quad\hbox{ for any }x,y\in\mathbb{R}^n\hbox{ with }x\neq y,
\end{equation}
for some constant $\tilde{\beta}>0$. We say that $u\in {W}^{s,2}(\mathbb{R}^n)$ is a weak solution of (\ref{4.+1}) if
\begin{eqnarray}
&&\int_{\mathbb{R}^n}\int_{\mathbb{R}^n}K(x,y)[u(x)-u(y)][\phi(x)-\phi(y)]dydx\nonumber\\
&=&\int_{\mathbb{R}^n}\int_{\mathbb{R}^n}D(x,y)[g(x)-g(y)][\phi(x)-\phi(y)]dydx+\int_{\mathbb{R}^n}f(x)L^s_b\phi(x)dx, \quad\forall\phi\in \mathcal{H}^s_0(\Omega),\nonumber
\end{eqnarray}
where   $g\in W^{s,2}(\mathbb{R}^n)$ and  $f\in L^2(\mathbb{R}^n)$.

\begin{thm}
Let the kernels $K$, $D$ and $b$ satisfy assumptions (\ref{1.2})-(\ref{1.3.1}) and (\ref{4.+2}).
Let $u\in {W}^{s,2}(\mathbb{R}^n)$ be a weak solution of (\ref{4.+1}) with $g\in W^{s,2}(\mathbb{R}^n)$ and  $f\in L^2(\mathbb{R}^n)$. For any ball $B=B_R(x_0)\subset\Omega$, if  $g\in W^{s+\delta_0,2}(B)$ and  $f\in L^{2+\delta_0}(B)$ for some  small $\delta_0\in(0,2s)$, then there exist  $\gamma=\gamma(n,s,\lambda,\Lambda,\tilde{\beta},\beta,\delta_0)\in(0,1-s)$   and $\delta=\delta(n,s,\lambda,\Lambda,\tilde{\beta},\beta,\delta_0)>0$ such that
   $ u\in {W}^{s+\gamma,2+\delta}(\frac{1}{2}B)$ and the following estimate holds:
\begin{eqnarray}
R^\gamma\Big(\fint_{\frac{1}{2}B}\int_{\frac{1}{2}B}\frac{|u(x)-u(y)|^{2+\delta}}{|x-y|^{n+(s+\gamma)(2+\delta)}}dydx\Big)^{\frac{1}{2+\delta}}
&\leq& C\Big[\Big(\fint_{B}(\nabla^{s}_{B}u)^{2}dx\Big)^{\frac{1}{2}}+R^{\delta_0}\Big(\fint_{B}(\nabla^{s+\delta_0}_{B}g)^{2}dx\Big)^{\frac{1}{2}}\nonumber\\
& &+R^{-s}\Big(\fint_{B}g^{2}dx\Big)^{\frac{1}{2}}+\Big(\fint_{B}|f|^{2+\delta_0}dx\Big)^{\frac{1}{2+\delta_0}}\nonumber\\
& &+\Big(\int_{\mathbb{R}^n\setminus B}\frac{[u(y)-\bar{u}_{B}]^2}{|y|^{n+2s}}dy\Big)^{\frac{1}{2}}
+\Big(\int_{\mathbb{R}^n\setminus B}\frac{|g(y)-\bar{g}_{B}|^2}{|y|^{n+2s}}dy\Big)^{\frac{1}{2}}\nonumber\\
& &
+R^{\frac{s}{2}}\Big(\int_{\mathbb{R}^n\setminus B}\frac{|f(y)|^2}{|y|^{n+s}}dy\Big)^{\frac{1}{2}}
\Big],\nonumber
\end{eqnarray}
where $C=C(n,s,\lambda,\Lambda,\tilde{\beta},\beta,\delta_0)$.
\end{thm}
We first prove the following Caccioppoli-type estimate for solutions of equation (\ref{4.+1}), following a procedure similar to that in \cite{KMS2015} and incorporating some additional ideas.
\begin{lem}\label{lemma-2}
Let $u\in W^{s,2}(\mathbb{R}^n)$ be a weak solution of
(\ref{4.+1})
with $g\in W^{s,2}(\mathbb{R}^n)$ and $f\in L^2(\mathbb{R}^n)$. For some $\delta_0>0$, let  $q\in (\max\{\frac{2n}{n+2s}, \frac{2n}{n+\delta_0}\}, 2)$  and  $B=B_r(x_0)\subset \Omega$ for some $r>0$.  If $g\in W^{s+\delta_0,2}(B)$, then  there exists a constant $C=C(n,s,\lambda,\Lambda,\alpha,\beta,\delta_0,q)$ such that
 \begin{eqnarray}
&&\int_{\frac{1}{2}B}\int_{\frac{1}{2}B}\frac{[u(x)-u(y)]^2}{|x-y|^{n+2s}}dydx\nonumber\\
&\leq&C\Big[\frac{1}{r^{2s}\zeta^2}\int_{B}u^2(x)dx+\zeta^2r^{n+2s}\Big(\int_{\mathbb{R}^n\setminus B} \frac{|u(x)|}{|x-x_0|^{n+2s}}dx\Big)^2+r^{n-2s+2\sigma-\frac{2n}{q}}\Big(\int_{B}\int_{B}\frac{|g(x)-g(y)|^q}{|x-y|^{n+q\sigma}}dydx\Big)^{\frac{2}{q}}\nonumber\\
& &+r^{n+2s}\Big(\int_{\mathbb{R}^n\setminus B}\frac{|g(y)-\bar{g}_B|}{|y-x_0|^{n+2s}}dy\Big)^2
+\int_Bf^2(x)dx+r^{n+s}\int_{\mathbb{R}^n\setminus B}\frac{f^2(x)}{|x-x_0|^{n+s}}dx\Big],\nonumber
\end{eqnarray}
for any $\zeta\in(0,1)$, where $\sigma=s+\frac{\delta_0}{2}$ and $\bar{g}_B=\fint_Bg(x)dx=\frac{1}{|B|}\int_Bg(x)dx$.
\end{lem}
\begin{proof} Let
$\psi\in C^\infty_0(\frac{3}{4}B)$ be a cut-off function
satisfying
\begin{equation}
0\leq\psi\leq1, \quad \psi\equiv1\hbox{ in }\frac{1}{2}B, \quad\psi\equiv0\hbox{ in }\mathbb{R}^n\setminus\frac{3}{4}B\quad\hbox{ and }\quad|D\psi|\leq\frac{c(n)}{r}.\nonumber\\
\end{equation}
Take $\phi=\psi^2u\in \mathcal{H}^{s}_0(\Omega)$ in (\ref{4.+2}), we have
\begin{eqnarray}\label{4.+3}
&&\int_{\mathbb{R}^n}\int_{\mathbb{R}^n}K(x,y)[u(x)-u(y)][\psi^2(x)u(x)-\psi^2(y)u(y)]dydx\nonumber\\
&=&\int_{\mathbb{R}^n}\int_{\mathbb{R}^n}D(x,y)[g(x)-g(y)][\psi^2(x)u(x)-\psi^2(y)u(y)]dydx+\int_{\mathbb{R}^n}f(x)L^s_b(\psi^2u)(x)dx.
\end{eqnarray}
We rewrite the left-hand side of (\ref{4.+3}) as
\begin{eqnarray}\label{4.+301}
&&\int_{\mathbb{R}^n}\int_{\mathbb{R}^n}K(x,y)\varphi(u(x)-u(y))[\psi^2(x)u(x)-\psi^2(y)u(y)]dydx\nonumber\\
&=&\int_{B}\int_{B}K(x,y)[u(x)-u(y)][\psi^2(x)u(x)-\psi^2(y)u(y)]dydx\nonumber\\
& &+\int_{B}\int_{\mathbb{R}^n\setminus B}K(x,y)[u(x)-u(y)]\psi^2(x)u(x)dydx\nonumber\\
& &-\int_{\mathbb{R}^n\setminus B}\int_{B}K(x,y)[u(x)-u(y)]\psi^2(y)u(y)dydx
\end{eqnarray}
If $\psi(x)<\psi(y)$, by assumption (\ref{1.2}), we get
\begin{eqnarray}
&&[u(x)-u(y)][\psi^2(x)u(x)-\psi^2(y)u(y)]\nonumber\\
&=&[u(x)-u(y)]^2\psi^2(y)+[u(x)-u(y)][\psi^2(x)-\psi^2(y)]u(x)\nonumber\\
&\geq&[u(x)-u(y)]^2\psi^2(y)-2|\psi(x)-\psi(y)||u(x)||u(x)-u(y)|\psi(x)\nonumber\\
&\geq&\frac{1}{2}[u(x)-u(y)]^2\psi^2(y)-2[\psi(x)-\psi(y)]^2u^2(x).\nonumber
\end{eqnarray}
Similarly, if  $\psi(x)\geq\psi(y)$, we get
\begin{eqnarray}
[u(x)-u(y)][\psi^2(x)u(x)-\psi^2(y)u(y)]
\geq\frac{1}{2}[u(x)-u(y)]^2\psi^2(x)-2[\psi(x)-\psi(y)]^2u^2(y).\nonumber
\end{eqnarray}
Therefore, by assumption (\ref{1.3}), we have
\begin{eqnarray}\label{4.+302}
&&\int_{B}\int_{B}K(x,y)[u(x)-u(y)][\psi^2(x)u(x)-\psi^2(y)u(y)]dydx\nonumber\\
&\geq&\int_{(B\times B)\cap\{\psi(x)>\psi(y)\}}\frac{\lambda}{2}\frac{[u(x)-u(y)]^2\psi^2(x)}{|x-y|^{n+2s}}
-2\Lambda\frac{[\psi(x)-\psi(y)]^2u^2(y)}{|x-y|^{n+2s}}dydx\nonumber\\
& &+\int_{(B\times B)\cap\{\psi(x)\leq\psi(y)\}}\frac{\lambda}{2}\frac{[u(x)-u(y)]^2\psi^2(y)}{|x-y|^{n+2s}}
-2\Lambda\frac{[\psi(x)-\psi(y)]^2u^2(x)}{|x-y|^{n+2s}}dydx\nonumber\\
&\geq&\frac{\lambda}{2}\int_{B}\int_{B}\frac{[u(x)-u(y)]^2\psi^2(x)}{|x-y|^{n+2s}}dydx
-4\Lambda\int_{B}\int_{B}\frac{[\psi(x)-\psi(y)]^2u^2(x)}{|x-y|^{n+2s}}dydx\nonumber\\
&\geq&\frac{\lambda}{2}\int_{B}\int_{B}\frac{[u(x)-u(y)]^2\psi^2(x)}{|x-y|^{n+2s}}dydx
-c(n,s,\Lambda)\frac{1}{r^{2s}}\int_{B} u^2(x)dx.
\end{eqnarray}
Using H\"{o}lder's inequality and Young's inequality,  we get
\begin{eqnarray}\label{4.+303}
&&\int_{B}\int_{\mathbb{R}^n\setminus B}K(x,y)[u(x)-u(y)]\psi^2(x)u(x)dydx\nonumber\\
&\geq&-\Lambda\int_{\frac{3}{4}B}|u(x)|\int_{\mathbb{R}^n\setminus B}\frac{|u(y)|}{|x-y|^{n+2s}}dydx\nonumber\\
&\geq&-c(n,s,\Lambda)\int_{\frac{3}{4}B}|u(x)|dx\int_{\mathbb{R}^n\setminus B}\frac{|u(y)|}{|y-x_0|^{n+2s}}dy\nonumber\\
&\geq&-c(n,s,\Lambda)r^{\frac{n}{2}}\Big(\int_{B}u^2(x)dx\Big)^{\frac{1}{2}}\int_{\mathbb{R}^n\setminus B}\frac{|u(y)|}{|y-x_0|^{n+2s}}dy\nonumber\\
&\geq&-c(n,s,\Lambda)\Big[\frac{1}{\zeta^2r^{2s}}\int_{B}u^2(x)dx+\zeta^2r^{n+2s}\Big(\int_{\mathbb{R}^n\setminus B}\frac{|u(y)|}{|y-x_0|^{n+2s}}dy\Big)^2\Big],
\end{eqnarray}
for any $\zeta\in(0,1)$, where we used the fact that
\begin{equation}
|y-x_0|\leq |y-x|+|x-x_0|=|y-x|\Big(1+\frac{|x-x_0|}{|x-y|}\Big)\leq 4|x-y|,
\nonumber
\end{equation}
for  any $y\in\mathbb{R}^n\setminus B$ and $x\in \frac{3}{4}B$. The estimate of the last term on the right-hand side of (\ref{4.+301}) is similar to that of (\ref{4.+303}). Therefore, we get
\begin{eqnarray}\label{4.+5}
&&\int_{\mathbb{R}^n}\int_{\mathbb{R}^n}K(x,y)\varphi(u(x)-u(y))[\psi^2(x)u(x)-\psi^2(y)u(y)]dydx\nonumber\\
&\geq&\frac{\lambda}{2}\int_{B}\int_{B}\frac{\psi^2(x)[u(x)-u(y)]^2}{|x-y|^{n+2s}}dydx\nonumber\\
& &-c(n,s,\Lambda)\Big[\frac{1}{\zeta^2 r^{2s}}\int_{B} u^2(x)dx
+\zeta^2r^{n+2s}\Big(\int_{\mathbb{R}^n\setminus B} \frac{|u(x)|}{|x-x_0|^{n+2s}}dx\Big)^2\Big].
\end{eqnarray}

Now, we estimate the first term on the right-hand side of (\ref{4.+3}). We rewrite it as
\begin{eqnarray}\label{4.+6}
&&\int_{\mathbb{R}^n}\int_{\mathbb{R}^n}D(x,y)[g(x)-g(y)][\psi^2(x)u(x)-\psi^2(y)u(y)]dydx\nonumber\\
&=&\int_{\mathbb{R}^n}\int_{\mathbb{R}^n}D(x,y)[g(x)-g(y)][\psi(x)u(x)-\psi(y)u(y)]\psi(x)dydx\nonumber\\
& &+\int_{\mathbb{R}^n}\int_{\mathbb{R}^n}D(x,y)[g(x)-g(y)][\psi(x)-\psi(y)]\psi(y)u(y)dydx\nonumber\\
&=&I+II.
\end{eqnarray}
and
\begin{eqnarray}
I&=&\int_{B}\int_{\mathbb{R}^n}D(x,y)[g(x)-g(y)][\psi(x)u(x)-\psi(y)u(y)]\psi(x)dydx\nonumber\\
&\leq&\tilde{\beta}\int_{B}\int_{B}\frac{|g(x)-g(y)||\psi(x)u(x)-\psi(y)u(y)|}{|x-y|^{n+2s}}dydx
+\tilde{\beta}\int_{B}\int_{\mathbb{R}^n\setminus B}\frac{|g(x)-g(y)||\psi(x)u(x)|}{|x-y|^{n+2s}}dydx\nonumber\\
&=&I_1+I_2.\nonumber
\end{eqnarray}
Notice that $\frac{2n}{n+\delta_0}<q<2$.  Take   $\sigma=s+\frac{\delta_0}{2}$, $\hat{s}=\frac{n}{s}\Big(\frac{2-q}{2q}\Big)$. Then it follows that
\begin{equation}
\frac{1}{q}+\frac{1-\hat{s}}{2}+\frac{\hat{s}}{2^*}=1.\nonumber\\
\end{equation}
By H\"{o}lder's inequality, we have
\begin{eqnarray}
I_1&=&\tilde{\beta}\int_{B}\int_{B}\frac{|g(x)-g(y)||\psi(x)u(x)-\psi(y)u(y)|}{|x-y|^{n+2s}}dydx\nonumber\\
&=&\tilde{\beta}\int_{B}\int_{B}\frac{|g(x)-g(y)|}{|x-y|^{\frac{n+q\sigma}{q}}}
\frac{|\psi(x)u(x)-\psi(y)u(y)|^{1-\hat{s}}}{|x-y|^{\frac{(n+2s)(1-\hat{s})}{2}}}
\frac{|\psi(x)u(x)-\psi(y)u(y)|^{\hat{s}}}{|x-y|^{s\hat{s}+s-\sigma+n\frac{\hat{s}}{2^*}}}dydx\nonumber\\
&\leq&\tilde{\beta}\Big(\int_{B}\int_{B}\frac{|g(x)-g(y)|^q}{|x-y|^{n+q\sigma}}dydx\Big)^{\frac{1}{q}}
\cdot\Big(\int_{B}\int_{B}\frac{|\psi(x)u(x)-\psi(y)u(y)|^2}{|x-y|^{n+2s}}dydx\Big)^{\frac{1-\hat{s}}{2}}\nonumber\\
& &\cdot
\Big(\int_{B}\int_{B}\frac{|\psi(x)u(x)-\psi(y)u(y)|^{2^*}}{|x-y|^{n+(s\hat{s}+s-\sigma)\frac{2^*}{\hat{s}}}}dydx\Big)^{\frac{\hat{s}}{2^*}}.\nonumber
\end{eqnarray}
Notice that by a direct calculation and fractional Sobolev inequality,
\begin{eqnarray}
\int_{B}\int_{B}\frac{|\psi(x)u(x)-\psi(y)u(y)|^{2^*}}{|x-y|^{n+(s\hat{s}+s-\sigma)\frac{2^*}{\hat{s}}}}dydx
&\leq&2^{2^*}\int_{B}\int_{B}\frac{|\psi(x)u(x)|^{2^*}}{|x-y|^{n+(s\hat{s}+s-\sigma)\frac{2^*}{\hat{s}}}}dydx\nonumber\\
&\leq&2^{2^*}\int_{B}|\psi(x)u(x)|^{2^*}dx\int_{2B}\frac{1}{|z|^{n+(s\hat{s}+s-\sigma)\frac{2^*}{\hat{s}}}}dz\nonumber\\
&\leq&c(n,s,q,\delta_0)r^{(\sigma-s-s\hat{s})\frac{2^*}{\hat{s}}}\int_{B}|\psi(x)u(x)|^{2^*}dx\nonumber\\
&\leq&C(n,s,q,\delta_0)r^{(\sigma-s-s\hat{s})\frac{2^*}{\hat{s}}}\Big(\int_{B}\int_{B}\frac{|\psi(x)u(x)-\psi(y)u(y)|^2}{|x-y|^{n+2s}}dydx\Big)^{\frac{2^*}{2}}.\nonumber
\end{eqnarray}
Therefore, we get
\begin{eqnarray}\label{4.+7}
I_1\leq C(n,s,\tilde{\beta},q,\delta_0)r^{\frac{n}{2}-\frac{n}{q}+\sigma-s}\Big(\int_{B}\int_{B}\frac{|g(x)-g(y)|^q}{|x-y|^{n+q\sigma}}dydx\Big)^{\frac{1}{q}}
\Big(\int_{B}\int_{B}\frac{|\psi(x)u(x)-\psi(y)u(y)|^2}{|x-y|^{n+2s}}dydx\Big)^{\frac{1}{2}}.
\end{eqnarray}
 We rewrite $I_2$ as
\begin{eqnarray}\label{4.+8}
I_2&=&\tilde{\beta}\int_{B}|\psi(x)u(x)|\int_{\mathbb{R}^n\setminus B}\frac{|g(x)-g(y)|}{|x-y|^{n+2s}}dydx\nonumber\\
&\leq&\tilde{\beta}\int_{B}|\psi(x)u(x)|\int_{\mathbb{R}^n\setminus B}\frac{|g(x)-\bar{g}_B|}{|x-y|^{n+2s}}dydx+\tilde{\beta}\int_{B}|\psi(x)u(x)|\int_{\mathbb{R}^n\setminus B}\frac{|g(y)-\bar{g}_B|}{|x-y|^{n+2s}}dydx.
\end{eqnarray}
Notice the fact that $2<\frac{nq}{n-q\sigma}$. By H\"{o}lder's inequality, fractional Poincar\'{e} inequality and fractional Sobolev-Poincar\'{e} inequality, we have
\begin{eqnarray}\label{4.+9}
&&\tilde{\beta}\int_{B}|\psi(x)u(x)|\int_{\mathbb{R}^n\setminus B}\frac{|g(x)-\bar{g}_B|}{|x-y|^{n+2s}}dydx\nonumber\\
&\leq&c(n,s,\tilde{\beta})r^{-2s}\int_{B}|\psi(x)u(x)||g(x)-\bar{g}_B|dx\nonumber\\
&\leq&c(n,s,\tilde{\beta})r^{-2s}\Big(\int_{B}|\psi(x)u(x)|^2dx\Big)^{\frac{1}{2}}\Big(\int_{B}|g(x)-\bar{g}_B|^2dx\Big)^{\frac{1}{2}}\nonumber\\
&\leq&c(n,s,\tilde{\beta},q,\delta_0)r^{\frac{n}{2}-\frac{n}{q}+\sigma-s}\Big(\int_{B}\int_B\frac{[\psi(x)u(x)-\psi(y)u(y)]^2}{|x-y|^{n+2s}}dx\Big)^{\frac{1}{2}}
\Big(\int_{B}\int_{B}\frac{|g(x)-g(y)|^q}{|x-y|^{n+q\sigma}}dydx\Big)^{\frac{1}{q}}.
\end{eqnarray}
By fractional Poincar\'{e} inequality and the fact that
$$
|x_0-y|\leq |x_0-x|+|x-y|\leq |x-y|\Big(1+\frac{|x_0-x|}{|x-y|}\Big)\leq 4|x-y|\quad\hbox{ for any }x\in \frac{3}{4}B\hbox{ and }y\in \mathbb{R}^n\setminus B,
$$
 we get
\begin{eqnarray}\label{4.+10}
&&\tilde{\beta}\int_{B}\psi(x)|u(x)|\int_{\mathbb{R}^n\setminus B}\frac{|g(y)-\bar{g}_B|}{|x-y|^{n+2s}}dydx\nonumber\\
&\leq&c(n,s,\tilde{\beta})\int_{B}\psi(x)|u(x)|dx\int_{\mathbb{R}^n\setminus B}\frac{|g(y)-\bar{g}_B|}{|y-x_0|^{n+2s}}dy\nonumber\\
&\leq&c(n,s,\tilde{\beta})r^{\frac{n}{2}}\Big(\int_{B}\psi^2(x)u^2(x)dx\Big)^{\frac{1}{2}}\int_{\mathbb{R}^n\setminus B}\frac{|g(y)-\bar{g}_B|}{|y-x_0|^{n+2s}}dy\nonumber\\
&\leq&c(n,s,\tilde{\beta})r^{\frac{n}{2}+s}\Big(\int_{B}\int_{B}\frac{[\psi (x)u(x)-\psi(y)u(y)]^2}{|x-y|^{n+2s}}dydx\Big)^{\frac{1}{2}}\int_{\mathbb{R}^n\setminus B}\frac{|g(y)-\bar{g}_B|}{|y-x_0|^{n+2s}}dy.
\end{eqnarray}
Combining (\ref{4.+7})-(\ref{4.+10}), we obtain
\begin{eqnarray}\label{4.+11}
I&\leq&\leq C(n,s,\tilde{\beta},q,\delta_0)\Big(\int_{B}\int_{B}\frac{|\psi(x)u(x)-\psi(y)u(y)|^2}{|x-y|^{n+2s}}dydx\Big)^{\frac{1}{2}} \Big[r^{\frac{n}{2}-\frac{n}{q}+\sigma-s}\Big(\int_{B}\int_{B}\frac{|g(x)-g(y)|^q}{|x-y|^{n+q\sigma}}dx\Big)^{\frac{1}{q}}\nonumber\\
& &+r^{\frac{n}{2}+s}\int_{\mathbb{R}^n\setminus B}\frac{|g(y)-\bar{g}_B|}{|y-x_0|^{n+2s}}dy\Big]\nonumber\\
&\leq&\epsilon\int_{B}\int_{B}\frac{|\psi(x)u(x)-\psi(y)u(y)|^2}{|x-y|^{n+2s}}dydx+\frac{ C(n,s,\tilde{\beta},q,\delta_0)}{\epsilon}\Big[r^{n-\frac{2n}{q}+2\sigma-2s}\Big(\int_{B}\int_{B}\frac{|g(x)-g(y)|^q}{|x-y|^{n+q\sigma}}dydx\Big)^{\frac{2}{q}}\nonumber\\
& &+r^{n+2s}\Big(\int_{\mathbb{R}^n\setminus B}\frac{|g(y)-\bar{g}_B|}{|y-x_0|^{n+2s}}dy\Big)^{2}\Big]\nonumber\\
&\leq&2\epsilon\int_{B}\int_{B}\frac{\psi^2(x)[u(x)-u(y)]^2}{|x-y|^{n+2s}}dydx+\epsilon\frac{c(n,s)}{r^{2s}}\int_B u^2(x)dx\nonumber\\
& &+\frac{ C}{\epsilon}\Big[r^{n-\frac{2n}{q}+2\sigma-2s}\Big(\int_{B}\int_{B}\frac{|g(x)-g(y)|^q}{|x-y|^{n+q\sigma}}dydx\Big)^{\frac{2}{q}}
+r^{n+2s}\Big(\int_{\mathbb{R}^n\setminus B}\frac{|g(y)-\bar{g}_B|}{|y-x_0|^{n+2s}}dy\Big)^{2}\Big],
\end{eqnarray}
for any $\epsilon>0$, where $C=C(n,s,\tilde{\beta},q,\delta_0)$.

We rewrite II as
\begin{eqnarray}
II&=&\int_{\mathbb{R}^n}\int_{B}A(x,y)[g(x)-g(y)][\psi(x)-\psi(y)]\psi(y)u(y)dydx\nonumber\\
&\leq&\tilde{\beta}\int_{B}\int_{B}\frac{|g(x)-g(y)||\psi(x)-\psi(y)|\psi(y)|u(y)|}{|x-y|^{n+2s}}dydx
+\tilde{\beta}\int_{\mathbb{R}^n\setminus B}\int_{B}\frac{|g(x)-g(y)|\psi(y)|u(y)|}{|x-y|^{n+2s}}dydx\nonumber\\
&=&II_1+II_2.\nonumber
\end{eqnarray}
Notice the fact that  $\frac{q}{q-1}<2^*$. By H\"{o}lder's inequality and  fractional Sobolev inequality, we get
\begin{eqnarray}
II_1&=&\tilde{\beta}\int_{B}\int_{B}\frac{|g(x)-g(y)||\psi(x)-\psi(y)|\psi(y)|u(y)|}{|x-y|^{n+2s}}dydx\nonumber\\
&\leq&\frac{c(n,\tilde{\beta})}{r}\int_{B}\int_{B}\frac{|g(x)-g(y)|\psi(y)|u(y)|}{|x-y|^{n+2s-1}}dydx\nonumber\\
&=&\frac{c(n,\tilde{\beta})}{r}\int_{B}\int_{B}\frac{|g(x)-g(y)|}{|x-y|^{\frac{n+q\sigma}{q}}}\frac{\psi(y)|u(y)|}{|x-y|^{n+2s-1-\frac{n+q\sigma}{q}}}dydx\nonumber\\
&\leq&\frac{c(n,\tilde{\beta})}{r}\Big(\int_{B}\int_{B}\frac{|g(x)-g(y)|^q}{|x-y|^{n+q\sigma}}dydx\Big)^{\frac{1}{q}}
\Big(\int_{B}\int_{B}\frac{|\psi(y)u(y)|^{\frac{q}{q-1}}}{|x-y|^{(n+2s-1-\frac{n+q\sigma}{q}){\frac{q}{q-1}}}}dydx\Big)^{\frac{q-1}{q}}\nonumber\\
&\leq&c(n,s,\tilde{\beta},q,\delta_0)r^{\sigma-2s}\Big(\int_{B}\int_{B}\frac{|g(x)-g(y)|^q}{|x-y|^{n+q\sigma}}dydx\Big)^{\frac{1}{q}}
\Big(\int_{B}|\psi(y)u(y)|^{\frac{q}{q-1}}dy\Big)^{\frac{q-1}{q}}\nonumber\\
&\leq&C(n,s,\tilde{\beta},q,\delta_0)r^{\frac{n}{2}-\frac{n}{q}+\sigma-s}\Big(\int_{B}\int_{B}\frac{|g(x)-g(y)|^q}{|x-y|^{n+q\sigma}}dydx\Big)^{\frac{1}{q}}
\Big(\int_{B}\int_B\frac{[\psi(x)u(x)-\psi(y)u(y)]^2}{|x-y|^{n+2s}}dx\Big)^{\frac{1}{2}}.\nonumber
\end{eqnarray}
By a similar process to the estimate for $I_2$, we get
\begin{eqnarray}
II_2&\leq&c(n,s,\tilde{\beta},q,\delta_0)\Big(\int_{B}\int_B\frac{[\psi(x)u(x)-\psi(y)u(y)]^2}{|x-y|^{n+2s}}dx\Big)^{\frac{1}{2}}\Big[r^{\frac{n}{2}-\frac{n}{q}+\sigma-s}
\Big(\int_{B}\int_{B}\frac{|g(x)-g(y)|^q}{|x-y|^{n+q\sigma}}dydx\Big)^{\frac{1}{q}}\nonumber\\
& &+r^{\frac{n}{2}+s}\int_{\mathbb{R}^n\setminus B}\frac{|g(y)-\bar{g}_B|}{|y-x_0|^{n+2s}}dy\Big].\nonumber
\end{eqnarray}
Therefore,we obtain
\begin{eqnarray}\label{4.+12}
II&\leq&C(n,s,\tilde{\beta},q,\delta_0)\Big(\int_{B}\int_B\frac{[\psi(x)u(x)-\psi(y)u(y)]^2}{|x-y|^{n+2s}}dx\Big)^{\frac{1}{2}}\Big[r^{\frac{n}{2}-\frac{n}{q}+\sigma-s}
\Big(\int_{B}\int_{B}\frac{|g(x)-g(y)|^q}{|x-y|^{n+q\sigma}}dydx\Big)^{\frac{1}{q}}\nonumber\\
& &+r^{\frac{n}{2}+s}\int_{\mathbb{R}^n\setminus B}\frac{|g(y)-\bar{g}_B|}{|y-x_0|^{n+2s}}dy\Big]\nonumber\\
&\leq&\epsilon\int_{B}\int_{B}\frac{|\psi(x)u(x)-\psi(y)u(y)|^2}{|x-y|^{n+2s}}dydx+\frac{ C(n,s,\tilde{\beta},q,\delta_0)}{\epsilon}\Big[r^{n-\frac{2n}{q}+2\sigma-2s}\Big(\int_{B}\int_{B}\frac{|g(x)-g(y)|^q}{|x-y|^{n+q\sigma}}dydx\Big)^{\frac{2}{q}}\nonumber\\
& &+r^{n+2s}\Big(\int_{\mathbb{R}^n\setminus B}\frac{|g(y)-\bar{g}_B|}{|y-x_0|^{n+2s}}dy\Big)^{2}\Big]\nonumber\\
&\leq&2\epsilon\int_{B}\int_{B}\frac{\psi^2(x)[u(x)-u(y)]^2}{|x-y|^{n+2s}}dydx+\epsilon\frac{c(n,s)}{r^{2s}}\int_B u^2(x)dx\nonumber\\
& &+\frac{ C}{\epsilon}\Big[r^{n-\frac{2n}{q}+2\sigma-2s}\Big(\int_{B}\int_{B}\frac{|g(x)-g(y)|^q}{|x-y|^{n+q\sigma}}dydx\Big)^{\frac{2}{q}}
+r^{n+2s}\Big(\int_{\mathbb{R}^n\setminus B}\frac{|g(y)-\bar{g}_B|}{|y-x_0|^{n+2s}}dy\Big)^{2}\Big],
\end{eqnarray}
for any $\epsilon>0$, where $C=C(n,s,\tilde{\beta},q,\delta_0)$.

For the second term on the right-hand side of (\ref{4.+3}), by using  H\"{o}lder's inequality, Cauchy inequality and assumption (\ref{1.3.1}), we have
\begin{eqnarray}\label{f}
&&\int_{\mathbb{R}^n}f(x)L^s_b(\psi^2u)(x)dx\nonumber\\
&=&\int_{\mathbb{R}^n}f(x)\int_{\mathbb{R}^n}b(x,y)[\psi^2(x)u(x)-\psi^2(y)u(y)]dydx\nonumber\\
&=&\int_{\mathbb{R}^n}f(x)\int_{\mathbb{R}^n}b(x,y)[\psi(x)-\psi(y)]\psi(y)u(y)dydx
+\int_{\mathbb{R}^n}f(x)\psi(x)\int_{\mathbb{R}^n}b(x,y)[\psi(x)u(x)-\psi(y)u(y)]dydx
\nonumber\\
&\leq&\int_{\mathbb{R}^n}f(x)\int_{B}b(x,y)[\psi(x)-\psi(y)]\psi(y)u(y)dydx+\|f\psi\|_{L^2(\mathbb{R}^n)}\|L^s_b(\psi u)\|_{L^2(\mathbb{R}^n)}
\nonumber\\
&\leq&\int_{\mathbb{R}^n}f(x)\int_{B}b(x,y)[\psi(x)-\psi(y)]\psi(y)u(y)dydx+c(n,s,\beta)\epsilon_1[\psi u]^2_{H^s(\mathbb{R}^n)}+\frac{1}{\epsilon_1}\|f\|^2_{L^2(B)}.\nonumber
\end{eqnarray}
Here $\epsilon_1\in(0,1)$. By a direct calculation,   we derive that
\begin{eqnarray}
[\psi u]^2_{H^s(\mathbb{R}^n)}
&=&\int_{\mathbb{R}^n}\int_{\mathbb{R}^n}\frac{[\psi(x)u(x)-\psi(y)u(y)]^2}{|x-y|^{n+2s}}dydx\nonumber\\
&\leq&2\Big[\int_{B}\int_{B}\psi^2(x)\frac{[u(x)-u(y)]^2}{|x-y|^{n+2s}}dydx
+\int_{B}\int_{B}\frac{[\psi(x)-\psi(y)]^2}{|x-y|^{n+2s}}u^2(y)dydx\nonumber\\
& &+\int_{B}\int_{\mathbb{R}^n\setminus B}\frac{\psi^2(x) u^2(x)}{|x-y|^{n+2s}}dydx\Big]\nonumber\\
&\leq&2\int_{B}\int_{B}\psi^2(x)\frac{[u(x)-u(y)]^2}{|x-y|^{n+2s}}dydx
+c(n,s)\frac{1}{r^{2s}}\int_{B}u^2(x)dx.\nonumber
\end{eqnarray}
We divide the following integral into two parts,
\begin{eqnarray}
&&\int_{\mathbb{R}^n}f(x)\int_{B}b(x,y)[\psi(x)-\psi(y)]\psi(y)u(y)dydx\nonumber\\
&=&\int_{\mathbb{R}^n}f(x)\int_{B_{\frac{3r}{4}}(x_0)}b(x,y)[\psi(x)-\psi(y)]\psi(y)u(y)dydx\nonumber\\
&=&\int_{\mathbb{R}^n}f(x)\int_{B_{\frac{3r}{4}}(x_0)\cap\{|x-y|<\frac{r}{8}\}}b(x,y)[\psi(x)-\psi(y)]\psi(y)u(y)dydx\nonumber\\
& &+\int_{\mathbb{R}^n}f(x)\int_{B_{\frac{3r}{4}}(x_0)\cap\{|x-y|\geq\frac{r}{8}\}}b(x,y)[\psi(x)-\psi(y)]\psi(y)u(y)dydx\nonumber
.\nonumber
\end{eqnarray}
For the first part,
by H\"{o}lder's inequality and Cauchy inequality, we can deduce that
\begin{eqnarray}
&&\int_{\mathbb{R}^n}f(x)\int_{B_{\frac{3r}{4}}(x_0)\cap\{|x-y|<\frac{r}{8}\}}b(x,y)[\psi(x)-\psi(y)]\psi(y)u(y)dydx\nonumber\\
&\leq&\beta\int_{\mathbb{R}^n}|f(x)|\int_{B_{\frac{3r}{4}}(x_0)\cap\{|x-y|<\frac{r}{8}\}}\frac{|\psi(x)-\psi(y)||u(y)|}{|x-y|^{n+s}}dydx\nonumber\\
&\leq&\frac{C(n,\beta)}{r}\int_{\mathbb{R}^n}|f(x)|\int_{B_{\frac{3r}{4}}(x_0)\cap\{|x-y|<\frac{r}{8}\}}\frac{|u(y)|}{|x-y|^{n+s-1}}dydx
\nonumber\\
&\leq&\frac{C(n,\beta)}{r}\int_{B_{\frac{r}{8}}(0)}\frac{1}{|z|^{n+s-1}}dz\int_{B_{\frac{7r}{8}}(x_0)}|f(x)||u(x-z)|dx\nonumber\\
&\leq&c(n,s,\beta)r^{-s}\Big(\int_{B}f^2(x)dx\Big)^{\frac{1}{2}}\Big(\int_{B}u^2(x)dx\Big)^{\frac{1}{2}}\nonumber\\
&\leq&\tilde{c}(n,s,\beta)\Big[\int_{B}f^2(x)dx+\frac{1}{r^{2s}}\int_{B}u^2(x)dx\Big],\nonumber
\end{eqnarray}
For the second part, by H\"{o}lder's inequality, we get
\begin{eqnarray}&&\int_{\mathbb{R}^n}f(x)\int_{B_{\frac{3r}{4}}(x_0)\cap\{|x-y|\geq\frac{r}{8}\}}b(x,y)[\psi(x)-\psi(y)]\psi(y)u(y)dydx\nonumber\\
&\leq&\beta\int_{\mathbb{R}^n}|f(x)|\int_{B_{\frac{3r}{4}}(x_0)\cap\{|x-y|\geq\frac{r}{8}\}}\frac{|\psi(x)-\psi(y)||u(y)|}{|x-y|^{n+s}}dydx
\nonumber\\
&=&\beta\int_{B}|f(x)|\int_{B_{\frac{3r}{4}}(x_0)\cap\{|x-y|\geq\frac{r}{8}\}}\frac{|\psi(x)-\psi(y)||u(y)|}{|x-y|^{n+s}}dydx\nonumber\\
& &
+\beta\int_{\mathbb{R}^n\setminus B}|f(x)|\int_{B_{\frac{3r}{4}}(x_0)\cap\{|x-y|\geq\frac{r}{8}\}}\frac{|\psi(x)-\psi(y)||u(y)|}{|x-y|^{n+s}}dydx\nonumber\\
&\leq&c(n,s,\beta)\Big[\frac{1}{r^{n+s}}\int_{B}|f(x)|dx\int_{B}|u(y)|dy
+\int_{\frac{3}{4}B}|u(y)|\int_{\mathbb{R}^n\setminus B}\frac{|f(x)|}{|x-x_0|^{n+s}}dxdy\Big]\nonumber\\
&\leq&c(n,s,\beta)\Big[\frac{1}{r^{s}}\Big(\int_{B}f^2(x)dx\Big)^{\frac{1}{2}}\Big(\int_{B}u^2(y)dy\Big)^{\frac{1}{2}}
+r^{{\frac{n}{2}-\frac{s}{2}}}\Big(\int_{\frac{3}{4}B}u^2(y)dy\Big)^{\frac{1}{2}}\Big(\int_{\mathbb{R}^n\setminus B}\frac{f^2(x)}{|x-x_0|^{n+s}}dx\Big)^{\frac{1}{2}}\Big]\nonumber\\
&\leq&2c(n,s,\beta)\Big[\int_{B}f^2(x)dx+\frac{1}{r^{2s}}\int_{B}u^2(y)dy
+r^{n+s}\int_{\mathbb{R}^n\setminus B}\frac{f^2(x)}{|x-x_0|^{n+s}}dx\Big].\nonumber
\end{eqnarray}
Therefore, we obtain
\begin{eqnarray}\label{4.+13}
\int_{\mathbb{R}^n}f(x)L^s_b(\psi^2u)(x)dx
&\leq&2c(n,s,\beta)\epsilon_1\int_{B}\int_{B}\psi^2(x)\frac{[u(x)-u(y)]^2}{|x-y|^{n+2s}}dydx
+\frac{1}{\epsilon_1}\int_{B}f^2(x)dx\nonumber\\
& &
+c_1(n,s,\beta)\Big[\frac{1}{r^{2s}}\int_{B}u^2(x)dy
+r^{n+s}\int_{\mathbb{R}^n\setminus B}\frac{f^2(x)}{|x-x_0|^{n+s}}dx\Big].
\end{eqnarray}
Combining (\ref{4.+3}), (\ref{4.+5}), (\ref{4.+6}), (\ref{4.+11}), (\ref{4.+12}) and  (\ref{4.+13}) with $\epsilon=\frac{\lambda}{32}$ and $\epsilon_1=\frac{\lambda}{16c(n,s,\beta)}$, we obtain
\begin{eqnarray}
&&\int_{\frac{1}{2}B}\int_{\frac{1}{2}B}\frac{[u(x)-u(y)]^2}{|x-y|^{n+2s}}dydx\nonumber\\
&\leq&\int_{B}\int_{B}\psi^2(x)\frac{[u(x)-u(y)]^2}{|x-y|^{n+2s}}dydx\nonumber\\
&\leq&C\Big[\frac{1}{r^{2s}\zeta^2}\int_{B}u^2(x)dx+\zeta^2r^{n+2s}\Big(\int_{\mathbb{R}^n\setminus B} \frac{|u(x)|}{|x-x_0|^{n+2s}}dx\Big)^2+r^{n-2s+2\sigma-\frac{2n}{q}}\Big(\int_{B}\int_{B}\frac{|g(x)-g(y)|^q}{|x-y|^{n+q\sigma}}dydx\Big)^{\frac{2}{q}}\nonumber\\
& &+r^{n+2s}\Big(\int_{\mathbb{R}^n\setminus B}\frac{|g(y)-\bar{g}_B|}{|x_0-y|^{n+2s}}dy\Big)^{2}
+\int_Bf^2(x)dx+r^{n+s}\int_{\mathbb{R}^n\setminus B}\frac{f^2(x)}{|x-x_0|^{n+s}}dx\Big],\nonumber
\end{eqnarray}
where $C=C(n,s,\lambda,\Lambda,\tilde{\beta},\beta,\delta_0,q)$. Thus, Lemma 4.1 is proved.
\end{proof}

Following the approach in Section 5 of \cite{KMS2015} and Section 3 of \cite{BKK2024}, and employing Lemma 4.1 along with techniques from \cite{BR2013}, we establish the following reverse H\"{o}lder inequality (Theorem 4.2).

Let $\epsilon\in(0,1)$ and let
$$
\mu(\mathcal{D})=\int_{\mathcal{D}}\frac{1}{|x-y|^{n-2\epsilon}}dydx\quad\hbox{ for any }\mathcal{D}\subset\mathbb{R}^{2n}.
$$
This implies that for any $\mathcal{B}=B_r\times B_r$ with $B_r\subset\mathbb{R}^n$,
$$
\mu(\mathcal{B})=\int_B\int_B\frac{dydx}{|x-y|^{n-2\epsilon}}=\frac{c(n,\epsilon)}{\epsilon}r^{n+2\epsilon},\hbox{ with }c_1(n)\leq c(n,\epsilon)\leq c_2(n).
$$
For any $(x,y)\in\mathbb{R}^n\times\mathbb{R}^n$, let
$$
U(x,y)=\frac{|u(x)-u(y)|}{|x-y|^{s+\epsilon}}, \quad
G(x,y)=\frac{|g(x)-g(y)|}{|x-y|^{\sigma+\frac{2\epsilon}{q}}},\quad\hbox{ and }\quad F(x,y)=|f(x)|.
$$
For any $q\in (\max\{\frac{2n}{n+2s}, \frac{2n}{n+\delta_0}\}, 2)$, there exists $\delta_1\in(0,1)$ such that $W^{s+\delta_0,2}(B_r)\hookrightarrow W^{\sigma(1+\delta_1),q(1+\delta_1)}(B_r)$
whenever
$s+\delta_0-\frac{n}{2}=\sigma(1+\delta_1)-\frac{n}{q(1+\delta_1)}$, where $\sigma=s+\frac{\delta_0}{2}$.
Actually, the existence of $\delta_1$ can be easily verified if we notice that the equation $(s+\frac{\delta_0}{2})qx^2-(s+\delta_0-\frac{n}{2})qx-n=0$ has a root between 1 and  $\frac{s+\delta_0}{s+\frac{\delta_0}{2}}$ when $\frac{2n}{n+\delta_0}<q<2$ and $\delta_0\leq n-2s$. While for the case $\delta_0>n-2s$, obviously, we can find $\hat{\delta}_0<n-2s$ and $\delta_1$ such that $W^{s+\delta_0,2}(B_r)\hookrightarrow W^{s+\hat{\delta}_0,2}(B_r)\hookrightarrow W^{\sigma(1+\delta_1),q(1+\delta_1)}(B_r)$. Then by a direct  calculation, we can find that for any $\epsilon\in(0,\frac{q\sigma}{2}(1+\delta_1))$,
\begin{equation}
U\in L^2(\mathbb{R}^{2n},\mu), \quad
G\in L^{q(1+\delta_1)}(\mathcal{B},\mu),  \quad
F\in L^{2+\delta_0}(\mathcal{B},\mu),\nonumber\\
\end{equation}
whenever $u\in W^{s,2}(\mathbb{R}^n)$, $g\in W^{s+\delta_0,2}(B_r)$ and $f\in L^{2+\delta_0}(B_r)$.
\begin{thm}[Reverse H\"{o}lder inequality] Let $q\in (\max\{\frac{2n}{n+2s}, \frac{2n}{n+\delta_0}\}, 2)$ and $\epsilon\in(0,\epsilon_0)$ with $\epsilon_0=\min\{\frac{s}{2},\frac{(2s-\sigma)q}{4},\\ \frac{q\sigma}{2}(1+\delta_1)\}$.
Suppose $u\in W^{s,2}(\mathbb{R}^n)$ is a weak solution of
(\ref{4.+1})
with   $g\in W^{s,2}(\mathbb{R}^n)$ and $f\in L^2(\mathbb{R}^n)$. Let   $\mathcal{B}=\mathcal{B}(x_0,\rho_0)$,  where $B=B_{\rho_0}(x_0)\subset\Omega$ with $\rho_0>0$. If   $g\in W^{s+\delta_0,2}(B)$ and $f\in L^{2+\delta_0}(B)$ for some small $0<\delta_0<2s$,
then  there exist $\delta=\delta(n,s,\lambda,\Lambda,\tilde{\beta},\beta,\delta_0,q)>0$ and ${C}={C}(n,s,\lambda,\Lambda,\tilde{\beta},\beta,\delta_0,q,\rho_0)$ such that
\begin{eqnarray}
\Big(\fint_{\frac{1}{2}\mathcal{B}}U^{2+\delta})^{\frac{1}{2+\delta}}&\leq& C\Big[\Big(\fint_{\mathcal{B}}U^2d\mu\Big)^{\frac{1}{2}}
+\rho^{\sigma-s+\epsilon\frac{2-q}{q}}_0\Big(\fint_{\mathcal{B}}G^{q(1+\delta_1)}d\mu\Big)^{\frac{1}{q(1+\delta_1)}}
+\rho^{-\epsilon}_0\Big(\fint_{\mathcal{B}}F^{2+\delta_0}d\mu\Big)^{\frac{1}{2+\delta_0}}\nonumber\\
& &+\rho^{-\epsilon}_0\Big(\int_{\mathbb{R}^n\setminus B}\frac{[u(y)-\bar{u}_{B}]^2}{|y|^{n+2s}}dy\Big)^{\frac{1}{2}}
+\rho^{-\epsilon}_0\Big(\int_{\mathbb{R}^n\setminus B}\frac{|g(y)-\bar{g}_{B}|^{2}}{|y|^{n+2s}}dy\Big)^{\frac{1}{2}}\nonumber\\
& &
+\rho^{\frac{s}{2}-\epsilon}_0\Big(\int_{\mathbb{R}^n\setminus B}\frac{|f(y)|^2}{|y|^{n+s}}dy\Big)^{\frac{1}{2}}
\Big].\nonumber
\end{eqnarray}
\end{thm}
In order to prove Theorem 4.2, we need some lemmas.
\begin{lem} Let $q\in (\max\{\frac{2n}{n+2s}, \frac{2n}{n+\delta_0}\}, 2)$ and $\epsilon\in(0,\epsilon_0)$. Let $0<r<\rho_0\leq1$ and take a positive integer $m$ such that $\frac{\rho_0}{2}\leq2^mr<\rho_0$. Let $\mathcal{B}=B\times B$ with $B=B_r(x_0)\subset\Omega$. Suppose $u\in W^{s,2}(\mathbb{R}^n)$ is a weak solution of
(\ref{4.+1})
with  $g\in W^{s,2}(\mathbb{R}^n)\cap W^{s+\delta_0,2}(B)$   and $f\in L^2(\mathbb{R}^n)$.
Then we have
\begin{eqnarray}
\Big(\fint_{\frac{1}{2}\mathcal{B}}U^{2}d\mu\Big)^{\frac{1}{2}}
&\leq& C\Big[\frac{1}{\zeta\epsilon^{\frac{s}{n+2\epsilon}}}\Big(\fint_{\mathcal{B}}U^{\tilde{p}}d\mu\Big)^{\frac{1}{\tilde{p}}}
+\frac{\zeta}{\epsilon^{\frac{s}{n+2\epsilon}}}\sum^m_{j=1}2^{-j(s-\epsilon)}\Big(\fint_{2^j\mathcal{B}}U^{\tilde{p}}d\mu\Big)^{\frac{1}{\tilde{p}}}\nonumber\\
& &
+\zeta\epsilon^{\frac{1}{2}}r^{-\epsilon}\Big(\int_{\mathbb{R}^n\setminus 2^{m}B}\frac{[u(x)-\bar{u}_{2^mB}]^2}{|x-x_0|^{n+2s}}dx\Big)^{\frac{1}{2}}+\frac{r^{\sigma-s+\epsilon\frac{2-q}{q}}}{\epsilon^{\frac{1}{q}-\frac{1}{2}}}\sum^m_{j=0}2^{-j(2s-\sigma-\frac{2\epsilon}{q})}\Big(\fint_{2^j\mathcal{B}}G^qd\mu\Big)^{\frac{1}{q}}
\nonumber\\
& &+\epsilon^{\frac{1}{2}}r^{\sigma-s-\epsilon}\Big(\int_{\mathbb{R}^n\setminus 2^mB}\frac{|g(y)-\bar{g}_{2^mB}|^q}{|x_0-y|^{n+q\sigma}}dy\Big)^{\frac{1}{q}}+r^{-\epsilon}\sum^m_{j=0}2^{-\frac{js}{2}}\Big(\fint_{2^j\mathcal{B}}F^2d\mu\Big)^{\frac{1}{2}}
\nonumber\\
& &+\epsilon^{\frac{1}{2}}r^{\frac{s}{2}-\epsilon}\Big(\int_{\mathbb{R}^n\setminus 2^mB}\frac{|f(x)|^2}{|x-x_0|^{n+s}}dx\Big)^{\frac{1}{2}}
\Big],\nonumber
\end{eqnarray}
for any $\zeta\in(0,1)$, where $\tilde{p}=2(n+2\epsilon)/[n+2(s+\epsilon)]\in(0,2)$, $C=C(n,s,\lambda,\Lambda,\tilde{\beta},\beta,\delta_0,q)$.
\end{lem}
\begin{proof}
By Lemma 4.1,  we have
\begin{eqnarray}\label{4.+14}
&&\int_{\frac{1}{2}B}\int_{\frac{1}{2}B}\frac{[u(x)-u(y)]^2}{|x-y|^{n+2s}}dydx\nonumber\\
&\leq&C\Big[\frac{1}{\zeta^2 r^{2s}}\int_{B}[u(x)-\bar{u}_{B}]^2dx+\zeta^2r^{n+2s}\Big(\int_{\mathbb{R}^n\setminus B} \frac{|u(x)-\bar{u}_B|}{|x-x_0|^{n+2s}}dx\Big)^2\nonumber\\
& &+r^{n-2s+2\sigma-\frac{2n}{q}}\Big(\int_{B}\int_{B}\frac{|g(x)-g(y)|^q}{|x-y|^{n+q\sigma}}dydx\Big)^{\frac{2}{q}}
+r^{n+2s}\Big(\int_{\mathbb{R}^n\setminus B}\frac{|g(y)-\bar{g}_B|}{|x_0-y|^{n+2s}}dy\Big)^2\nonumber\\
& &+\int_Bf^2(x)dx+r^{n+s}\int_{\mathbb{R}^n\setminus B}\frac{f^2(x)}{|x-x_0|^{n+s}}dx\Big],
\end{eqnarray}
for any $\zeta\in(0,1)$, where $C=C(n,s,\lambda,\Lambda,\tilde{\beta},\beta,\delta_0,q)$.

For the left hand side of inequality (\ref{4.+14}), we have
\begin{equation}\label{4.+15}
\int_{\frac{1}{2}B}\int_{\frac{1}{2}B}\frac{[u(x)-u(y)]^2}{|x-y|^{n+2s}}dydx
=\int_{\frac{1}{2}\mathcal{B}}U^2d\mu.
\end{equation}

For the first term on the right hand side of (\ref{4.+14}),
applying Lemma 2.2 for $u\in W^{s,2}(B)\hookrightarrow W^{\tilde{s},\tilde{p}}(B)\hookrightarrow L^2(B)$ with $\tilde{s}=s+\epsilon-2\epsilon/{\tilde{p}}$ and $\tilde{p}=2(n+2\epsilon)/[n+2(s+\epsilon)]$, we have
\begin{eqnarray}\label{4.+16}
\frac{1}{r^{2s}}\int_B(u-\bar{u}_{B})^2dx&\leq& \frac{c(n,s)}{r^{2s}}\Big(\int_B\int_B\frac{|u(x)-u(y)|^{\tilde{p}}}{|x-y|^{n+\tilde{p}\tilde{s}}}dydx\Big)^{\frac{2}{\tilde{p}}}\nonumber\\
&=&\frac{c(n,s)}{r^{2s}}\Big(\int_{\mathcal{B}}U^{\tilde{p}}d\mu\Big)^{\frac{2}{\tilde{p}}}\nonumber\\
&\leq&\frac{C(n,s)\mu(\mathcal{B})}{\epsilon^{\frac{2s}{n+2\epsilon}}}\Big(\fint_{\mathcal{B}}U^{\tilde{p}}d\mu\Big)^{\frac{2}{\tilde{p}}}.
\end{eqnarray}
We estimate the second term on the right hand side of (\ref{4.+14}) by a similar process as that in \cite{BKK2024}. In fact, we have
\begin{eqnarray}
r^{n+2s}\Big(\int_{\mathbb{R}^n\setminus {B}}\frac{|u(x)-\bar{u}_B|}{|x-x_0|^{n+2s}}dx\Big)^2
=r^{n+2s}\Big(\sum_{i=1}^{m}\int_{2^iB\setminus 2^{i-1}B}\frac{|u(x)-\bar{u}_B|}{|x-x_0|^{n+2s}}dx
+\int_{\mathbb{R}^n\setminus 2^{m}B}\frac{|u(x)-\bar{u}_B|}{|x-x_0|^{n+2s}}dx\Big)^2.\nonumber
\end{eqnarray}
By a direct calculation, we get
\begin{eqnarray}
&&\int_{2^iB\setminus 2^{i-1}B}\frac{|u(x)-\bar{u}_B|}{|x-x_0|^{n+2s}}dx\nonumber\\
&\leq&\tilde{c}(n,s)(2^{i}r)^{-2s}\fint_{2^iB}|u(x)-\bar{u}_B|dx\nonumber\\
&\leq&c_0(n,s)(2^{i}r)^{-2s}\Big(\fint_{2^iB}|u(x)-\bar{u}_{2^iB}|dx+\sum^{i}_{j=1}|\bar{u}_{2^jB}-\bar{u}_{2^{j-1}B}|\Big)\nonumber\\
&\leq&c_0(n,s)(2^{i}r)^{-2s}\Big(\fint_{2^iB}|u(x)-\bar{u}_{2^iB}|dx+\sum^{i}_{j=1}\fint_{2^{j-1}B}[|u(x)-\bar{u}_{2^{j}B}|dx\Big)\nonumber\\
&\leq&c_0(n,s)(2^{i}r)^{-2s}\Big(\fint_{2^iB}|u(x)-\bar{u}_{2^iB}|dx+2^n\sum^{i}_{j=1}\Big(\fint_{2^{j}B}|u(x)-\bar{u}_{2^{j}B}|dx\Big)\nonumber\\
&\leq&c_1(n,s)(2^{i}r)^{-2s}\sum^{i}_{j=1}\fint_{2^{j}B}|u(x)-\bar{u}_{2^{j}B}|dx.\nonumber
\end{eqnarray}
Then by applying H\"{o}lder's inequality and Fractional Sobolev-Poincar\'{e} inequality, we obtain
\begin{eqnarray}
&&\sum_{i=1}^{m}\int_{2^iB\setminus 2^{i-1}B}\frac{|u(x)-\bar{u}_B|}{|x-x_0|^{n+2s}}dx\nonumber\\
&\leq&\sum_{i=1}^{m}c_1(n,s)(2^{i}r)^{-2s}\sum^{i}_{j=1}\fint_{2^{j}B}|u(x)-\bar{u}_{2^{j}B}|dx\nonumber\\
&\leq&\sum_{i=1}^{m}c_1(n,s)(2^{i}r)^{-2s}\sum^{i}_{j=1}\Big(\fint_{2^{j}B}|u(x)-\bar{u}_{2^{j}B}|^2dx\Big)^{\frac{1}{2}}\nonumber\\
&\leq&\sum_{i=1}^{m}c_1(n,s,\epsilon)(2^{i}r)^{-2s}\sum^{i}_{j=1}(2^jr)^{-\frac{n}{2}}\Big(\int_{2^jB}\int_{2^jB}\frac{|u(x)-u(y)|^{\tilde{p}}}{|x-y|^{n+\tilde{p}\tilde{s}}}dydx\Big)^{\frac{1}{\tilde{p}}}\nonumber\\
&\leq&C(n,s,\epsilon)\frac{\mu(\mathcal{B})^{\frac{1}{\tilde{p}}}}{r^{\frac{n}{2}+2s}}\sum^{m}_{j=1}2^{-j(s-\epsilon)}\Big(\fint_{2^j\mathcal{B}}U^{\tilde{p}}d\mu\Big)^{\frac{1}{\tilde{p}}}.\nonumber
\end{eqnarray}
Similarly, we can deduce that
\begin{eqnarray}
&&\int_{\mathbb{R}^n\setminus 2^{m}B}\frac{|u(x)-\bar{u}_B|}{|x-x_0|^{n+2s}}dx\nonumber\\
&\leq&\int_{\mathbb{R}^n\setminus 2^{m}B}\frac{|u(x)-\bar{u}_{2^mB}|}{|x-x_0|^{n+2s}}dx
+\int_{\mathbb{R}^n\setminus 2^{m}B}\frac{|\bar{u}_{2^mB}-\bar{u}_B|}{|x-x_0|^{n+2s}}dx\nonumber\\
&\leq&\int_{\mathbb{R}^n\setminus 2^{m}B}\frac{|u(x)-\bar{u}_{2^mB}|}{|x-x_0|^{n+2s}}dx
+c(n,s)(2^mr)^{-2s}|\bar{u}_{2^mB}-\bar{u}_B|\nonumber\\
&\leq&\int_{\mathbb{R}^n\setminus 2^{m}B}\frac{|u(x)-\bar{u}_{2^mB}|}{|x-x_0|^{n+2s}}dx
+C(n,s,\epsilon)\frac{\mu(\mathcal{B})^{\frac{1}{\tilde{p}}}}{r^{\frac{n}{2}+2s}}\sum^{m}_{j=1}2^{-j(s-\epsilon)}\Big(\fint_{2^j\mathcal{B}}U^{\tilde{p}}d\mu\Big)^{\frac{1}{\tilde{p}}}.\nonumber
\end{eqnarray}
Therefore, we obtain that
\begin{eqnarray}\label{4.+17}
&&r^{n+2s}\Big(\int_{\mathbb{R}^n\setminus {B}}\frac{|u(x)-\bar{u}_B|}{|x-x_0|^{n+2s}}dx\Big)^2\nonumber\\
&\leq&r^{n+2s}\Big[\int_{\mathbb{R}^n\setminus 2^{m}B}\frac{|u(x)-\bar{u}_{2^mB}|}{|x-x_0|^{n+2s}}dx
+C(n,s,\epsilon)\frac{\mu(\mathcal{B})^{\frac{1}{\tilde{p}}}}{r^{\frac{n}{2}+2s}}\sum^{m}_{j=1}2^{-j(s-\epsilon)}\Big(\fint_{2^j\mathcal{B}}U^{\tilde{p}}d\mu\Big)^{\frac{1}{\tilde{p}}}\Big]^2\nonumber\\
&\leq&c(n)\frac{\epsilon\mu(\mathcal{B})}{r^{2\epsilon}}\int_{\mathbb{R}^n\setminus2^m{B}}\frac{|u(x)-\bar{u}_{2^m{B}}|^2}{|x-x_0|^{n+2s}}dx
+C_1(n,s,\epsilon)\frac{\mu(\mathcal{B})}{\epsilon^{\frac{2s}{n+2\epsilon}}}\Big[\sum^{m}_{j=1}2^{-j(s-\epsilon)}\Big(\fint_{2^j\mathcal{B}}U^{\tilde{p}}d\mu\Big)^{\frac{1}{\tilde{p}}}\Big]^2.
\end{eqnarray}
Obviously,
\begin{eqnarray}\label{4.+19}
r^{n+2\sigma-2s-\frac{2n}{q}}\Big(\int_{B}\int_{B}\frac{|g(x)-g(y)|^q}{|x-y|^{n+q\sigma}}dydx\Big)^{\frac{2}{q}}
&=&r^{n+2\sigma-2s-\frac{2n}{q}}\Big(\int_{\mathcal{B}}G^qd\mu\Big)^{\frac{2}{q}}\nonumber\\
&\leq& C(n,s)\frac{r^{2\sigma-2s-2\epsilon+\frac{4\epsilon}{q}}}{\epsilon^{\frac{2}{q}-1}}\mu(\mathcal{B})\Big(\fint_{\mathcal{B}}G^qd\mu\Big)^{\frac{2}{q}}.
\end{eqnarray}
Similar to (\ref{4.+17}), we can derive that
\begin{eqnarray}\label{4.+20}
r^{n+2s}\Big(\int_{\mathbb{R}^n\setminus B}\frac{|g(y)-\bar{g}_B|}{|x_0-y|^{n+2s }}dy\Big)^2
&\leq&c(q)\epsilon r^{2\sigma-2s-2\epsilon}\mu(\mathcal{B})\Big(\int_{\mathbb{R}^n\setminus 2^mB}\frac{|g(y)-\bar{g}_{2^mB}|^q}{|x_0-y|^{n+q\sigma}}dy\Big)^{\frac{2}{q}}\nonumber\\
& &+C(n,q,\sigma)\frac{r^{2\sigma-2s-2\epsilon+\frac{4\epsilon}{q}}}{\epsilon^{\frac{2}{q}-1}}\mu(\mathcal{B})\Big[\sum^m_{j=1}2^{-j(2s-\sigma-\frac{2\epsilon}{q})}\Big(\fint_{2^j\mathcal{B}}G^qd\mu\Big)^{\frac{1}{q}}\Big]^2.
\end{eqnarray}
\begin{equation}\label{4.+21}
\int_{B}f^2(x)dx\leq c(n)r^n\fint_B\fint_{B}f^2(x)dydx\leq\frac{\hat{C}(n)\mu(\mathcal{B})}{r^{2\epsilon}}\fint_{\mathcal{B}}
{F}^2d\mu,
\end{equation}
and
\begin{eqnarray}\label{4.+22}
r^{n+s}\int_{\mathbb{R}^n\setminus B}\frac{|f(x)|^2}{|x-x_0|^{n+s}}dx
 \leq C(n)\epsilon r^{s-2\epsilon}\mu(\mathcal{B})\int_{\mathbb{R}^n\setminus 2^mB}\frac{|f(x)|^2}{|x-x_0|^{n+s}}dx
 +c(n,s)\frac{\mu(\mathcal{B})}{r^{2\epsilon}}\sum^m_{j=1}2^{-js}\fint_{2^j\mathcal{B}}
{F}^2d\mu.
\end{eqnarray}
Combining (\ref{4.+14})-(\ref{4.+22}), we obtain
\begin{eqnarray}
\Big(\fint_{\frac{1}{2}\mathcal{B}}U^{2}d\mu\Big)^{\frac{1}{2}}
&\leq& C\Big[\frac{1}{\zeta\epsilon^{\frac{s}{n+2\epsilon}}}\Big(\fint_{\mathcal{B}}U^{\tilde{p}}d\mu\Big)^{\frac{1}{\tilde{p}}}
+\frac{\zeta}{\epsilon^{\frac{s}{n+2\epsilon}}}\sum^m_{j=1}2^{-j(s-\epsilon)}\Big(\fint_{2^j\mathcal{B}}U^{\tilde{p}}d\mu\Big)^{\frac{1}{\tilde{p}}}\nonumber\\
& &
+\zeta\epsilon^{\frac{1}{2}}r^{-\epsilon}\Big(\int_{\mathbb{R}^n\setminus 2^{m}B}\frac{[u(x)-\bar{u}_{2^mB}]^2}{|x-x_0|^{n+2s}}dx\Big)^{\frac{1}{2}}+\frac{r^{\sigma-s+\epsilon\frac{2-q}{q}}}{\epsilon^{\frac{1}{q}-\frac{1}{2}}}\sum^m_{j=0}2^{-j(2s-\sigma-\frac{2\epsilon}{q})}\Big(\fint_{2^j\mathcal{B}}G^qd\mu\Big)^{\frac{1}{q}}
\nonumber\\
& &+\epsilon^{\frac{1}{2}}r^{\sigma-s-\epsilon}\Big(\int_{\mathbb{R}^n\setminus 2^mB}\frac{|g(y)-\bar{g}_{2^mB}|^q}{|x_0-y|^{n+q\sigma}}dy\Big)^{\frac{1}{q}}+r^{-\epsilon}\sum^m_{j=0}2^{-\frac{js}{2}}\Big(\fint_{2^j\mathcal{B}}F^2d\mu\Big)^{\frac{1}{2}}
\nonumber\\
& &+\epsilon^{\frac{1}{2}}r^{\frac{s}{2}-\epsilon}\Big(\int_{\mathbb{R}^n\setminus 2^mB}\frac{|f(x)|^2}{|x-x_0|^{n+s}}dx\Big)^{\frac{1}{2}}
\Big],\nonumber
\end{eqnarray}
where $C=C(n,s,\lambda,\Lambda,\tilde{\beta},\beta,\delta_0,q)$.
\end{proof}

Now we define that
\begin{eqnarray}\label{4.001}
\Phi_{\theta_1,\theta_2}(x,\rho)&=&\rho^{s+2\epsilon}\Big(\fint_{\mathcal{B}(x,\rho)}U^{2}d\mu\Big)^{\frac{1}{2}}
+\theta_1\frac{\rho^{\sigma+\epsilon\frac{(2+q)}{q}}}{\epsilon^{\frac{1}{q}-\frac{1}{2}}}\Big(\fint_{\mathcal{B}(x,\rho)}G^qd\mu\Big)^{\frac{1}{q}}
+\theta_2\rho^{s+\epsilon}\Big(\fint_{\mathcal{B}(x,\rho)}F^2d\mu\Big)^{\frac{1}{2}},
\end{eqnarray}
\begin{equation}\label{4.002}
\phi(x,\rho)=\Big(\fint_{\mathcal{B}(x,\rho)}F^{2+\delta_0}d\mu\Big)^{\frac{1}{2+\delta_0}}
+\Big(\fint_{\mathcal{B}(x,\rho)}G^{q(1+\delta_1)}d\mu\Big)^{\frac{1}{q(1+\delta_1)}},
\end{equation}
and
\begin{eqnarray}\label{4.003}
T(x,\rho)&=&\rho^{s+2\epsilon}\sum^m_{j=1}2^{-j(s-\epsilon)}\Big(\fint_{\mathcal{B}(x,2^j\rho)}U^{\tilde{p}}d\mu\Big)^{\frac{1}{\tilde{p}}}
+\epsilon^{\frac{1}{2}}\rho^{s+\epsilon}\Big(\int_{\mathbb{R}^n\setminus 2^{m}B_{\rho}(x)}\frac{[u(y)-\bar{u}_{2^mB_{\rho}(x)}]^2}{|y-x|^{n+2s}}dy\Big)^{\frac{1}{2}}
\nonumber\\
& &+\theta_1\frac{\rho^{\sigma+\epsilon\frac{2+q}{q}}}{\epsilon^{\frac{1}{q}-\frac{1}{2}}}\sum^m_{j=0}2^{-j(2s-\sigma-\frac{2\epsilon}{q})}\Big(\fint_{\mathcal{B}(x,2^j\rho)}G^qd\mu\Big)^{\frac{1}{q}}
+\theta_1\epsilon^{\frac{1}{2}}\rho^{\sigma+\epsilon}\Big(\int_{\mathbb{R}^n\setminus 2^mB_{\rho}(x)}\frac{|g(y)-\bar{g}_{2^mB_{\rho}(x)}|^q}{|y-x|^{n+q\sigma}}dy\Big)^{\frac{1}{q}}\nonumber\\
& &+\theta_2\rho^{s+\epsilon}\sum^m_{j=0}2^{-\frac{js}{2}}\Big(\fint_{2^j\mathcal{B}(x,2^j\rho)}F^2d\mu\Big)^{\frac{1}{2}}
+\theta_2\epsilon^{\frac{1}{2}}\rho^{\frac{3s}{2}+\epsilon}\Big(\int_{\mathbb{R}^n\setminus 2^mB_{\rho}(x)}\frac{|f(y)|^2}{|y-x|^{n+s}}dy\Big)^{\frac{1}{2}},
\end{eqnarray}
 for any $x\in\mathbb{R}^n$ and $\rho>0$, where $\theta_i$, $i=1,2$, are positive constants, $q\in (\max\{\frac{2n}{n+2s}, \frac{2n}{n+\delta_0}\}, 2)$ and $\epsilon\in(0,\epsilon_0)$.
\begin{lem}
Let $q\in (\max\{\frac{2n}{n+2s}, \frac{2n}{n+\delta_0}\}, 2)$, $\epsilon\in(0,\epsilon_0)$,  and $\mathcal{B}(x_0,\rho_0)=B_{\rho_0}(x_0)\times B_{\rho_0}(x_0)$ for some $B_{\rho_0}(x_0)\subset\Omega$ with $\rho_0\in(0,1]$. Let $\rho_0<\rho_1<\rho_2<3\rho_0/2$ and $t_0=\frac{c}{\epsilon}\Big(\frac{\rho_0}{\rho_2-\rho_1}\Big)^{2n}(\Phi_{\theta_1,\theta_2}+\phi+T)(x_0,2\rho_0)$ with $c=c(n,s,q,\delta_0)$.
Suppose $u\in W^{s,2}(\mathbb{R}^n)$  is a weak solution of (\ref{4.+1}) with $g\in W^{s,2}(\mathbb{R}^n) \cap W^{s+\delta_0,2}(B_{\rho_0}(x_0))$ and $f\in L^{2}(\mathbb{R}^n)\cap L^{2+\delta_0}(B_{\rho_0}(x_0))$ for some $\delta_0>0$. Then
there exist $\tau_u$, $\tau_g$, $\tau_{\bar{g}}$, $\tau_f\in(0,1)$  depending on $n$, $s$, $\lambda$, $\Lambda$, $\tilde{\beta}$, $\beta$,  $\delta_0$, $q$ and $\epsilon$  such that for any $t\geq t_0$, the following inequality holds
\begin{eqnarray}
\int_{\mathcal{B}_{\rho_1}\cap\{U>t\}}U^2d\mu
&\leq&C_1\epsilon^{-\frac{4s}{n+2\epsilon}}t^{2-\tilde{p}}\int_{\mathcal{B}_{\rho_2}\cap\{U> t\}}U^{\tilde{p}}d\mu
+C_2\Big[t^{2-q-\frac{q(1+\delta_1)[q(\sigma+\epsilon)+2\epsilon]}{n-q(\sigma+\epsilon)}}
t^{\frac{q(1+\delta_1)[q(\sigma+\epsilon)+2\epsilon]}{n-q(\sigma+\epsilon)}}_0\int_{\mathcal{B}_{\rho_2}\cap\{G>\tau_g t\}}G^qd\mu\nonumber\\
& &
+t^{-\frac{2(s+\epsilon)(2+\delta_0)}{n-2s}}t^{\frac{2(s+\epsilon)(2+\delta_0)}{n-2s}}_0\int_{\mathcal{B}_{\rho_2}\cap\{F>\tau_f t\}}F^{2}d\mu\Big],\nonumber
\end{eqnarray}
where $C_1=C_1(n,s,\lambda,\Lambda,\tilde{\beta},\beta,\delta_0,q)$ and $C_2=C_2(n,s,\lambda,\Lambda,\tilde{\beta},\beta,\delta_0,q,\epsilon)$.
\end{lem}
For  any $\tau\in(0,1]$, let
\begin{equation}\label{4.+01}
\tilde{t}_0=\frac{1}{\tau}\sup_{\rho\in[\frac{\rho_2-\rho_1}{40^n},\frac{\rho_0}{2}]}\sup_{x\in B_{\rho_1}(x_0)}\{(\Phi_{\theta_1,\theta_2}+\phi+T)(x,\rho)\}.
\end{equation}
For  any $t\geq \tilde{t}_0$, let
$$
D_{\tau t}=\{(x,x)\in\mathcal{B}(x_0,\rho_1)|\quad \sup_{0<\rho<\frac{\rho_2-\rho_1}{40^n}}\Phi_{\theta_1,\theta_2}(x,\rho)>\tau t \}.
$$
Then for any $(x,x)\in D_{\tau t}$, there exists $\rho(x)\in(0,\frac{\rho_2-\rho_1}{40^n})$ such that
\begin{equation}\label{4.+02}
\Phi_{\theta_1,\theta_2}(x,\rho(x))\geq\tau t\quad\hbox{ and }\sup_{\rho(x)<\rho<\frac{\rho_2-\rho_1}{40^n}}\Phi_{\theta_1,\theta_2}(x,\rho)\leq\tau t.
\end{equation}
Collecting all the balls $\{\mathcal{B}(x,2\rho(x))\}_{(x,x)\in D_{\tau t}}$ and applying Vitali covering lemma,  we obtain a family of mutually disjoint balls $\{\mathcal{B}(x_j,2\rho(x_j))\}_{J_D}$  such that
$$
\bigcup_{(x,x)\in D_{\tau t}}\mathcal{B}(x,2\rho(x))\subset\bigcup_{j\in J_D}\mathcal{B}(x_j,10\rho(x_j))\subset \mathcal{B}(x_0,\rho_2)
$$
where $J_D$ is a countable set. Denote $\mathcal{B}_j=\mathcal{B}(x_j,\rho(x_j))$ for any $j\in J_D$.
Let
$$
k_0=\Big[-log_2\Big(\frac{\rho_2-\rho_1}{n10^{10n}}\Big)\Big]+1
$$
where $[\cdot]$ denotes the integer part of a given number. For any $k\geq k_0$, let $\triangle_k$ be the disjoint collection of half-open cubes of side length $2^{-k}$ whose closures are touching $\bar{B}_{\frac{\rho_2+\rho_1}{2}}(x_0)$, i.e.,
$$
\triangle_k=\{x_0+2^{-k}v+[0,2^{-k})^n|\quad v\in\mathbb{Z}^n,x_0+2^{-k}v+[0,2^{-k}]^n\cap\bar{B}_{\frac{\rho_2+\rho_1}{2}}(x_0)\neq \emptyset\}.
$$
Obviously, there holds $B_{\rho_1}(x_0)\subset \cup_{K\in\triangle_k}K\subset B_{\rho_2}(x_0)$ for any $k\geq k_0$. Denote $\Xi_k=\{\mathcal{K}=K_1\times K_2|\quad K_1, K_2\in \triangle_k\}$ for any $k\geq k_0$.
Let
$$
t_1=\max\Big\{\tilde{t}_0,\sup_{\mathcal{K}\in\Xi_{k_0}}\Big(\fint_{\mathcal{K}}U^2d\mu\Big)^{\frac{1}{2}}\Big\}
\hbox{ and }\tau_0=c(n,s)\epsilon^{\frac{1}{\tilde{p}}}
$$
be  as the same as that in sections 5 of \cite{KMS2015}.

With the help of the Calder\'{o}n-Zygmund decomposition lemma, the following lemma was established in \cite{KMS2015}.
\begin{prop}(\cite{KMS2015})
For any $t\geq t_1$,
$$
\int_{\mathcal{B}_{\rho_1}\cap\{U>t\}}U^2d\mu
\leq (10n)^{n+2}\tau^2_0t^2\sum_{j\in J_D}\mu(\mathcal{B}_j)+c(n,s)t^{2-\tilde{p}}\int_{\mathcal{B}_{\rho_2}\cap\{U>\tau_0 t\}}U^{\tilde{p}}d\mu.
$$
where  $\mathcal{B}_{\rho_1}=\mathcal{B}(x_0,\rho_1)$, $\mathcal{B}_{\rho_2}=\mathcal{B}(x_0,\rho_2)$.
\end{prop}
\begin{proof}[Proof of Lemma 4.3.] For any $t\geq t_1$, by Proposition 4.1, we mainly need to estimate $\sum_{j\in J_D}\mu(\mathcal{B}_j)$.
By the definition  of $\mathcal{B}_j=\mathcal{B}(x_j,\rho(x_j))$, $\Phi_{\theta_1,\theta_2}(x_j,\rho(x_j))\geq\tau_0 t$, which implies that at least one of the following inequalities hold:
\begin{equation}\label{4.+23}
r^{s+2\epsilon}_j\Big(\fint_{\mathcal{B}_j}U^2d\mu\Big)^{\frac{1}{2}}\geq\frac{1}{3}\tau_0 t,
\end{equation}
\begin{equation}\label{4.+24}
\theta_1\frac{r^{\sigma+\epsilon\frac{(2+q)}{q}}_j}{\epsilon^{\frac{1}{q}-\frac{1}{2}}}\Big(\fint_{\mathcal{B}_j}G^qd\mu\Big)^{\frac{1}{q}}\geq\frac{1}{3}\tau_0 t,
\end{equation}
\begin{equation}\label{4.+25}
\theta_2r^{s+\epsilon}_j\Big(\fint_{\mathcal{B}_j}F^2d\mu\Big)^{\frac{1}{2}}\geq\frac{1}{3}\tau_0 t,
\end{equation}
where $r_j=\rho(x_j)$.

Suppose (\ref{4.+23}) holds. By Lemma 4.2, we have
\begin{eqnarray}
\frac{1}{3}\tau_0 t&\leq&r^{s+2\epsilon}_j\Big(\fint_{\mathcal{B}_j}U^{2}d\mu\Big)^{\frac{1}{2}}\nonumber\\
&\leq& C_0\Big[\frac{1}{\zeta\epsilon^{\frac{s}{n+2\epsilon}}}r^{s+2\epsilon}_j\Big(\fint_{2\mathcal{B}_j}U^{\tilde{p}}d\mu\Big)^{\frac{1}{\tilde{p}}}
+\frac{\zeta}{\epsilon^{\frac{s}{n+2\epsilon}}}r^{s+2\epsilon}_j\sum^{m+1}_{i=2}2^{-i(s-\epsilon)}\Big(\fint_{2^i\mathcal{B}_j}U^{\tilde{p}}d\mu\Big)^{\frac{1}{\tilde{p}}}
\nonumber\\
& &+\zeta\epsilon^{\frac{1}{2}}r^{s+\epsilon}_j\Big(\int_{\mathbb{R}^n\setminus 2^{m_j+1}B_j}\frac{[u(x)-\bar{u}_{2^{m_j+1}B}]^2}{|x-x_j|^{n+2s}}dx\Big)^{\frac{1}{2}}
+\frac{r^{\sigma+\epsilon\frac{2+q}{q}}_j}{\epsilon^{\frac{1}{q}-\frac{1}{2}}}\sum^{m_j+1}_{i=1}2^{-i(2s-\sigma-\frac{2\epsilon}{q})}\Big(\fint_{2^i\mathcal{B}_j}G^qd\mu\Big)^{\frac{1}{q}}
\nonumber\\
& &+\epsilon^{\frac{1}{2}}r^{\sigma+\epsilon}_j\Big(\int_{\mathbb{R}^n\setminus 2^{m_j+1}B_j}\frac{|g(x)-\bar{g}_{2^{m_j+1}B_j}|^q}{|x-x_j|^{n+q\sigma}}dx\Big)^{\frac{1}{q}}
+r^{s+\epsilon}_j\sum^{m_j+1}_{i=1}2^{-\frac{is}{2}}\Big(\fint_{2^i\mathcal{B}_j}F^2d\mu\Big)^{\frac{1}{2}}\nonumber\\
& &+\epsilon^{\frac{1}{2}}r^{\frac{3s}{2}+\epsilon}_j\Big(\int_{\mathbb{R}^n\setminus 2^{m_j+1}B_j}\frac{|f(x)|^2}{|x-x_j|^{n+s}}dx\Big)^{\frac{1}{2}}
\Big],\nonumber
\end{eqnarray}
for some positive integer $m_j$ such that $\frac{\rho_0}{2}\leq 2^{m_j+1}r_j<\rho_0$, where $C_0=C_0(n,s,\lambda,\Lambda,\tilde{\beta},\beta,\delta_0,q)$.
Notice that since $r_j\leq\frac{\rho_2-\rho_1}{40^n}$, we have $m_j\geq3$ and $\frac{\rho_2-\rho_1}{40^n}\leq\frac{\rho_0}{4}\leq2^{m_j}r_j<\frac{\rho_0}{2}$. Therefore, by H\"{o}lder's inequality and (\ref{4.+02}), and (\ref{4.+01}), respectively, we get
\begin{eqnarray}
\Big(\fint_{2^i\mathcal{B}_j}U^{\tilde{p}}d\mu\Big)^{\frac{1}{\tilde{p}}}\leq \tau_0t\quad\hbox{ for }i=2,\cdots,m_j-1,\nonumber
\end{eqnarray}
and
\begin{eqnarray}
r^{s+2\epsilon}_j\sum^{m_j+1}_{i=m_j}2^{-i(s-\epsilon)}\Big(\fint_{2^i\mathcal{B}_j}U^{\tilde{p}}d\mu\Big)^{\frac{1}{\tilde{p}}}
+\epsilon^{\frac{1}{2}}r^{s+\epsilon}_j\Big(\int_{\mathbb{R}^n\setminus 2^{m_j+1}B_j}\frac{[u(x)-\bar{u}_{2^{m_j+1}B}]^2}{|x-x_j|^{n+2s}}dx\Big)^{\frac{1}{2}}
&\leq&T(x_j,2^{m_j}r_j)\leq\tau_0\tilde{t}_0\leq\tau_0t,\nonumber
\end{eqnarray}
Similarly, we  can also get
\begin{eqnarray} \frac{r^{\sigma+\frac{\epsilon(2+q)}{q}}_j}{\epsilon^{\frac{1}{q}-\frac{1}{2}}}\Big(\fint_{2^i\mathcal{B}_j}G^qd\mu\Big)^{\frac{1}{q}}\leq\frac{\tau_0 t}{\theta_1}\quad\hbox{ for }i=1,2,\cdots,m_j-1,\nonumber
\end{eqnarray}
\begin{eqnarray}
\frac{r^{\sigma+\epsilon\frac{2+q}{q}}_j}{\epsilon^{\frac{1}{q}-\frac{1}{2}}}\sum^{m_j+1}_{i=m_j}2^{-i(2s-\sigma-\frac{2\epsilon}{q})}\Big(\fint_{2^i\mathcal{B}_j}G^qd\mu\Big)^{\frac{1}{q}}
+\epsilon^{\frac{1}{2}}r^{\sigma+\epsilon}_j\Big(\int_{\mathbb{R}^n\setminus 2^{m_j+1}B_j}\frac{|g(x)-\bar{g}_{2^{m_j+1}B_j}|^q}{|x-x_j|^{n+q\sigma}}dx\Big)^{\frac{1}{q}}
&\leq&\frac{T(x_j,2^{m_j}r_j)}{\theta_1}\leq\frac{\tau_0\tilde{t}_0}{\theta_1}\leq\frac{\tau_0t}{\theta_1},\nonumber
\end{eqnarray}
and
\begin{eqnarray}
r^{s+\epsilon}_j\Big(\fint_{2^i\mathcal{B}_j}F^2d\mu\Big)^{\frac{1}{2}}\leq\frac{\tau_0 t}{\theta_2}\quad\hbox{ for }i=1,2,\cdots,m_j-1,\nonumber
\end{eqnarray}
\begin{eqnarray}
r^{s+\epsilon}_j\sum^{m_j+1}_{i=m_j}2^{-\frac{is}{2}}\Big(\fint_{2^i\mathcal{B}_j}F^2d\mu\Big)^{\frac{1}{2}}
+\epsilon^{\frac{1}{2}}r^{\frac{3s}{2}+\epsilon}_j\Big(\int_{\mathbb{R}^n\setminus 2^{m_j+1}B_j}\frac{|f(x)|^2}{|x-x_j|^{n+s}}dx\Big)^{\frac{1}{2}}
&\leq&\frac{T(x_j,2^{m_j}r_j)}{\theta_2}\leq\frac{\tau_0\tilde{t}_0}{\theta_2}\leq\frac{\tau_0t}{\theta_2},\nonumber
\end{eqnarray}
Combining the above inequalities, we can derive that
\begin{eqnarray}
\frac{1}{3}\tau_0 t&\leq&C_0\frac{1}{\zeta\epsilon^{\frac{s}{n+2\epsilon}}}r^{s+2\epsilon}_j\Big(\fint_{2\mathcal{B}_j}U^{\tilde{p}}d\mu\Big)^{\frac{1}{\tilde{p}}}
+\frac{C\zeta}{\epsilon^{\frac{s}{n+2\epsilon}}}\tau_0t+\frac{C_1\tau_0 t}{\theta_1}+\frac{C_2\tau_0 t}{\theta_2}
\nonumber
\end{eqnarray}
where $C$ and $C_i$, $i=0,1,2$, are positive constants depending on $n$, $s$, $\lambda$, $\Lambda$, $\tilde{\beta}$, $\beta$, $\delta_0$ and $q$.
Taking
\begin{equation}
\zeta=\frac{\epsilon^{\frac{s}{n+2\epsilon}}}{18C},\quad \theta_1=18C_1,\quad \theta_2=18C_2,\nonumber
\end{equation}
we deduce that
\begin{equation}
\frac{1}{6}\tau_0 t\leq \hat{C}\epsilon^{-\frac{2s}{n+2\epsilon}}\Big(\fint_{2\mathcal{B}_j}U^{\tilde{p}}d\mu\Big)^{\frac{1}{\tilde{p}}},\nonumber
\end{equation}
which implies that
\begin{eqnarray}
\mu(\mathcal{B}_j)&\leq& \frac{\tilde{C}\epsilon^{-\frac{2s\tilde{p}}{n+2\epsilon}}}{(\tau_0 t)^{\tilde{p}}}\int_{2\mathcal{B}_j}U^{\tilde{p}}d\mu\nonumber\\
&=&\frac{\tilde{C}\epsilon^{-\frac{2s\tilde{p}}{n+2\epsilon}}}{(\tau_0 t)^{\tilde{p}}}\Big[\int_{2\mathcal{B}_j\cap\{U\leq \tilde{\tau}\tau_0 t\}}U^{\tilde{p}}d\mu+\int_{\mathcal{B}_j\cap\{U>\tilde{\tau}\tau_0 t\}}U^{\tilde{p}}d\mu\Big]\nonumber\\
&\leq&\frac{\tilde{C}\epsilon^{-\frac{2s\tilde{p}}{n+2\epsilon}}}{(\tau_0 t)^{\tilde{p}}}\Big[\mu(2\mathcal{B}_j)(\tilde{\tau}\tau_0 t)^{\tilde{p}}+\int_{2\mathcal{B}_j\cap\{U>\tilde{\tau}\tau_0 t\}}U^{\tilde{p}}d\mu\Big],\nonumber
\end{eqnarray}
 where $\tilde{C}=\tilde{C}(n,s,\lambda,\Lambda,\tilde{\beta},\beta,\delta_0,q)$. Take $\tilde{\tau}=2^{-\frac{n+3}{\tilde{p}}}\tilde{C}^{-\frac{1}{\tilde{p}}}\epsilon^{\frac{2s}{n+2\epsilon}}$. Then we get
\begin{equation}
\mu(\mathcal{B}_j)\leq\frac{2\tilde{C}\epsilon^{-\frac{2s\tilde{p}}{n+2\epsilon}}}{(\tau_0 t)^{\tilde{p}}}\int_{2\mathcal{B}_j\cap\{U>\tilde{\tau}\tau_0 t\}}U^{\tilde{p}}d\mu.\nonumber
\end{equation}

If (\ref{4.+24}) is true,  we can derive that
\begin{eqnarray}
\mu(\mathcal{B}_j)^{\frac{n-(\sigma+\epsilon)q}{n+2\epsilon}}
&\leq&
 \Big(\frac{\epsilon}{c(n,\epsilon)}\Big)^{\frac{\sigma q+\epsilon(2+q)}{n+2\epsilon}}\frac{1}{\epsilon^{1-\frac{q}{2}}}\Big(\frac{3\theta_1}{\tau_0 t}\Big)^{q}
\int_{\mathcal{B}_j}G^qd\mu\nonumber\\
&=&\Big(\frac{\epsilon}{c(n,\epsilon)}\Big)^{\frac{\sigma q+\epsilon(2+q)}{n+2\epsilon}}\frac{1}{\epsilon^{1-\frac{q}{2}}}\Big(\frac{3\theta_1}{\tau_0 t}\Big)^{q}\Big[\int_{\mathcal{B}_j\cap\{G\leq \tau_1\tau_0 t\}}G^qd\mu
+\int_{\mathcal{B}_j\cap\{G>\tau_1\tau_0 t\}}G^qd\mu\Big]\nonumber\\
&\leq& \Big(\frac{\epsilon}{c(n,\epsilon)}\Big)^{\frac{\sigma q+\epsilon(2+q)}{n+2\epsilon}}\frac{1}{\epsilon^{1-\frac{q}{2}}}\Big(\frac{3\theta_1}{\tau_0 t}\Big)^{q}
\Big[\mu(\mathcal{B}_j)(\tau_1\tau_0 t)^{q}
+\int_{\mathcal{B}_j\cap\{G>\tau_1\tau_0 t\}}G^qd\mu\Big]\nonumber\\
&\leq& \frac{(3\theta_1)^{q}}{\epsilon^{1-\frac{q}{2}}}
\tau_1^{q}\mu(\mathcal{B}_j)^{\frac{n-(\sigma+\epsilon)q}{n+2\epsilon}}
+\Big(\frac{\epsilon}{c(n,\epsilon)}\Big)^{\frac{\sigma q+\epsilon(2+q)}{n+2\epsilon}}\frac{1}{\epsilon^{1-\frac{q}{2}}}\Big(\frac{3\theta_1}{\tau_0 t}\Big)^{q}\int_{\mathcal{B}_j\cap\{G>\tau_1\tau_0 t\}}G^qd\mu
,\nonumber
\end{eqnarray}
where we used the fact $\mu(\mathcal{B}_j)\leq \mu(\mathcal{B}(x_0,\rho_0))\leq c(n,\epsilon)\epsilon^{-1}$. Taking $\tau_1=\epsilon^{\frac{1}{q}-\frac{1}{2}}2^{-\frac{1}{q}}(3\theta_1)^{-1}$,  we get
\begin{equation}
\mu(\mathcal{B}_j)^{\frac{n-(\sigma+\epsilon)q}{n+2\epsilon}}\leq\frac{c(n)}{\epsilon^{1-\frac{q}{2}}}\Big(\frac{3\theta_1}{\tau_0 t}\Big)^{q}\int_{\mathcal{B}_j\cap\{G>\tau_1\tau_0 t\}}G^qd\mu,\nonumber\\
\end{equation}
that is,
\begin{equation}
\mu(\mathcal{B}_j)\leq \Big(\frac{c(n)}{\epsilon^{1-\frac{q}{2}}}\Big)^{\frac{n+2\epsilon}{n-(\sigma+\epsilon)q}}\Big(\frac{3\theta_1}{\tau_0 t}\Big)^{\frac{q(n+2\epsilon)}{n-(\sigma+\epsilon)q}}\Big(\int_{\mathcal{B}_j\cap\{G>\tau_1\tau_0 t\}}G^qd\mu\Big)^{\frac{n+2\epsilon}{n-(\sigma+\epsilon)q}}.\nonumber\\
\end{equation}
Notice that by the definition of $\mathcal{B}_j$ and $\tilde{t}_0$,  there exists positive integer $m$ such that $\frac{\rho_2-\rho_1}{40^n}\leq 2^mr_j\leq\frac{\rho_0}{2}$. Therefore,
\begin{eqnarray}
\int_{\mathcal{B}_j\cap\{G>\tau_1\tau_0 t\}}G^qd\mu
&=&\int_{\mathcal{B}_j\cap\{G>\tau_1\tau_0 t\}}\Big(\frac{G}{\tau_1\tau_0 t}\Big)^q(\tau_1\tau_0 t)^qd\mu\nonumber\\
&\leq&\int_{\mathcal{B}_j\cap\{G>\tau_1\tau_0 t\}}\Big(\frac{G}{\tau_1\tau_0 t}\Big)^{q(1+\delta_1)}(\tau_1\tau_0 t)^qd\mu\nonumber\\
&\leq&(\tau_1\tau_0 t)^{-q\delta_1}\int_{\mathcal{B}_j}G^{q(1+\delta_1)}d\mu\nonumber\\
&\leq&(\tau_1\tau_0 t)^{-q\delta_1}\mu(2^{m}\mathcal{B}_j)\fint_{2^{m}\mathcal{B}_j}G^{q(1+\delta_1)}d\mu\nonumber\\
&\leq&(\tau_1\tau_0 t)^{-q\delta_1}\mu(2^{m}\mathcal{B}_j)\tilde{t}^{q(1+\delta_1)}_0.\nonumber
\end{eqnarray}
Therefore, we obtain
\begin{eqnarray}
\mu(\mathcal{B}_j)&\leq&\Big(\frac{c(n)}{\epsilon^{1-\frac{q}{2}}}\Big)^{\frac{n+2\epsilon}{n-(\sigma+\epsilon)q}}\Big(\frac{3\theta_1}{\tau_0 t}\Big)^{\frac{q(n+2\epsilon)}{n-(\sigma+\epsilon)q}}\int_{\mathcal{B}_j\cap\{G>\tau_1\tau_0 t\}}G^qd\mu\Big(\int_{\mathcal{B}_j\cap\{G>\tau_1\tau_0 t\}}G^qd\mu\Big)^{\frac{\sigma q+\epsilon(2+q)}{n-(\sigma+\epsilon)q}}\nonumber\\
&\leq& C_1(n,q,\sigma,\theta_3,\epsilon)(\tau_0 t)^{-q-\frac{q(1+\delta_1)[q(\sigma+\epsilon)+2\epsilon]}{n-q(\sigma+\epsilon)}}
\tilde{t}^{\frac{q(1+\delta_1)[q(\sigma+\epsilon)+2\epsilon]}{n-q(\sigma+\epsilon)}}_0\int_{\mathcal{B}_j\cap\{G>\tau_1\tau_0 t\}}G^qd\mu.\nonumber
\end{eqnarray}
Similarly, when (\ref{4.+25}) holds, we can derive that
\begin{equation}
\mu(\mathcal{B}_j)\leq C_2(n,s,\theta_2)(\tau_0 t)^{-2-\frac{2(s+\epsilon)(2+\delta_0)}{n-2s}}\tilde{t}^{\frac{2(s+\epsilon)(2+\delta_0)}{n-2s}}_0\int_{\mathcal{B}_j\cap\{F>\tau_2\tau_0 t\}}F^2d\mu.\nonumber
\end{equation}
where  $\tau_2=2^{-\frac{1}{2}}(3\theta_2)^{-1}$.
Consequently, we always have
\begin{eqnarray}
\mu(\mathcal{B}_j)&\leq& C\epsilon^{-\frac{2s\tilde{p}}{n+2\epsilon}}(\tau_0 t)^{-\tilde{p}}\int_{2\mathcal{B}_j\cap\{U>\tilde{\tau}\tau_0 t\}}U^{\tilde{p}}d\mu
+\tilde{C}\Big[(\tau_0 t)^{-q-\frac{q(1+\delta_1)[q(\sigma+\epsilon)+2\epsilon]}{n-q(\sigma+\epsilon)}}
\tilde{t}^{\frac{q(1+\delta_1)[q(\sigma+\epsilon)+2\epsilon]}{n-q(\sigma+\epsilon)}}_0\int_{\mathcal{B}_j\cap\{G>\tau_1\tau_0 t\}}G^qd\mu\nonumber\\
& &
+(\tau_0 t)^{-2-\frac{2(s+\epsilon)(2+\delta_0)}{n-2s}}\tilde{t}^{\frac{2(s+\epsilon)(2+\delta_0)}{n-2s}}_0\int_{\mathcal{B}_j\cap\{F>\tau_2\tau_0 t\}}F^{2}d\mu
\Big],\nonumber
\end{eqnarray}
where $C=C(n,s,\lambda,\Lambda,\tilde{\beta},\beta,\delta_0,q)$ and $\tilde{C}=\tilde{C}(n,s,\lambda,\Lambda,\tilde{\beta},\beta,\delta_0,q,\epsilon)$. This together with the facts that $\{2\mathcal{B}_j\}$ is a disjoint family  and all members belong to $\mathcal{B}_{\rho_2}$ imply that
\begin{eqnarray}
\sum_{j\in J_D}\mu(\mathcal{B}_j)&\leq& C\epsilon^{-\frac{2s\tilde{p}}{n+2\epsilon}}(\tau_0 t)^{-\tilde{p}}\int_{\mathcal{B}_{\rho_2}\cap\{U>\tilde{\tau}\tau_0 t\}}U^{\tilde{p}}d\mu
+\tilde{C}\Big[(\tau_0 t)^{-q-\frac{q(1+\delta_1)[q(\sigma+\epsilon)+2\epsilon]}{n-q(\sigma+\epsilon)}}
\tilde{t}^{\frac{q(1+\delta_1)[q(\sigma+\epsilon)+2\epsilon]}{n-q(\sigma+\epsilon)}}_0\int_{\mathcal{B}_{\rho_2}\cap\{G>\tau_1\tau_0 t\}}G^qd\mu\nonumber\\
& &
+(\tau_0 t)^{-2-\frac{2(s+\epsilon)(2+\delta_0)}{n-2s}}\tilde{t}^{\frac{2(s+\epsilon)(2+\delta_0)}{n-2s}}_0\int_{\mathcal{B}_{\rho_2}\cap\{F>\tau_2\tau_0 t\}}F^{2}d\mu\Big].\nonumber
\end{eqnarray}
Applying Proposition 4.1, we obtain that
\begin{eqnarray}
\int_{\mathcal{B}_{\rho_1}\cap\{U>t\}}U^2d\mu&\leq& C_1\epsilon^{-\frac{4s}{n+2\epsilon}}(\tilde{\tau}\tau_0 t)^{2-\tilde{p}}\int_{\mathcal{B}_{\rho_2}\cap\{U>\tilde{\tau}\tau_0 t\}}U^{\tilde{p}}d\mu\nonumber\\
& &+\tilde{C}_1\Big[(\tau_0 t)^{2-q-\frac{q(1+\delta_1)[q(\sigma+\epsilon)+2\epsilon]}{n-q(\sigma+\epsilon)}}
\tilde{t}^{\frac{q(1+\delta_1)[q(\sigma+\epsilon)+2\epsilon]}{n-q(\sigma+\epsilon)}}_0\int_{\mathcal{B}_{\rho_2}\cap\{G>\tau_1\tau_0 t\}}G^qd\mu\nonumber\\
& &
+(\tau_0 t)^{-\frac{2(s+\epsilon)(2+\delta_0)}{n-2s}}\tilde{t}^{\frac{2(s+\epsilon)(2+\delta_0)}{n-2s}}_0\int_{\mathcal{B}_{\rho_2}\cap\{F>\tau_2\tau_0 t\}}F^{2}d\mu\Big],\nonumber
\end{eqnarray}
which together with the following inequality
\begin{equation}
\int_{\mathcal{B}_{\rho_1}\cap\{U>\tilde{\tau}\tau_0 t\}}U^2d\mu
\leq t^{2-\tilde{p}}\int_{\mathcal{B}_{\rho_1}\cap\{U>\tilde{\tau}\tau_0 t\}}U^{\tilde{p}}d\mu+\int_{\mathcal{B}_{\rho_1}\cap\{U> t\}}U^2d\mu,\nonumber\\
\end{equation}
implies that
\begin{eqnarray}
\int_{\mathcal{B}_{\rho_1}\cap\{U>\tilde{\tau}\tau_0 t\}}U^2d\mu&\leq& \hat{C}_1\epsilon^{-\frac{4s}{n+2\epsilon}}(\tilde{\tau}\tau_0 t)^{2-\tilde{p}}\int_{\mathcal{B}_{\rho_2}\cap\{U>\tilde{\tau}\tau_0 t\}}U^{\tilde{p}}d\mu\nonumber\\
& &+\hat{\tilde{C}}_1\Big[(\tau_0 t)^{2-q-\frac{q(1+\delta_1)[q(\sigma+\epsilon)+2\epsilon]}{n-q(\sigma+\epsilon)}}\tilde{t}^{\frac{q(1+\delta_1)[q(\sigma+\epsilon)+2\epsilon]}{n-q(\sigma+\epsilon)}}_0\int_{\mathcal{B}_{\rho_2}\cap\{G>\tau_1\tau_0 t\}}G^qd\mu\nonumber\\
& &
+(\tau_0 t)^{-\frac{2(s+\epsilon)(2+\delta_0)}{n-2s}}\tilde{t}^{\frac{2(s+\epsilon)(2+\delta_0)}{n-2s}}_0\int_{\mathcal{B}_{\rho_2}\cap\{F>\tau_2\tau_0 t\}}F^{2}d\mu\Big]
,\nonumber
\end{eqnarray}
for any $t\geq t_1$. This can be reformulated as
\begin{eqnarray}
\int_{\mathcal{B}_{\rho_1}\cap\{U>t\}}U^2d\mu
&\leq&C_1\epsilon^{-\frac{4s}{n+2\epsilon}}t^{2-\tilde{p}}\int_{\mathcal{B}_{\rho_2}\cap\{U> t\}}U^{\tilde{p}}d\mu
+C_2\Big[t^{2-q-\frac{q(1+\delta_1)[q(\sigma+\epsilon)+2\epsilon]}{n-q(\sigma+\epsilon)}}
\tilde{t}^{\frac{q(1+\delta_1)[q(\sigma+\epsilon)+2\epsilon]}{n-q(\sigma+\epsilon)}}_0\int_{\mathcal{B}_{\rho_2}\cap\{G>\tau_g t\}}G^qd\mu
\nonumber\\
& &
+t^{-\frac{2(s+\epsilon)(2+\delta_0)}{n-2s}}\tilde{t}^{\frac{2(s+\epsilon)(2+\delta_0)}{n-2s}}_0\int_{\mathcal{B}_{\rho_2}\cap\{F>\tau_f t\}}F^{2}d\mu\Big],\nonumber
\end{eqnarray}
for any $t\geq \tilde{\tau}\tau_0 t_1$, where $\tau_{g}=\tau_1/\tilde{\tau}$, $\tau_f=\tau_2/\tilde{\tau}$, $C_1=C_1(n,s,\lambda,\Lambda,\tilde{\beta},\beta,\delta_0,q)$ and
 $C_2=C_2(n,s,\lambda,\Lambda,\tilde{\beta},\beta,\delta_0,q,\epsilon)$.

Now by a similar process as that in \cite{KMS2015}, we can show that the above inequality holds for any $t\geq t_0$.
Notice that for any $x\in B(x_0,\rho_1)$ and $\rho\in[\frac{\rho_2-\rho_1}{40^n},\frac{\rho_0}{2}]$, we have $B(x,\rho)\subset B(x_0,2\rho_0)$.
Therefore, if $H$ is a $\mu$-integrable function, then
\begin{equation}
\fint_{\mathcal{B}(x,\rho)}Hd\mu\leq\frac{\mu(\mathcal{B}(x_0,2\rho_0))}{\mu(\mathcal{B}(x,\rho))}\fint_{\mathcal{B}(x_0,2\rho)}Hd\mu\leq c(n)\Big(\frac{\rho_0}{\rho_2-\rho_1}\Big)^{n+2\epsilon}\fint_{\mathcal{B}(x_0,2\rho_0)}Hd\mu.\nonumber
\end{equation}
By a direct calculation, we get
\begin{equation}
(\Phi_{\theta_1,\theta_2}+\phi)(x,\rho)\leq c(n)\Big(\frac{\rho_0}{\rho_2-\rho_1}\Big)^{\frac{n+2\epsilon}{q}}(\Phi_{\theta_1,\theta_2}+\phi)(x_0,2\rho_0).\nonumber\\
\end{equation}
and
\begin{equation}
\rho^{s+2\epsilon}\sum^m_{j=1}2^{-j(s-\epsilon)}\Big(\fint_{\mathcal{B}(x,2^j\rho)}U^{\tilde{p}}d\mu\Big)^{\frac{1}{\tilde{p}}}
\leq c(n,s)\Big(\frac{\rho_0}{\rho_2-\rho_1}\Big)^{\frac{n+2\epsilon}{\tilde{p}}}(2\rho_0)^{s+2\epsilon}\Big(\fint_{\mathcal{B}(x_0,2\rho_0)}U^{\tilde{p}}d\mu\Big)^{\frac{1}{\tilde{p}}}
\end{equation}
By a similar process as that in the proof of Lemma 4.2, we have
\begin{eqnarray}
&&\int_{\mathbb{R}^n\setminus 2^{m}B_{\rho}(x)}\frac{[u(y)-\bar{u}_{2^mB_{\rho}(x)}]^2}{|y-x|^{n+2s}}dy\nonumber\\
&\leq&\int_{(\mathbb{R}^n\setminus 2^{m}B_{\rho}(x))\cap B_{2\rho_0}(x_0)}\frac{[u(y)-\bar{u}_{2^mB_{\rho}(x)}]^2}{|y-x|^{n+2s}}dy
+\int_{(\mathbb{R}^n\setminus B_{\rho}(x))\cap(\mathbb{R}^n\setminus B_{2\rho_0}(x_0))}\frac{[u(y)-\bar{u}_{2^mB_{\rho}(x)}]^2}{|y-x|^{n+2s}}dy\nonumber\\
&\leq&2\int_{(\mathbb{R}^n\setminus 2^{m}B_{\rho}(x))\cap B_{2\rho_0}(x_0)}\frac{[u(y)-\bar{u}_{B_{2\rho_0}(x_0)}]^2}{|y-x|^{n+2s}}dy
+2\int_{(\mathbb{R}^n\setminus 2^{m}B_{\rho}(x))\cap B_{2\rho_0}(x_0)}\frac{[\bar{u}_{B_{2\rho_0}(x_0)}-\bar{u}_{2^mB_{\rho}(x)}]^2}{|y-x|^{n+2s}}dy\nonumber\\
& &+2\int_{(\mathbb{R}^n\setminus B_{\rho}(x))\cap(\mathbb{R}^n\setminus B_{2\rho_0}(x_0))}\frac{[u(y)-\bar{u}_{B_{2\rho_0}(x_0)}]^2}{|y-x|^{n+2s}}dy
+2\int_{(\mathbb{R}^n\setminus B_{\rho}(x))\cap(\mathbb{R}^n\setminus B_{2\rho_0}(x_0))}\frac{[\bar{u}_{B_{2\rho_0}(x_0)}-\bar{u}_{2^mB_{\rho}(x)}]^2}{|y-x|^{n+2s}}dy\nonumber\\
&\leq&\frac{c(n,s)}{(2^m\rho)^{n+2s}}\int_{B_{2\rho_0}(x_0)}[u(y)-\bar{u}_{B_{2\rho_0}(x_0)}]^2dy+\frac{c(n,s)}{(2^m\rho)^{2s}}[\bar{u}_{2^mB_\rho(x)}-\bar{u}_{B_{2\rho_0}(x_0)}]^2\nonumber\\
& &+c(n,s)\Big(\frac{\rho_0}{\rho_2-\rho_1}\Big)^{n+2s}\int_{\mathbb{R}^n\setminus B_{2\rho_0}(x_0)}\frac{[u(y)-\bar{u}_{B_{2\rho}(x_0)}]^2}{|y-x_0|^{n+2s}}dy\nonumber\\
&\leq&c(n,s)\Big(\frac{1}{(2^m\rho)^{n+2s}}\int_{ B_{2\rho_0}(x_0)}[u(y)-\bar{u}_{B_{2\rho_0}(x_0)}]^2dy
+\Big(\frac{\rho_0}{\rho_2-\rho_1}\Big)^{n+2s}\int_{\mathbb{R}^n\setminus B_{2\rho_0}(x_0)}\frac{[u(y)-\bar{u}_{B_{2\rho_0}(x_0)}]^2}{|y-x_0|^{n+2s}}dy\Big)\nonumber\\
&\leq&c(n,s)\Big(\frac{\rho_0}{\rho_2-\rho_1}\Big)^{n+2s}\Big[\Big(\fint_{\mathcal{B}_{2\rho_0}}U^{\tilde{p}}d\mu\Big)^{\frac{2}{\tilde{p}}}
+\int_{\mathbb{R}^n\setminus B_{2\rho_0}(x_0)}\frac{[u(y)-\bar{u}_{B_{2\rho_0}(x_0)}]^2}{|y-x_0|^{n+2s}}dy\Big].\nonumber
\end{eqnarray}
We can also get the similar estimates for the integrals involving $G$ and $F$ in $T(x,\rho)$,
\begin{eqnarray}
&&\theta_1\frac{\rho^{\sigma+\epsilon\frac{2+q}{q}}}{\epsilon^{\frac{1}{q}-\frac{1}{2}}}\sum^m_{j=0}2^{-j(2s-\sigma-\frac{2\epsilon}{q})}\Big(\fint_{\mathcal{B}(x,2^j\rho)}G^qd\mu\Big)^{\frac{1}{q}}
+\theta_1\epsilon^{\frac{1}{2}}\rho^{\sigma+\epsilon}\Big(\int_{\mathbb{R}^n\setminus 2^mB_{\rho}(x)}\frac{|g(y)-\bar{g}_{2^mB_{\rho}(x)}|^q}{|y-x|^{n+q\sigma}}dy\Big)^{\frac{1}{q}}\nonumber\\
&\leq&c(n,q,\sigma)\Big(\frac{\rho_0}{\rho_2-\rho_1}\Big)^{\frac{n+q\sigma}{q}}\Big[\theta_1\frac{(2\rho_0)^{\sigma+\epsilon\frac{2+q}{q}}}{\epsilon^{\frac{1}{q}-\frac{1}{2}}}\Big(\fint_{\mathcal{B}(x_0,2\rho_0)}G^qd\mu\Big)^{\frac{1}{q}}
+\theta_1\epsilon^{\frac{1}{2}}(2\rho_0)^{\sigma+\epsilon}\Big(\int_{\mathbb{R}^n\setminus B_{2\rho_0}(x_0)}\frac{|g(y)-\bar{g}_{B_{2\rho_0}(x_0)}|^q}{|y-x_0|^{n+q\sigma}}dy\Big)^{\frac{1}{q}}\Big],\nonumber
\end{eqnarray}
and
\begin{eqnarray}
&&\theta_2\rho^{s+\epsilon}\sum^m_{j=0}2^{-\frac{js}{2}}\Big(\fint_{2^j\mathcal{B}(x,2^j\rho)}F^2d\mu\Big)^{\frac{1}{2}}
+\theta_2\epsilon^{\frac{1}{2}}\rho^{\frac{3s}{2}+\epsilon}\Big(\int_{\mathbb{R}^n\setminus 2^mB_{\rho}(x)}\frac{|f(y)|^2}{|y-x|^{n+s}}dy\Big)^{\frac{1}{2}}\nonumber\\
& &\leq c(n,s)\Big(\frac{\rho_0}{\rho_2-\rho_1}\Big)^{\frac{n+s}{2}}\Big[\theta_2(2\rho_0)^{s+\epsilon}\Big(\fint_{\mathcal{B}(x_0,2\rho_0)}F^2d\mu\Big)^{\frac{1}{2}}
+\theta_2\epsilon^{\frac{1}{2}}(2\rho_0)^{\frac{3s}{2}+\epsilon}\Big(\int_{\mathbb{R}^n\setminus B_{2\rho_0}(x_0)}\frac{|f(y)|^2}{|y-x_0|^{n+s}}dy\Big)^{\frac{1}{2}}.\nonumber
\end{eqnarray}
These imply that
\begin{equation}
T(x,\rho)\leq c(n,s,q,\sigma)\Big(\frac{\rho_0}{\rho_2-\rho_1}\Big)^{\max\{\frac{n+2s}{2},\frac{n+q\sigma}{q},\frac{n+s}{2}\}}T(x_0,2\rho_0).
\end{equation}
Therefore, we obtain
\begin{eqnarray}
\frac{1}{\tau_0}(\Phi_{\theta_1,\theta_2}+\phi+T)(x,\rho)
&\leq&\frac{c(n,s,q,\sigma)}{\epsilon^{\frac{1}{\tilde{p}}}}\Big(\frac{\rho_0}{\rho_2-\rho_1}\Big)^{\max\{\frac{n+2\epsilon}{q},\frac{n+2\epsilon}{\tilde{p}},\frac{n+2s}{2},\frac{n+q\sigma}{q}\}}(\Phi_{\theta_1,\theta_2}+\phi+T)(x_0,2\rho_0).\nonumber
\end{eqnarray}
If $\mathcal{K}=K_1\times K_2\in\Xi_{k_0}$ with $k_0=\left[-\log_2\left(\frac{\rho_2-\rho_1}{n\cdot 10^{10n}}\right)\right]+1$, then $\mathcal{K}\subset\mathcal{B}(x_0,2\rho_0)$. Therefore, we have
\begin{eqnarray}
\mu(\mathcal{K})&=&\int_{K_1}\int_{K_2}\frac{1}{|x-y|^{n-2\epsilon}}dydx\nonumber\\
&\geq&\frac{1}{(2\rho_0)^{n-2\epsilon}}\int_{K_1}\int_{K_2}dydx\nonumber\\
&\geq&\frac{c(n)(\rho_2-\rho_1)^{2n}}{\rho^{n-2\epsilon}_0}.\nonumber
\end{eqnarray}

Consequently, we get
\begin{equation}
\Big(\fint_{\mathcal{K}}U^2d\mu\Big)^{\frac {1}{2}}\leq \left(\frac{\mu(\mathcal{B}(x_0,2\rho_0))}{\mu(\mathcal{K})}\fint_{\mathcal{B}(x_0,2\rho)}U^2d\mu\right)^{\frac {1}{2}}\leq \frac{c(n)}{\epsilon}\Big(\frac{\rho_0}{\rho_2-\rho_1}\Big)^{2n}\Big(\fint_{\mathcal{B}(x_0,2\rho_0)}U^2d\mu\Big)^{\frac {1}{2}},\nonumber
\end{equation}
Then by the definition of  $\tilde{t}_0$ and $t_1$, we obtain that
$$
\tilde{t}_0, t_1\leq
\frac{c(n,s,q,\sigma)}{\epsilon}\Big(\frac{\rho_0}{\rho_2-\rho_1}\Big)^{2n}(\Phi_{\theta_1,\theta_2}+\phi+T)(x_0,2\rho_0),
$$
which implies that
$$
\tilde{\tau}\tau_0 t_1\leq t_0=
\frac{c(n,s,q,\sigma)}{\epsilon}\Big(\frac{\rho_0}{\rho_2-\rho_1}\Big)^{2n}(\Phi_{\theta_1,\theta_2}+\phi+T)(x_0,2\rho_0).
$$
Consequently, we obtain
\begin{eqnarray}
\int_{\mathcal{B}_{\rho_1}\cap\{U>t\}}U^2d\mu
&\leq&C_1\epsilon^{-\frac{4s}{n+2\epsilon}}t^{2-\tilde{p}}\int_{\mathcal{B}_{\rho_2}\cap\{U> t\}}U^{\tilde{p}}d\mu
+C_2\Big[t^{2-q-\frac{q(1+\delta_1)[q(\sigma+\epsilon)+2\epsilon]}{n-q(\sigma+\epsilon)}}
t^{\frac{q(1+\delta_1)[q(\sigma+\epsilon)+2\epsilon]}{n-q(\sigma+\epsilon)}}_0\int_{\mathcal{B}_{\rho_2}\cap\{G>\tau_g t\}}G^qd\mu
\nonumber\\
& &
+t^{-\frac{2(s+\epsilon)(2+\delta_0)}{n-2s}}t^{\frac{2(s+\epsilon)(2+\delta_0)}{n-2s}}_0\int_{\mathcal{B}_{\rho_2}\cap\{F>\tau_f t\}}F^{2}d\mu\Big],\nonumber
\end{eqnarray}
for any $t\geq  t_0$, where $C_1=C_1(n,s,\lambda,\Lambda,\tilde{\beta},\beta,\delta_0,q)$ and
 $C_2=C_2(n,s,\lambda,\Lambda,\tilde{\beta},\beta,\delta_0,q,\epsilon)$.
Thus, Lemma 4.3 is proved.
\end{proof}


\begin{proof}[Proof of Theorem 4.2.]Assume $\rho_0=1$.
 Let $d\nu=U^2d\mu$ and let $U_m=\min\{U,m\}$ for any positive integer $m$. For any $\delta\in(0,1)$,
 \begin{eqnarray}\label{4.4.1}
 \int_{\mathcal{B}_{\rho_1}}U^\delta_mU^2d\mu&=&\int_{\mathcal{B}_{\rho_1}}U^\delta_md\nu\nonumber\\
 &=&\delta\int^\infty_0t^{\delta-1}\nu(\mathcal{B}_{\rho_1}\cap\{U_m>t\})dt\nonumber\\
 &=&\delta\int^\infty_0t^{\delta-1}\int_{\mathcal{B}_{\rho_1}\cap\{U_m>t\}}U^2d\mu dt\nonumber\\
  &=&\delta\int^m_0t^{\delta-1}\int_{\mathcal{B}_{\rho_1}\cap\{U>t\}}U^2d\mu dt\nonumber\\
 &=&\delta\int^{t_0}_0t^{\delta-1}\int_{\mathcal{B}_{\rho_1}\cap\{U>t\}}U^2d\mu dt
 +\delta\int^m_{t_0}t^{\delta-1}\int_{\mathcal{B}_{\rho_1}\cap\{U>t\}}U^2d\mu dt\nonumber\\
&\leq&t^{\delta}_0\int_{\mathcal{B}_{\rho_1}}U^2d\mu
 +\delta\int^m_{t_0}t^{\delta-1}\int_{\mathcal{B}_{\rho_1}\cap\{U>t\}}U^2d\mu dt.
 \end{eqnarray}
 For the first term on the right hand side of (\ref{4.4.1}), we have
\begin{eqnarray}\label{4.4.2}
t^{\delta}_0\int_{\mathcal{B}_{\rho_1}}U^2d\mu
\leq t^{\delta}_0\mu(\mathcal{B}_2)\fint_{\mathcal{B}_2}U^2d\mu
\leq C_0t^{2+\delta}_0\mu(\mathcal{B}_{\rho_1}),
\end{eqnarray}
where $C_0=C_0(n,s,\epsilon)$.
For the second term on the right hand side of (\ref{4.4.1}), by Lemma 4.3, we  can derive
\begin{eqnarray}\label{4.4.3}
&&\delta\int^m_{t_0}t^{\delta-1}\int_{\mathcal{B}_{\rho_1}\cap\{U>t\}}U^2d\mu dt\nonumber\\
&\leq&C_1\epsilon^{-\frac{4s}{n+2\epsilon}}\delta\int^m_{t_0}t^{\delta+1-\tilde{p}}\int_{\mathcal{B}_{\rho_2}\cap\{U> t\}}U^{\tilde{p}}d\mu dt
+C_2\Big[\delta \int^m_{t_0}t^{\delta+1-q}\int_{\mathcal{B}_{\rho_2}\cap\{G>\tau_g t\}}G^qd\mu dt\nonumber\\
& &
+\delta \int^m_{t_0}t^{\delta-1}\int_{\mathcal{B}_{\rho_2}\cap\{F>\tau_f t\}}F^{2}d\mu dt].
\end{eqnarray}
By a direct calculation,
\begin{eqnarray}\label{4.4.4}
C_1\epsilon^{-\frac{4s}{n+2\epsilon}}\delta\int^m_{t_0}t^{\delta+1-\tilde{p}}\int_{\mathcal{B}_{\rho_2}\cap\{U> t\}}U^{\tilde{p}}d\mu dt
&\leq&C_1\epsilon^{-\frac{4s}{n+2\epsilon}}\delta\int^\infty_{0}t^{\delta+1-\tilde{p}}\int_{\mathcal{B}_{\rho_2}\cap\{U_m> t\}}U^{\tilde{p}}d\mu dt\nonumber\\
&\leq&C_1\epsilon^{-\frac{4s}{n+2\epsilon}}\delta\int^{U_m}_{0}t^{\delta+1-\tilde{p}}dt\int_{\mathcal{B}_{\rho_2}}U^{\tilde{p}}d\mu\nonumber\\
&\leq&\frac{C_1\epsilon^{-\frac{4s}{n+2\epsilon}}\delta}{\delta+2-\tilde{p}}\int_{\mathcal{B}_{\rho_2}}U^{\tilde{p}}U^{\delta+2-\tilde{p}}_md\mu\nonumber\\
&\leq&\frac{C_1\epsilon^{-\frac{4s}{n+2\epsilon}}\delta}{\delta+2-\tilde{p}}\int_{\mathcal{B}_{\rho_2}}U^{2}U^{\delta}_md\mu.
\end{eqnarray}
By using H\"{o}lder's inequality,   for any $\delta\in(0, q(1+\delta_1)-2]$, and for $\tilde{G}=\frac {G}{\tau_g}$, we have
\begin{eqnarray}\label{4.4.5}
C_2\delta \int^m_{t_0}t^{\delta+1-q}\int_{\mathcal{B}_{\rho_2}\cap\{G>\tau_g t\}}G^qd\mu dt
&\leq&\tilde{C}_2\delta \int^m_{t_0}t^{\delta+1-q}\int_{\mathcal{B}_{\rho_2}\cap\{\tilde{G}>t\}}\tilde{G}^qd\mu dt\nonumber\\
&\leq&\tilde{C}_2 \delta \int^{\infty}_{0}t^{\delta+1-q}\int_{\mathcal{B}_{\rho_2}\cap\{\tilde{G}>t\}}\tilde{G}^qd\mu dt\nonumber\\
&\leq&\tilde{C}_2 \delta\int_{\mathcal{B}_{\rho_2}}\tilde{G}^{q}d\mu\int^{\tilde{G}}_{0}t^{\delta+1-q}dt\nonumber\\
&\leq&\frac{\tilde{C}_2 \delta}{\delta+2-q}\int_{\mathcal{B}_{\rho_2}}\tilde{G}^{2+\delta}d\mu\nonumber\\
&\leq&\frac{\hat{C}_2 \delta}{\delta+2-q}\int_{\mathcal{B}_{\rho_2}}{G}^{2+\delta}d\mu\nonumber\\
&\leq&\frac{\hat{C}_2 \delta\mu(\mathcal{B}_{\rho_2})}{\delta+2-q}\Big(\fint_{\mathcal{B}_{\rho_2}}G^{q(1+\delta_1)}d\mu\Big)^{\frac{2+\delta}{q(1+\delta_1)}}\nonumber\\
&\leq&\frac{\hat{\tilde{C}}_2 \delta\mu(\mathcal{B}_{\rho_1})}{\delta+2-q}t^{2+\delta}_0.
\end{eqnarray}
where $\hat{\tilde{C}}_2=\hat{\tilde{C}}_2(n,s,\lambda,\Lambda,\tilde{\beta},\beta,\delta_0,q,\epsilon)$.
Similarly,  for any $\delta\in(0, \delta_0]$, we have
\begin{eqnarray}\label{4.4.6}
C_2\delta \int^m_{t_0}t^{\delta-1}dt\int_{\mathcal{B}_{\rho_2}\cap\{F>\tau_f t\}}F^{2}d\mu
\leq \hat{C}_2\mu(\mathcal{B}_{\rho_1})t^{2+\delta}_0,
\end{eqnarray}
where $\hat{C}_2=\hat{C}_2(n,s,\lambda,\Lambda,\tilde{\beta},\beta,\delta_0,q,\epsilon)$.
Combining (\ref{4.4.1})-(\ref{4.4.6}), we obtain
\begin{equation}
\Big(\fint_{\mathcal{B}_{\rho_1}}U^\delta_mU^2d\mu\Big)^{\frac{1}{2+\delta}}
 \leq \Big(\frac{C_1\epsilon^{-\frac{4s}{n+2\epsilon}}\delta}{2+\delta-\tilde{p}}\Big)^{\frac{1}{2+\delta}}
\Big(\fint_{\mathcal{B}_{\rho_2}}U^{2}U^{\delta}_m d\mu\Big)^{\frac{1}{2+\delta}}
 +\tilde{C}t_0,\nonumber\\
\end{equation}
for any $\delta\in(0,\hat{\delta}]$, where $\hat{\delta}=\min\Big\{q(1+\delta_1)-2,\delta_0\Big\}$.
Take $\delta\in(0,\hat{\delta}]$ small such that $\delta<\frac{(2-\tilde{p})\epsilon^{\frac{4s}{n+2\epsilon}}}{8C_1}$,   then  we have
\begin{equation}
\Big(\fint_{\mathcal{B}_{\rho_1}}U^\delta_mU^2d\mu\Big)^{\frac{1}{2+\delta}}
 \leq \frac{1}{2}
\Big(\fint_{\mathcal{B}_{\rho_2}}U^{2}U^{\delta}_m d\mu\Big)^{\frac{1}{2+\delta}}
 +\tilde{C}t_0.\nonumber\\
\end{equation}
After an iteration process  and a limit argument as that in\cite{KMS2015}, we obtain
\begin{eqnarray}
\Big(\fint_{\mathcal{B}_{1}}U^{2+\delta}d\mu\Big)^{\frac{1}{2+\delta}}&\leq& \tilde{C}(\Phi_{\theta_1,\theta_2}+\phi+T)(\mathcal{B}_2).\nonumber
\end{eqnarray}
Using (\ref{4.001})-(\ref{4.003}) and  H\"{o}lder's inequality, we obtain that
\begin{eqnarray}
\Big(\fint_{\mathcal{B}_{1}}U^{2+\delta}d\mu\Big)^{\frac{1}{2+\delta}}&\leq&\tilde{\tilde{C}}\Big[\Big(\fint_{\mathcal{B}_2}U^2d\mu\Big)^{\frac{1}{2}}
+\Big(\fint_{\mathcal{B}_2}G^{q(1+\delta_1)}d\mu\Big)^{\frac{1}{q(1+\delta_1)}}
+\Big(\fint_{\mathcal{B}_2}F^{2+\delta_0}d\mu\Big)^{\frac{1}{2+\delta_0}}\nonumber\\
& &
+\Big(\int_{\mathbb{R}^n\setminus B_2}\frac{[u(y)-\bar{u}_{B_2}]^2}{|y|^{n+2s}}dy\Big)^{\frac{1}{2}}
+\Big(\int_{\mathbb{R}^n\setminus B_2}\frac{|g(y)-\bar{g}_{B_2}|^2}{|y|^{n+2s}}dy\Big)^{\frac{1}{2}}\nonumber\\
& &
+\Big(\int_{\mathbb{R}^n\setminus B_2}\frac{|f(y)|^2}{|y|^{n+s}}dy\Big)^{\frac{1}{2}}
\Big],\nonumber
\end{eqnarray}
where $\tilde{\tilde{\tilde{C}}}=\tilde{\tilde{\tilde{C}}}(n,s,\lambda,\Lambda,\tilde{\beta},\beta,\delta_0,q)$. Then, by a standard rescaling argument, we obtain the conclusion of Theorem 4.2.
\end{proof}

\begin{proof}[Proof of Theorem 4.1] The proof is a direct consequence of Theorem 4.2. Let $\mathcal{B}=B\times B$ with $B=B_1\subset\Omega$. Apply theorem 4.1 to $u-\bar{u}_{B}$, we have
\begin{eqnarray}
\Big(\fint_{\frac{1}{2}\mathcal{B}}U^{2+\delta}d\mu\Big)^{\frac{1}{2+\delta}}
&\leq&C\Big[\Big(\fint_{\mathcal{B}}U^2d\mu\Big)^{\frac{1}{2}}
+\Big(\fint_{\mathcal{B}}G^{q(1+\delta_1)}d\mu\Big)^{\frac{1}{q(1+\delta_1)}}
+\Big(\fint_{\mathcal{B}}F^{2+\delta_0}d\mu\Big)^{\frac{1}{2+\delta_0}}\nonumber\\
& &+\Big(\int_{\mathbb{R}^n\setminus B}\frac{[u(y)-\bar{u}_{B}]^2}{|y|^{n+2s}}dy\Big)^{\frac{1}{2}}
+\Big(\int_{\mathbb{R}^n\setminus B}\frac{|g(y)-\bar{g}_{B}|^2}{|y|^{n+2s}}dy\Big)^{\frac{1}{2}}\nonumber\\
& &
+\Big(\int_{\mathbb{R}^n\setminus B}\frac{|f(y)|^2}{|y|^{n+s}}dy\Big)^{\frac{1}{2}}
\Big],\nonumber
\end{eqnarray}
where $U$, $F$, $G$ are same to that in section 4,  $C=C(n,s,\lambda,\Lambda,\tilde{\beta},\beta,\delta_0,q)$.

Notice that by a direct calculation, we get
\begin{equation}
\Big(\fint_{\frac{1}{2}\mathcal{B}}U^{2+\delta}d\mu\Big)^{\frac{1}{2+\delta}}
=c(n,\epsilon)\Big(\frac{1}{2}\int_{\frac{1}{2}B}\int_{\frac{1}{2}B}\frac{|u(x)-u(y)|^{2+\delta}}{|x-y|^{n+(2+\delta)(s+\frac{\epsilon\delta}{2+\delta})}}dydx\Big)^{\frac{1}{2+\delta}},\nonumber
\end{equation}
\begin{equation}
\Big(\fint_{\mathcal{B}}U^{2}d\mu\Big)^{\frac{1}{2}}=c(n,\epsilon)\Big(\int_{B}\int_{B}\frac{|u(x)-u(y)|^{2}}{|x-y|^{n+2s}}dydx\Big)^{\frac{1}{2}}.\nonumber
\end{equation}
By the fractional Sobolev-Poincar\'{e} inequality, we have
\begin{eqnarray}
\Big(\fint_{\mathcal{B}}G^{q(1+\delta_1)}d\mu\Big)^{\frac{1}{q(1+\delta_1)}}
&\leq&c(n,\epsilon)\Big(\int_{B}\int_{B}\frac{|g(x)-g(y)|^{q(1+\delta_1)}}{|x-y|^{n+\sigma q(1+\delta_1)+2\epsilon\delta_1}}dydx\Big)^{\frac{1}{q(1+\delta_1)}}\nonumber\\
&\leq&c(n,\epsilon)\Big(\int_{B}\int_{B}\frac{|g(x)-g(y)|^{q(1+\delta_1)}}{|x-y|^{n+\sigma q(1+\delta_1)^2}}dydx\Big)^{\frac{1}{q(1+\delta_1)}}\nonumber
\end{eqnarray}
\begin{eqnarray}
\Big(\fint_{\mathcal{B}}F^{2+\delta_0}d\mu\Big)^{\frac{1}{2+\delta_0}}
&\leq&c(n,\epsilon)\Big(\int_{B}\int_{B}\frac{|f|^{2+\delta_0}(x)}{|x-y|^{n-2\epsilon}}dydx\Big)^{\frac{1}{2+\delta_0}}\nonumber\\
&\leq&c(n,\delta_0,\epsilon)\Big(\int_{B}|f(x)|^{2+\delta_0}dx\Big)^{\frac{1}{2+\delta_0}}.\nonumber
\end{eqnarray}
Then by the fact that $W^{s+\delta_0,2}(B)\hookrightarrow W^{\sigma(1+\delta_1),q(1+\delta_1)}(B)$, we deduce
\begin{eqnarray}
\Big(\int_{\frac{1}{2}B}\int_{\frac{1}{2}B}\frac{|u(x)-u(y)|^{2+\delta}}{|x-y|^{n+(s+\gamma)(2+\delta)}}dydx
\Big)^{\frac{1}{2+\delta}}
&\leq& C\Big[\Big(\int_{B}(\nabla^{s}_{B}u)^{2}dx\Big)^{\frac{1}{2}}+\Big(\int_{B}(\nabla^{s+\delta_0}_{B}g)^{2}dx\Big)^{\frac{1}{2}}\nonumber\\
& &+\Big(\int_{B}g^{2}dx\Big)^{\frac{1}{2}}+\Big(\int_{B}|f|^{2+\delta_0}dx\Big)^{\frac{1}{2+\delta_0}}\nonumber\\
& &+\Big(\int_{\mathbb{R}^n\setminus B}\frac{[u(y)-\bar{u}_{B}]^2}{|y|^{n+2s}}dy\Big)^{\frac{1}{2}}
+\Big(\int_{\mathbb{R}^n\setminus B}\frac{|g(y)-\bar{g}_{B}|^2}{|y|^{n+2s}}dy\Big)^{\frac{1}{2}}\nonumber\\
& &
+\Big(\int_{\mathbb{R}^n\setminus B}\frac{|f(y)|^2}{|y|^{n+s}}dy\Big)^{\frac{1}{2}},\nonumber
\end{eqnarray}
where $\gamma=\frac{\epsilon\delta}{2+\delta}$ and $C=C(n,s,\lambda,\Lambda,\tilde{\beta},\beta,\delta_0)$. This finishes the proof of Theorem 4.1 by a standard rescaling argument.
\end{proof}

\section{Interior $H^{s,p}$ estimate}\setcounter{equation}{0}
 In this section, we provide the proof of Theorem 1.1, i.e.,  the interior $H^{s,p}$ estimate  for the nonlocal elliptic equation(\ref{1.1})
with  $K$ and $b$ satisfying assumptions (\ref{1.2})-(\ref{1.3.1}).

\begin{lem}
Let  $K(x,y)$, $K_1(x,y)$  and $b(x,y)$ be kernels satisfying assumptions (\ref{1.2})-(\ref{1.3.1}) and
\begin{equation}\label{5.1.1}
|K_1(x,y)|\leq\frac{\mu}{|x-y|^{n+2s}}\hbox{ for any }x,y\in\mathbb{R}^n\hbox{ with }x\neq y,\\
\end{equation}
for some $\mu>0$.
 Then for any $h$, $g_1$ and $g_2$ satisfying  $\nabla^s_{\mathbb{R}^n}h,\nabla^s_{\mathbb{R}^n}g_1,g_2\in L^2(\mathbb{R}^n)$, there exists a unique solution $u\in H^s(\mathbb{R}^n)$ of  the Dirichlet problem
\begin{equation}\label{5.2}
\left\{
\begin{array}{lll}
 L^{2s}_Ku(x)=L^{2s}_{K_1}g_1((x)+L^s_bg_2(x),\quad &\hbox{ in }B_4,\\
 u=h,\quad &\hbox{ in }\mathbb{R}^n\setminus B_4.
\end{array}\right.
\end{equation}
Moreover, there holds
\begin{equation}
\|\nabla^s_{\mathbb{R}^n}u\|_{L^2(B_4)}\leq C\Big(\|\nabla^s_{\mathbb{R}^n}h\|_{L^2(B_4)}+\|\nabla^s_{\mathbb{R}^n}g_1\|_{L^2(B_4)}+\|g_2\|_{L^2(B_4)}+\int_{\mathbb{R}^n\setminus B_4}\frac{|g_2(x)|}{|x|^{n+s}}dx\Big),\nonumber\\
\end{equation}
where $C$ is a constant depending on $n$, $s$, $\lambda$, $\Lambda$, $\mu$ and $\beta$.
\end{lem}
\begin{proof} Let $\mathcal{H}=\mathcal{H}^s_0(B_4)$.
Define the bilinear form
$$
\varepsilon^{2s}_K(\omega,\phi)=\int_{\mathbb{R}^n}\int_{\mathbb{R}^n}K(x,y)[\omega(x)-\omega(y)][\phi(x)-\phi(y)]dxdy \quad\forall\omega, \phi\in \mathcal{H},
$$
and the functional
$$
F(\phi)=-\varepsilon^{2s}_K(h,\phi)+\varepsilon^{2s}_{K_1}(g_1,\phi)+\int_{\mathbb{R}^n}g_2(x)L^s_b\phi(x)dx\quad\forall\phi\in\mathcal{H}.
$$
Then the bilinear form $\varepsilon^{2s}_K$ is bounded and coercive on $\mathcal{H}$, and $F\in\mathcal{H}^*$, since
$$
\lambda\|\omega\|^2_{\mathcal{H}}\leq\varepsilon^{2s}_K(\omega,\omega)\quad\hbox{ and }\quad\varepsilon^{2s}_K(\omega,\phi)\leq\Lambda\|\omega\|_{\mathcal{H}}\|\phi\|_{\mathcal{H}}\quad\forall\omega,\phi\in\mathcal{H},
$$
and
\begin{eqnarray}
|F(\phi)|&\leq&\Lambda\|h\|_{\mathcal{H}}\|\phi\|_{\mathcal{H}}+\mu\|g_1\|_{\mathcal{H}}\|\phi\|_{\mathcal{H}}+\|g_2\|_{L^{2}(\mathbb{R}^n)}\|L^s_b\phi\|_{L^2(\mathbb{R}^n)}\nonumber\\
&\leq& C(n,s,\Lambda,\mu,\beta)[\|h\|_{\mathcal{H}} +\|g_1\|_{\mathcal{H}}+\|g_2\|_{L^2(\mathbb{R}^n)}]\|\phi\|_{\mathcal{H}}\quad\forall\phi\in\mathcal{H},\nonumber
\end{eqnarray}
due to assumptions (\ref{1.2})-(\ref{1.3.1}) and (\ref{5.1.1}). Then by the Lax-Milgram Theorem, there  exists a unique $\omega\in\mathcal{H}$ such that
$$
\varepsilon^{2s}_K(\omega,\phi)=F(\phi),\quad\forall\phi\in\mathcal{H}.
$$
Let $u=\omega+h$. Then we get
$$
\varepsilon^{2s}_K(u,\phi)=\varepsilon^{2s}_{K_1}(g_1,\phi)+\int_{\mathbb{R}^n}g_2(x)L^s_b\phi(x)dx,\quad\forall \phi\in \mathcal{H},
$$
which means that $u$ is the unique solution of (\ref{5.2}). Moreover, by H\"{o}lder's inequality, we have
\begin{eqnarray}
&&\|\nabla^s_{\mathbb{R}^n}\omega\|^2_{L^2(B_4)}\nonumber\\
&=&\int_{B_4}\int_{\mathbb{R}^n}\frac{[\omega(x)-\omega(y)]^2}{|x-y|^{n+2s}}dydx\nonumber\\
&\leq&\int_{\mathbb{R}^n}\int_{\mathbb{R}^n}\frac{[\omega(x)-\omega(y)]^2}{|x-y|^{n+2s}}dydx\nonumber\\
&\leq&\frac{1}{\lambda}\varepsilon^{2s}_K(\omega,\omega)\nonumber\\
&=&\frac{1}{\lambda}\Big[-\varepsilon^{2s}_K(h,\omega)+\varepsilon^{2s}_{K_1}(g_1,\omega)+\int_{\mathbb{R}^n}g_2(x)L^s_b\omega(x)dx\Big]\nonumber\\
&\leq&\frac{2\Lambda}{\lambda}\int_{B_4}\int_{\mathbb{R}^n}\frac{|h(x)-h(y)||\omega(x)-\omega(y)|}{|x-y|^{n+2s}}dydx\nonumber\\
& &+\frac{2\mu}{\lambda}\int_{B_4}\int_{\mathbb{R}^n}\frac{|g_1(x)-g_1(y)||\omega(x)-\omega(y)|}{|x-y|^{n+2s}}dydx
+\frac{1}{\lambda}\int_{\mathbb{R}^n}g_2(x)L^s_b\omega(x)dx\nonumber\\
&\leq&\frac{2\Lambda}{\lambda}\|\nabla^s_{\mathbb{R}^n}h\|_{L^2(B_4)}\|\nabla^s_{\mathbb{R}^n}\omega\|_{L^2(B_4)}
+\frac{2\mu}{\lambda}\|\nabla^s_{\mathbb{R}^n}g_1\|_{L^2(B_4)}\|\nabla^s_{\mathbb{R}^n}\omega\|_{L^2(B_4)}\nonumber\\
& &
+\frac{1}{\lambda}c(n,s,\beta)\Big(\int_{\mathbb{R}^n\setminus B_8}\frac{|g_2(x)|}{|x|^{n+s}}dx+\|g_2\|_{L^2(B_8)}\Big)\|\nabla^s_{\mathbb{R}^n}\omega\|_{L^2(B_4)}\nonumber\\
&\leq&C(n,s,\lambda,\Lambda,\mu,\beta)\Big(\|\nabla^s_{\mathbb{R}^n}h\|_{L^2(B_4)}
+\|\nabla^s_{\mathbb{R}^n}g_1\|_{L^2(B_4)}+\int_{\mathbb{R}^n\setminus B_8}\frac{|g_2(x)|}{|x|^{n+s}}dx+\|g_2\|_{L^2(B_8)}\Big)\|\nabla^s_{\mathbb{R}^n}\omega\|_{L^2(B_4)}.\nonumber
\end{eqnarray}
where we used the fact that
\begin{eqnarray}
\int_{\mathbb{R}^n}g_2(x)L^s_b\omega(x)dx&=&\int_{\mathbb{R}^n}g_2(x)\int_{\mathbb{R}^n}b(x,y)[\omega(x)-\omega(y)]dydx\nonumber\\
&=&\int_{B_8}g_2(x)\int_{\mathbb{R}^n}b(x,y)[\omega(x)-\omega(y)]dydx-\int_{\mathbb{R}^n\setminus B_8}g_2(x)\int_{B_4}b(x,y)\omega(y)dydx\nonumber\\
&=&\int_{\mathbb{R}^n}\hat{g}_2(x)\int_{\mathbb{R}^n}b(x,y)[\omega(x)-\omega(y)]dydx-\int_{\mathbb{R}^n\setminus B_8}g_2(x)\int_{B_4}b(x,y)\omega(y)dydx\nonumber\\
&\leq&\|\hat{g}_2\|_{L^2(\mathbb{R}^n)}\|L^s_b\omega\|_{L^2(\mathbb{R}^n)}+\beta\int_{\mathbb{R}^n\setminus B_8}\int_{B_4}\frac{|g_2(x)|}{|x-y|^{n+s}}|\omega(y)|dydx\nonumber\\
&\leq&c(n,s,\beta)\|{g}_2\|_{L^2(B_8)}\|\nabla^s_{\mathbb{R}^n}\omega\|_{L^2(\mathbb{R}^n)}+c(n,s,\beta)\int_{\mathbb{R}^n\setminus B_8}\frac{|g_2(x)|}{|x|^{n+s}}dx\Big(\int_{B_4}|\omega(y)|^2dy\Big)^{\frac{1}{2}}\nonumber\\
&\leq&c(n,s,\beta)\|{g}_2\|_{L^2(B_8)}\|\nabla^s_{\mathbb{R}^n}\omega\|_{L^2(\mathbb{R}^n)}+c(n,s,\beta)\int_{\mathbb{R}^n\setminus B_8}\frac{|g_2(x)|}{|x|^{n+s}}dx\|\nabla^s_{\mathbb{R}^n}\omega\|_{L^2(B_4)}\nonumber
\end{eqnarray}
where
\begin{equation}
\hat{g_2}(x)=\left\{
\begin{array}{lll}
g_2(x),\quad &\hbox{ in }B_8,\nonumber\\
0,\quad & \hbox{ in }\mathbb{R}^n\setminus B_8.\nonumber\\
\end{array}\right.
\end{equation}
Therefore, we get
\begin{eqnarray}
\|\nabla^s_{\mathbb{R}^n}\omega\|_{L^2(B_4)}&\leq& C(n,s,\lambda,\Lambda,\mu,\beta)\Big(\|\nabla^s_{\mathbb{R}^n}h\|_{L^2(B_4)}
+\|\nabla^s_{\mathbb{R}^n}g_1\|_{L^2(B_4)}+\|g_2\|_{L^2(B_8)}+\int_{\mathbb{R}^n\setminus B_8}\frac{|g_2(x)|}{|x|^{n+s}}dx\Big)\nonumber\\
&\leq&\tilde{C}(n,s,\lambda,\Lambda,\mu,\beta)\Big(\|\nabla^s_{\mathbb{R}^n}h\|_{L^2(B_4)}
+\|\nabla^s_{\mathbb{R}^n}g_1\|_{L^2(B_4)}+\|g_2\|_{L^2(B_4)}+\int_{\mathbb{R}^n\setminus B_4}\frac{|g_2(x)|}{|x|^{n+s}}dx\Big).\nonumber
\end{eqnarray}
This completes the proof of  Lemma 5.1.
\end{proof}

\begin{lem}\label {Lemma-1}
Let $0<s<1$ and $\Omega$ be a bounded domain in $\mathbb{R}^n$. Let  the kernels $K$, $K_1$ and $b$ satisfy assumptions (\ref{1.2})-(\ref{1.3.1}) and (\ref{5.1.1}).
Let $u\in {H}^{s}(\mathbb{R}^n)$ be a weak solution of
\begin{equation}
 L^{2s}_Ku(x)=L^{2s}_{K_1}g(x)+L^s_bf(x),\quad \hbox{ in }\Omega\subset\mathbb{R}^n,\nonumber\\
\end{equation}
 with $g\in H^{s}(\mathbb{R}^n)$ and  $f\in L^2(\mathbb{R}^n)$.  Then for any $B=B_r(x_0)\subset\Omega$, we have
 \begin{eqnarray}
\int_{\frac{1}{2}B}\Big(\nabla^s_{\mathbb{R}^n}u\Big)^2dx
&\leq&C\Big[\frac{1}{r^{2s}}\int_{B}u^2(x)dx+r^n\int_{\mathbb{R}^n\setminus{B}}\frac{u^2(x)}{|x_0-x|^{n+2s}}dx
+\int_{B}\int_{B}\frac{|g(x)-g(y)|^2}{|x-y|^{n+2s}}dydx\nonumber\\
& &+r^n\int_{\mathbb{R}^n\setminus B}\frac{|g(x)-\bar{g}_B|^2}{|x_0-x|^{n+2s}}dx
+\int_Bf^2(x)dx+r^{n+s}\int_{\mathbb{R}^n\setminus{B}}\frac{f^2(x)}{|x_0-x|^{n+s}}dx\Big],\nonumber
\end{eqnarray}
where $C=C(n,s,\lambda,\Lambda,\mu,\beta)$.
\end{lem}

\begin{proof} The proof of Lemma \ref{Lemma-1} is very similar to that of Lemma \ref{lemma-2}. For the convenience of the readers, we present it in detail here.

Let $\psi\in C^\infty_0(\frac{3}{4}B)$ be a cut-off function
satisfying
\begin{equation}
0\leq\psi\leq1, \quad \psi\equiv1\hbox{ in }\frac{1}{2}B, \quad\psi\equiv0\hbox{ in }\mathbb{R}^n\setminus\frac{3}{4}B\quad\hbox{ and }\quad|D\psi|\leq\frac{c(n)}{r}.\nonumber\\
\end{equation}
Take $\phi=\psi^2u\in \mathcal{H}^{s}_0(\Omega)$, we have
\begin{eqnarray}\label{5.2.1}
&&\int_{\mathbb{R}^n}\int_{\mathbb{R}^n}K(x,y)[u(x)-u(y)][\psi^2(x)u(x)-\psi^2(y)u(y)]dydx\nonumber\\
&=&\int_{\mathbb{R}^n}\int_{\mathbb{R}^n}K_1(x,y)[g(x)-g(y)][\psi^2(x)u(x)-\psi^2(y)u(y)]dydx+\int_{\mathbb{R}^n}f(x)L^s_b(\psi^2u)(x)dx.
\end{eqnarray}
Similar to Lemma 4.1, we have
\begin{eqnarray}
&&\int_{\mathbb{R}^n\times\mathbb{R}^n}K(x,y)[u(x)-u(y)][\psi^2(x)u(x)-\psi^2(y)u(y)]dydx\nonumber\\
&=&\int_{(\mathbb{R}^n\times\mathbb{R}^n)\cap\{\psi(x)>\psi(y)\}}K(x,y)[u(x)-u(y)][\psi^2(x)u(x)-\psi^2(y)u(y)]dydx\nonumber\\
& &+\int_{(\mathbb{R}^n\times\mathbb{R}^n)\cap\{\psi(x)\leq\psi(y)\}}K(x,y)[u(x)-u(y)][\psi^2(x)u(x)-\psi^2(y)u(y)]dydx\nonumber\\
&\geq&\int_{(\mathbb{R}^n\times\mathbb{R}^n)\cap\{\psi(x)>\psi(y)\}}\frac{\lambda}{2}\frac{[u(x)-u(y)]^2\psi^2(x)}{|x-y|^{n+2s}}
-2\Lambda\frac{[\psi(x)-\psi(y)]^2u^2(y)}{|x-y|^{n+2s}}dydx\nonumber\\
& &+\int_{(\mathbb{R}^n\times\mathbb{R}^n)\cap\{\psi(x)\leq\psi(y)\}}\frac{\lambda}{2}\frac{[u(x)-u(y)]^2\psi^2(y)}{|x-y|^{n+2s}}
-2\Lambda\frac{[\psi(x)-\psi(y)]^2u^2(x)}{|x-y|^{n+2s}}dydx\nonumber\\
&\geq&\frac{\lambda}{2}\int_{\mathbb{R}^n\times\mathbb{R}^n}\frac{[u(x)-u(y)]^2\psi^2(x)}{|x-y|^{n+2s}}dydx
-2\Lambda\int_{\mathbb{R}^n\times\mathbb{R}^n}\frac{[\psi(x)-\psi(y)]^2u^2(x)}{|x-y|^{n+2s}}dydx.\nonumber
\end{eqnarray}
In addition, note that
\begin{eqnarray}
&&\int_{\mathbb{R}^n\times\mathbb{R}^n}\frac{[\psi(x)-\psi(y)]^2u^2(x)}{|x-y|^{n+2s}}dydx\nonumber\\
&=&\int_{B}\int_{B}\frac{[\psi(x)-\psi(y)]^2 u(x)}{|x-y|^{n+2s}}dydx
+\int_{\mathbb{R}^n\setminus B}\int_{B}\frac{[\psi(x)-\psi(y)]^2 u^2(x)}{|x-y|^{n+2s}}dydx\nonumber\\
& &
+\int_{B}\int_{\mathbb{R}^n\setminus B}\frac{[\psi(x)-\psi(y)]^2 u^2(x)}{|x-y|^{n+2s}}dydx\nonumber\\
&\leq&\frac{c(n)}{r^2}\int_{B} u^2(x)dx\int_B\frac{1}{|x-y|^{n+2s-2}}dy+\int_{\mathbb{R}^n\setminus B} u^2(x)dx\int_{\frac{3}{4}B}\frac{1}{|x-y|^{n+2s}}dy\nonumber\\
& &+\int_{\frac{3}{4}B} u^2(x)dx\int_{\mathbb{R}^n\setminus B}\frac{1}{|x-y|^{n+2s}}dy\nonumber\\
&\leq&c(n,s)\Big(\frac{1}{r^{2s}}\int_{B} u^2(x)dx\
+r^n\int_{\mathbb{R}^n\setminus B} \frac{u^2(x)}{|x-x_0|^{n+2s}}dx\Big).\nonumber
\end{eqnarray}
Consequently, we get
\begin{eqnarray}\label{5.2.2}
&&\int_{\mathbb{R}^n}\int_{\mathbb{R}^n}K(x,y)[(u(x)-u(y)][\psi^2(x)u(x)-\psi^2(y)u(y)]dydx\nonumber\\
&\geq&\frac{\lambda}{2}\int_{\mathbb{R}^n}\int_{\mathbb{R}^n}\frac{[u(x)-u(y)]^2\psi^2(x)}{|x-y|^{n+2s}}dydx\nonumber\\
& &-c(n,s,\Lambda)\Big(\frac{1}{r^{2s}}\int_{B} u^2(x)dx\
+r^n\int_{\mathbb{R}^n\setminus B} \frac{u^2(x)}{|x-x_0|^{n+2s}}dx\Big).
\end{eqnarray}
The estimate for the first term on the right-hand side of (\ref{5.2.1}) is similar to that of  \cite{KMS2015}. We rewrite it as
\begin{eqnarray}
&&\int_{\mathbb{R}^n}\int_{\mathbb{R}^n}K_1(x,y)[g(x)-g(y)][\psi^2(x)u(x)-\psi^2(y)u(y)]dydx\nonumber\\
&=&\int_{B}\int_{B}K_1(x,y)[g(x)-g(y)][\psi^2(x)u(x)-\psi^2(y)u(y)]dydx\nonumber\\
& &+\int_{B}\int_{\mathbb{R}^n\setminus B}K_1(x,y)[g(x)-g(y)]\psi^2(x)u(x)dydx\nonumber\\
& &-\int_{\mathbb{R}^n\setminus B}\int_{B}K_1(x,y)[g(x)-g(y)]\psi^2(y)u(y)dydx\nonumber\\
& &=I+II+III.\nonumber
\end{eqnarray}
For $I$, we have
\begin{eqnarray}
I&=&\int_{B}\int_{B}K_1(x,y)[g(x)-g(y)][\psi^2(x)u(x)-\psi^2(y)u(y)]dydx\nonumber\\
&=&\int_{B}\int_{B}K_1(x,y)[g(x)-g(y)][\psi(x)u(x)-\psi(y)u(y)]\psi(y)dydx\nonumber\\
& &+\int_{B}\int_{B}K_1(x,y)[g(x)-g(y)]\psi(x)u(x)[\psi(x)-\psi(y)]dydx\nonumber\\
&\leq&\mu\int_{B}\int_{B}\frac{|g(x)-g(y)||\psi(x)u(x)-\psi(y)u(y)|}{|x-y|^{n+2s}}dydx\nonumber\\
& &+\mu\int_{B}\int_{B}\frac{|g(x)-g(y)|\psi(x)|u(x)||\psi(x)-\psi(y)|}{|x-y|^{n+2s}}dydx\nonumber\\
&=&I_1+I_2.\nonumber
\end{eqnarray}
By H\"{o}lder's inequality and Cauchy's inequality, we get
\begin{eqnarray}
I_1&=&\mu\int_{B}\int_{B}\frac{|g(x)-g(y)||\psi(x)u(x)-\psi(y)u(y)|}{|x-y|^{n+2s}}dydx\nonumber\\
&\leq&\mu\Big(\int_{B}\int_{B}\frac{|g(x)-g(y)|^2}{|x-y|^{n+2s}}dydx\Big)^{\frac{1}{2}}\Big(\int_{B}\int_{B}\frac{|\psi(x)u(x)-\psi(y)u(y)|^2}{|x-y|^{n+2s}}dydx\Big)^{\frac{1}{2}}\nonumber\\
&\leq&\epsilon\int_{B}\int_{B}\frac{|\psi(x)u(x)-\psi(y)u(y)|^2}{|x-y|^{n+2s}}dydx
+\frac{\mu^2}{\epsilon}\int_{B}\int_{B}\frac{|g(x)-g(y)|^2}{|x-y|^{n+2s}}dydx\nonumber\\
&\leq&2\epsilon\int_{B}\int_{B}\frac{\psi^2(x)|u(x)-u(y)|^2}{|x-y|^{n+2s}}dydx+\frac{2\epsilon c(n,s)}{r^{2s}}\int_Bu^2(x)dx
+\frac{\mu^2}{\epsilon}\int_{B}\int_{B}\frac{|g(x)-g(y)|^2}{|x-y|^{n+2s}}dydx,\nonumber
\end{eqnarray}
for any $\epsilon>0$,  where we used the facts that
\begin{equation}
|\psi(x)u(x)-\psi(y)u(y)|^2\leq 2\psi^2(x)[u(x)-u(y)]^2+2[\psi(x)-\psi(y)]^2u^2(y),\nonumber\\
\end{equation}
and
\begin{equation}
\int_{B}\int_{B}\frac{|\psi(x)-\psi(y)|^2u(y)^2}{|x-y|^{n+2s}}dydx
\leq\frac{c(n)}{r^2}\int_{B}\int_{B}\frac{u(y)^2}{|x-y|^{n+2s-2}}dydx
\leq\frac{ c(n,s)}{r^{2s}}\int_Bu^2(x)dx.\nonumber
\end{equation}
By H\"{o}lder's inequality, fractional Poincar\'{e} inequality and Cauchy's inequality, we have
\begin{eqnarray}
I_2&=&\mu\int_{B}\int_{B}\frac{|g(x)-g(y)|\psi(x)|u(x)||\psi(x)-\psi(y)|}{|x-y|^{n+2s}}dydx\nonumber\\
&\leq&\frac{\mu c(n)}{r}\int_{B}\int_{B}\frac{|g(x)-g(y)|\psi(x)|u(x)|}{|x-y|^{n+2s-1}}dydx\nonumber\\
&\leq&\frac{\mu c(n)}{r}\Big(\int_{B}\int_{B}\frac{|g(x)-g(y)|^2}{|x-y|^{n+3s-1}}dydx\Big)^{\frac{1}{2}}\Big(\int_{B}\int_{B}\frac{\psi^2(x)u^2(x)}{|x-y|^{n+s-1}}dydx\Big)^{\frac{1}{2}}\nonumber\\
&\leq&\frac{\mu c(n,s)}{r^s}\Big(\int_{B}\int_{B}\frac{|g(x)-g(y)|^2}{|x-y|^{n+2s}}dydx\Big)^{\frac{1}{2}}\Big(\int_{B}\psi^2(x)u^2(x)dx\Big)^{\frac{1}{2}}\nonumber\\
&\leq&\mu c_1(n,s)\Big(\int_{B}\int_{B}\frac{|g(x)-g(y)|^2}{|x-y|^{n+2s}}dydx\Big)^{\frac{1}{2}}\Big(\int_B\int_{B}\frac{[\psi(x)u(x)-\psi(y)u(y)]^2}{|x-y|^{n+2s}}dydx\Big)^{\frac{1}{2}}\nonumber\\
&\leq&\epsilon\int_B\int_{B}\frac{[\psi(x)u(x)-\psi(y)u(y)]^2}{|x-y|^{n+2s}}dydx+\frac{\mu^2 c^2_1(n,s)}{\epsilon}\int_{B}\int_{B}\frac{|g(x)-g(y)|^2}{|x-y|^{n+2s}}dydx\nonumber\\
&\leq&2\epsilon\int_B\int_{B}\frac{\psi^2(x)[u(x)-u(y)]^2}{|x-y|^{n+2s}}dydx+\frac{2\epsilon c(n,s)}{r^{2s}}\int_Bu^2(x)dx
+\frac{\mu^2 c^2_1(n,s)}{\epsilon}\int_{B}\int_{B}\frac{|g(x)-g(y)|^2}{|x-y|^{n+2s}}dydx.\nonumber
\end{eqnarray}
Therefore, we obtain
\begin{eqnarray}
I&\leq&4\epsilon\int_B\int_{B}\frac{\psi^2(x)[u(x)-u(y)]^2}{|x-y|^{n+2s}}dydx+\frac{4\epsilon C(n,s)}{r^{2s}}\int_Bu^2(x)dx
+\frac{ C(n,s,\mu)}{\epsilon}\int_{B}\int_{B}\frac{|g(x)-g(y)|^2}{|x-y|^{n+2s}}dydx,\nonumber
\end{eqnarray}
for any $\epsilon>0$. We rewrite $II$ as
\begin{eqnarray}
II&=&\int_{B}\int_{\mathbb{R}^n\setminus B}K_1(x,y)[g(x)-g(y)]\psi^2(x)u(x)dydx\nonumber\\
&\leq&\mu\int_{B}\int_{\mathbb{R}^n\setminus B}\frac{|g(x)-\bar{g}_B|\psi(x)|u(x)|}{|x-y|^{n+2s}}dydx
+\mu\int_{B}\int_{\mathbb{R}^n\setminus B}\frac{|g(y)-\bar{g}_B|\psi(x)|u(x)|}{|x-y|^{n+2s}}dydx\nonumber\\
&=&II_1+II_2.\nonumber
\end{eqnarray}
By H\"{o}lder's inequality, fractional Poincar\'{e} inequality and Cauchy's inequality, we have
\begin{eqnarray}\label{g21}
II_1&=&\mu\int_{B}\int_{\mathbb{R}^n\setminus B}\frac{|g(x)-\bar{g}_B|\psi(x)|u(x)|}{|x-y|^{n+2s}}dydx\nonumber\\
&\leq&\mu\sup_{x\in supp\psi}\Big(\int_{\mathbb{R}^n\setminus B}\frac{1}{|x-y|^{n+2s}}dy\Big)\Big(\int_B\psi^2(x)u^2(x)dx\Big)^{\frac{1}{2}}\Big(\int_B[g(x)-\bar{g}_B]^2dx\Big)^{\frac{1}{2}}\nonumber\\
&\leq&c(n,s,\mu)\Big(\int_B\int_{B}\frac{[\psi(x)u(x)-\psi(y)u(y)]^2}{|x-y|^{n+2s}}dydx\Big)^{\frac{1}{2}}\Big(\int_{B}\int_{B}\frac{|g(x)-g(y)|^2}{|x-y|^{n+2s}}dydx\Big)^{\frac{1}{2}}\nonumber\\
&\leq&\epsilon\int_B\int_{B}\frac{[\psi(x)u(x)-\psi(y)u(y)]^2}{|x-y|^{n+2s}}dydx+\frac{ c^2(n,s,\mu)}{\epsilon}\int_{B}\int_{B}\frac{|g(x)-g(y)|^2}{|x-y|^{n+2s}}dydx\nonumber\\
&\leq&2\epsilon\int_{B}\int_{B}\frac{\psi^2 (x)[u(x)-u(y)]^2}{|x-y|^{n+2s}}dydx+\frac{2\epsilon c(n,s)}{r^{2s}}\int_Bu^2(x)dx
+\frac{c^2(n,s,\mu)}{\epsilon}\int_{B}\int_{B}\frac{|g(x)-g(y)|^2}{|x-y|^{n+2s}}dydx.\nonumber
\end{eqnarray}
By H\"{o}lder's inequality, fractional Poincar\'{e} inequality and Cauchy's inequality,
\begin{eqnarray}\label{g22}
II_2&=&\mu\int_{B}\int_{\mathbb{R}^n\setminus B}\frac{|g(y)-\bar{g}_B|\psi(x)|u(x)|}{|x-y|^{n+2s}}dydx\nonumber\\
&\leq&c(n,s,\mu)\int_{B}\psi(x)|u(x)|dx\int_{\mathbb{R}^n\setminus B}\frac{|g(y)-\bar{g}_B|}{|x_0-y|^{n+2s}}dy\nonumber\\
&\leq&c_1(n,s,\mu)r^{\frac{n}{2}-s}\Big(\int_{B}\psi^2(x)u^2(x)dx\Big)^{\frac{1}{2}}\Big(\int_{\mathbb{R}^n\setminus B}\frac{|g(y)-\bar{g}_B|^2}{|x_0-y|^{n+2s}}dy\Big)^{\frac{1}{2}}\nonumber\\
&\leq&c_2(n,s,\mu)r^{\frac{n}{2}}\Big(\int_{B}\int_{B}\frac{[\psi (x)u(x)-\psi(y)u(y)]^2}{|x-y|^{n+2s}}dydx\Big)^{\frac{1}{2}}\Big(\int_{\mathbb{R}^n\setminus B}\frac{|g(y)-\bar{g}_B|^2}{|x_0-y|^{n+2s}}dy\Big)^{\frac{1}{2}}\nonumber\\
&\leq&\epsilon\int_{B}\int_{B}\frac{[\psi (x)u(x)-\psi(y)u(y)]^2}{|x-y|^{n+2s}}dydx
+\frac{c^2_2(n,s,\mu)}{\epsilon}r^n\Big(\int_{\mathbb{R}^n\setminus B}\frac{|g(y)-\bar{g}_B|^2}{|x_0-y|^{n+2s}}dy\Big)^{\frac{1}{2}}\nonumber\\
&\leq&2\epsilon\int_{B}\int_{B}\frac{\psi^2 (x)[u(x)-u(y)]^2}{|x-y|^{n+2s}}dydx+\frac{2\epsilon c(n,s)}{r^{2s}}\int_Bu^2(x)dx
+\frac{c^2_2(n,s,\mu)}{\epsilon}r^{n}\int_{\mathbb{R}^n\setminus B}\frac{|g(y)-\bar{g}_B|^2}{|x_0-y|^{n+2s}}dy.\nonumber
\end{eqnarray}
Therefore, we get
\begin{eqnarray}
II&\leq&4\epsilon\int_{B}\int_{B}\frac{\psi^2 (x)[u(x)-u(y)]^2}{|x-y|^{n+2s}}dydx+\frac{4\epsilon c(n,s)}{r^{2s}}\int_Bu^2(x)dx\nonumber\\
& &+\frac{C_1(n,s,\mu)}{\epsilon}\int_{B}\int_{B}\frac{|g(x)-g(y)|^2}{|x-y|^{n+2s}}dydx
+\frac{C_2(n,s,\mu)}{\epsilon}r^{n}\int_{\mathbb{R}^n\setminus B}\frac{|g(y)-\bar{g}_B|^2}{|x_0-y|^{n+2s}}dy,\nonumber
\end{eqnarray}
for any $\epsilon>0$.
The estimate for III is the same to that for II. Combining  the above inequalities, we obtain
\begin{eqnarray}\label{5.2.5}
&&\int_{\mathbb{R}^n}\int_{\mathbb{R}^n}K_1(x,y)[g(x)-g(y)][\psi^2(x)u(x)-\psi^2(y)u(y)]dydx\nonumber\\
&\leq&12\epsilon\int_{B}\int_{B}\frac{\psi^2 (x)[u(x)-u(y)]^2}{|x-y|^{n+2s}}dydx+\frac{12\epsilon c(n,s)}{r^{2s}}\int_Bu^2(x)dx\nonumber\\
& &+\frac{C_1(n,s,\mu)}{\epsilon}\int_{B}\int_{B}\frac{|g(x)-g(y)|^2}{|x-y|^{n+2s}}dydx
+\frac{C_2(n,s,\mu)}{\epsilon}r^{n}\int_{\mathbb{R}^n\setminus B}\frac{|g(y)-\bar{g}_B|^2}{|x_0-y|^{n+2s}}dy,
\end{eqnarray}
 for any $\epsilon>0$.
Similar to (\ref{4.+13}), we have
\begin{eqnarray}\label{5.2.6}
\int_{\mathbb{R}^n}f(x)L^s_b(\psi^2u)(x)dx
&\leq&2c(n,s,\beta)\epsilon_1\int_{\mathbb{R}^n}\int_{\mathbb{R}^n}\psi^2(x)\frac{[u(x)-u(y)]^2}{|x-y|^{n+2s}}dydx
+\frac{1}{\epsilon_1}\int_{B}f^2(x)dx\nonumber\\
& &
+c_1(n,s,\beta)\Big[\frac{1}{r^{2s}}\int_{B}u^2(x)dy
+r^{n+s}\int_{\mathbb{R}^n\setminus B}\frac{f^2(x)}{|x-x_0|^{n+s}}dx\Big],
\end{eqnarray}
for any $\epsilon_1>0$.
Combining (\ref{5.2.1})-(\ref{5.2.6})  with $\epsilon=\frac{\lambda}{96}$ and $\epsilon_1=\frac{\lambda}{16c(n,s,\beta)}$, we obtain
\begin{eqnarray}
\int_{\frac{1}{2}B}\int_{\mathbb{R}^n}\frac{[u(x)-u(y)]^2}{|x-y|^{n+2s}}dydx
&\leq&\int_{\mathbb{R}^n}\int_{\mathbb{R}^n}\psi^2(x)\frac{[u(x)-u(y)]^2}{|x-y|^{n+2s}}dydx\nonumber\\
&\leq&C\Big[\frac{1}{r^{2s}}\int_{B}u^2(x)dx+r^{n}\int_{\mathbb{R}^n\setminus B} \frac{u^2(x)}{|x-x_0|^{n+2s}}dx+\int_{B}\int_{B}\frac{|g(x)-g(y)|^2}{|x-y|^{n+2s}}dydx\nonumber\\
& &+r^{n}\int_{\mathbb{R}^n\setminus B}\frac{|g(y)-\bar{g}_B|^2}{|x_0-y|^{n+2s}}dy
+\int_Bf^2(x)dx+r^{n+s}\int_{\mathbb{R}^n\setminus B}\frac{f^2(x)}{|x-x_0|^{n+s}}dx\Big],\nonumber
\end{eqnarray}
where $C=C(n,s,\lambda,\Lambda,\mu,\beta)$. Thus, Lemma 5.2 is proved.
\end{proof}

Using Lemmas 5.1 and 5.2, together with arguments from \cite{BR2013}, we obtain the following key lemma.

\begin{lem}For any $p>2$, there exists  $N_0>0$ such that for any $\epsilon>0$, there is $\delta=\delta(n,s,\lambda, \Lambda,\beta, p,\epsilon)>0$
such that if  $u\in {H}^s(\mathbb{R}^n)$  is  a weak solution of (\ref{1.1})
in $\Omega\supset B_8$ with $A$ $(\delta,8)$-vanishing and $f\in L^2(\mathbb{R}^n)\cap L^p(B_8)$ such that if there exists $x_0\in B_1$ satisfying
\begin{equation}\label{5.3.1}
\nabla^s_{\mathbb{R}^n}u(x_0)\leq1,\quad\mathcal{M}_{B_8}(|f|^{p_0})(x_0)\leq \delta^{p_0}\hbox{ and }
\int_{\mathbb{R}^n\setminus B_5(x_0)}\frac{|f|^{2}}{|x-x_0|^{n+s}}dx\leq \delta^{2}\quad\hbox{ for some }p_0\in(2,p),
\end{equation}
then there exists a weak solution $v\in\mathcal{H}^s_0(B_4)$ of
\begin{eqnarray}
L^{2s}_{K_0}v=g,\hbox{ in }B_4,\nonumber
\end{eqnarray}
where
$$
K_0=\frac{{A}_{0}(x,y)}{|x-y|^{n+2s}}
$$
with
\begin{equation}
A_0(x,y)=\left\{
\begin{array}{lll}
\bar{A}_{B_3\times B_3},\quad &\hbox{ if }(x,y)\in B_3\times B_3,\nonumber\\
A(x,y),\quad & \hbox{ if }(x,y)\in(\mathbb{R}^n\times\mathbb{R}^n)\setminus (B_3\times B_3),
\end{array}\right.
\end{equation}
and
\begin{equation}
 g(x)=L^s_bf\chi_{\mathbb{R}^n\setminus B_4}(x)-L^{2s}_K[(1-\eta)(u-u(x_0))](x),\nonumber\\
\end{equation}
such that
\begin{equation}
\int_{B_4}|\nabla^s_{\mathbb{R}^n}(\eta (u-u(x_0))-v)|^2(x)dx\leq \epsilon^2\quad\hbox{ and }\quad\|\nabla^s_{B_1(x)}v(x)\|_{L^\infty(B_1)}\leq N_0,\nonumber\\
\end{equation}
where $\eta\in C^\infty_0(B_4)$ satisfying
\begin{equation}
0\leq\eta\leq1,\quad \eta\equiv1\hbox{ in }B_3,\quad\eta\equiv0\hbox{ in }\mathbb{R}^n\setminus B_4\hbox{ and }|\nabla\eta|\leq c(n).\nonumber\\
\end{equation}
\end{lem}

\begin{proof}
For any cut-off function $\eta\in C^\infty_0(B_4)$ satisfying
\begin{equation}
0\leq\eta\leq1,\quad \eta\equiv1\hbox{ in }B_3,\quad\eta\equiv0\hbox{ in }\mathbb{R}^n\setminus B_4\hbox{ and }|\nabla\eta|\leq c(n),\nonumber\\
\end{equation}
we have $\eta u\in\mathcal{H}^s_0(B_4)$ satisfying
\begin{equation}
 L^{2s}_{K}(\eta (u-u(x_0)))(x)=L^s_bf\chi_{B_4}(x)+g(x),\quad \hbox{ in }B_4,\nonumber\\
\end{equation}
where $g(x)=-L^{2s}_K[(1-\eta)(u-u(x_0))](x)+L^{s}_{b}(f\chi_{\mathbb{R}^n\setminus B_4})(x)$.
By Lemma 5.1, there exists a weak solution $v\in \mathcal{H}^s_0(B_4)$  of
\begin{equation}\label{5.3.2}
 L^{2s}_{K_0}v(x)=g(x),\quad \hbox{ in }B_4.\\
\end{equation}

For any $p>2$, take $p_0\in(2,p)$;  then we have $f\in L^{p_0}(B_8)\cap L^2(\mathbb{R}^n)$. By applying Theorem 4.1 to equation (\ref{1.1}) with $\frac{1}{2}B$ and $B$ replaced by  $B_5$, and $B_{\frac{11}{2}}$, respectively, we have $u\in {W}^{s+\varrho_0,2+\tau_0}(B_5)$ for some $\varrho_0$, $\tau_0\in(0,1-s)$. Then, applying Lemma 5.2 to $u-u(x_0)$ with $\frac{1}{2}B$ and $B$ replaced by  $B_{\frac{11}{2}}$, and $B_{6}$, respectively, we get
\begin{eqnarray}
[u]_{{W}^{s+\varrho_0,2+\tau_0}(B_5)}&=&\|\nabla^{s+\varrho_0}_{B_5}u\|_{L^{2+\tau_0}(B_5)}\nonumber\\
&\leq &C\Big(\|\nabla^s_{B_{\frac{11}{2}}}u\|_{L^2(B_{\frac{11}{2}})}+\Big(\int_{\mathbb{R}^n\setminus B_{\frac{11}{2}}} \frac{[u(x)-\bar{u}_{B_{\frac{11}{2}}}]^2}{|y|^{n+2s}}dy\Big)^{\frac{1}{2}}+
\|f\|_{L^{p_0}(B_{\frac{11}{2}})}+\Big(\int_{\mathbb{R}^n\setminus B_{\frac{11}{2}}}\frac{f^2(x)}{|x|^{n+s}}dx\Big)^{\frac{1}{2}}\Big)\nonumber\\
&\leq &C\Big(\|u-u(x_0)\|_{L^2(B_{6})}+\Big(\int_{\mathbb{R}^n\setminus B_6} \frac{[u(x)-u(x_0)]^2}{|x|^{n+2s}}dx\Big)^{\frac{1}{2}}+\|f\|_{L^{p_0}(B_{6})}+\Big(\int_{\mathbb{R}^n\setminus B_{\frac{11}{2}}}\frac{f^2(x)}{|x|^{n+s}}dx\Big)^{\frac{1}{2}}\Big).\nonumber
\end{eqnarray}
By a direct calculation, we have
\begin{eqnarray}
\|u-u(x_0)\|_{L^2(B_{6})}
\leq 7^{\frac{n+2s}{2}}\Big(\int_{B_6}\frac{|u(x)-u(x_0)|^2}{|x-x_0|^{n+2s}}dx\Big)^{\frac{1}{2}}
\leq 7^{\frac{n+2s}{2}}\nabla^s_{\mathbb{R}^n}u(x_0),\nonumber
\end{eqnarray}
\begin{eqnarray}
\Big(\int_{\mathbb{R}^n\setminus B_6} \frac{[u(x)-u(x_0)]^2}{|x|^{n+2s}}dx\Big)^{\frac{1}{2}}
\leq2^{\frac{n+2s}{2}}\Big(\int_{\mathbb{R}^n\setminus B_6} \frac{[u(x)-u(x_0)]^2}{|x-x_0|^{n+2s}}dx\Big)^{\frac{1}{2}}
\leq2^{\frac{n+2s}{2}}\nabla^s_{\mathbb{R}^n}u(x_0),\nonumber
\end{eqnarray}
\begin{eqnarray}
\|f\|_{L^{p_0}(B_{6})}
\leq |B_7(x_0)|^{\frac{1}{p_0}}\Big(\fint_{B_7(x_0)\cap B_8}|f|^{p_0}(x)dx\Big)^{\frac{1}{p_0}}
\leq |B_7(x_0)|^{\frac{1}{2}}\mathcal{M}^{\frac{1}{p_0}}_{B_8}(|f|^{p_0})(x_0),\nonumber
\end{eqnarray}
and
\begin{eqnarray}
\Big(\int_{\mathbb{R}^n\setminus B_{\frac{11}{2}}}\frac{f^2(x)}{|x|^{n+s}}dx\Big)^{\frac{1}{2}}
&\leq&c(n,s)\Big[\|f\|_{L^{2}(B_6)}+\Big(\int_{\mathbb{R}^n\setminus B_6}\frac{|f|^{2}(x)}{|x|^{n+s}}dx\Big)^{\frac{1}{2}}\Big]\nonumber\\
&\leq&c_1(n,s)\Big[\|f\|_{L^{p_0}(B_6)}+\Big(\int_{\mathbb{R}^n\setminus B_5(x_0)}\frac{|f|^{2}(x)}{|x-x_0|^{n+s}}dx\Big)^{\frac{1}{2}}\Big]\nonumber\\
&\leq&C(n,s)\Big[\mathcal{M}^{\frac{1}{p_0}}_{B_8}(|f|^{p_0})(x_0)+\Big(\int_{\mathbb{R}^n\setminus B_5(x_0)}\frac{|f|^{2}(x)}{|x-x_0|^{n+s}}dx\Big)^{\frac{1}{2}}\Big].\nonumber
\end{eqnarray}
Here we have used the fact that 
$$
|x-x_0|\leq |x|+|x_0|\leq |x|(1+\frac{|x_0|}{|x|})\leq 2|x|, ~~~\text{for any } x\in \mathbb{R}^n\setminus B_{\frac {11}{2}}.
$$
Therefore, we obtain
\begin{eqnarray}
[u]_{{W}^{s+\varrho_0,2+\tau_0}(B_5)}
&\leq &C\Big[\nabla^s_{\mathbb{R}^n}u(x_0)+\mathcal{M}^{\frac{1}{p_0}}_{B_8}(|f|^{p_0})(x_0)+\Big(\int_{\mathbb{R}^n\setminus B_5(x_0)}\frac{|f|^{2}(x)}{|x-x_0|^{n+s}}dx\Big)^{\frac{1}{2}}\Big].\nonumber
\end{eqnarray}
Then, by utilizing  Theorem 4.1 to $v-\bar{v}_{B_4}$ of equation (\ref{5.3.2}) with $\frac{1}{2}B$ and $B$ replaced by  $B_3$ and $B_{4}$, respectively, we get $v\in{W}^{s+\varrho,2+\tau}(B_3)$ for some $\varrho\in (0,1-s)$ and $\tau>0$. Moreover,  by  Lemma 5.1, we have
\begin{eqnarray}\label{4.00}
[v]_{{W}^{s+\varrho,2+\tau}(B_3)}&=&\|\nabla^{s+\varrho}_{B_3}v\|_{L^{2+\tau}(B_3)}\nonumber\\
&\leq& C\Big(\|\nabla^s_{B_4}v\|_{L^2(B_4)}+\|\nabla^{s+\varrho_0}_{B_4}(1-\eta)[u-u(x_0)]\|_{L^2(B_4)}
+\|(1-\eta)(u-u(x_0))\|_{L^2(B_4)}+\|f\|_{L^{p_0}(B_4)}\nonumber\\
& &+\Big(\int_{\mathbb{R}^n\setminus B_4}\frac{|v-\bar{v}_{B_4}|^2}{|x|^{n+2s}}dx\Big )^{\frac 12}+ \Big ( \int_{\mathbb{R}^n\setminus B_4}\frac{|[1-\eta(x)][u(x)-u(x_0)]-\overline{(1-\eta)[u-u(x_0)]}_{B_4}|^2}{|x|^{n+2s}}dx\nonumber\Big)^{\frac 12}\\
& &+\Big(\int_{\mathbb{R}^n\setminus B_4}\frac{|f|^2(x)}{|x|^{n+s}}dx  \Big)^{\frac 12})\nonumber\\
&\leq&C\Big(\|\nabla^s_{\mathbb{R}^n}v\|_{L^2(B_4)}+\|\nabla^{s+\varrho_0}_{B_4}(1-\eta)[u-u(x_0)]\|_{L^2(B_4)}
+\|u-u(x_0)\|_{L^2(B_4)}+\|f\|_{L^{p_0}(B_4)}\nonumber\\
& &+\|\nabla^{s}_{\mathbb{R}^n}(1-\eta)[u-u(x_0)]\|_{L^2(B_4)}+\Big(\int_{\mathbb{R}^n\setminus B_4}\frac{|f|^2(x)}{|x|^{n+s}}dx \Big)^{\frac 12}\Big)\nonumber\\
&\leq& C\Big(\|\nabla^s_{\mathbb{R}^n}(1-\eta)[u-u(x_0)]\|_{L^2(B_4)}+\|\nabla^{s+\rho_0}_{B_4}(1-\eta)(u-u(x_0))\|_{L^2(B_4)}
\nonumber\\
& &+\|u-u(x_0)\|_{L^2(B_4)}+\|f\|_{L^{p_0}(B_4)}+\Big(\int_{\mathbb{R}^n\setminus B_4}\frac{|f|^2(x)}{|x|^{n+s}}dx \Big)^{\frac{1}{2}}\Big)\nonumber\\
&\leq &C\Big[\nabla^s_{\mathbb{R}^n}u(x_0)+\mathcal{M}^{\frac{1}{p_0}}_{B_8}(|f|^{p_0})(x_0)+\Big(\int_{\mathbb{R}^n\setminus B_5(x_0)}\frac{|f|^{2}(x)}{|x-x_0|^{n+s}}dx\Big)^{\frac{1}{2}}\Big],\nonumber
\end{eqnarray}
where we used the facts that
\begin{eqnarray}
&&\|\nabla^s_{\mathbb{R}^n}(1-\eta)[u-u(x_0)]\|^2_{L^2(B_4)}\nonumber\\
&=&\int_{B_4}\int_{\mathbb{R}^n}\frac{[(1-\eta(x))(u(x)-u(x_0))-(1-\eta(y))(u(y)-u(x_0))]^2}{|x-y|^{n+2s}}dydx\nonumber\\
&\leq&2\int_{B_4}\int_{\mathbb{R}^n}\frac{[1-\eta(y)]^2[u(x)-u(y)]^2}{|x-y|^{n+2s}}dydx
+2\int_{B_4}\int_{\mathbb{R}^n}\frac{[\eta(x)-\eta(y)]^2[u(x)-u(x_0)]^2}{|x-y|^{n+2s}}dydx\nonumber\\
&\leq&2\|\nabla^s_{\mathbb{R}^n}u\|^2_{L^2(B_4)}+c(n,s)\|u-u(x_0)\|^2_{L^2(B_4)}\nonumber\\
&\leq &C\Big[\nabla^s_{\mathbb{R}^n}u(x_0)+\mathcal{M}^{\frac{1}{p_0}}_{B_8}(|f|^{p_0})(x_0)+\Big(\int_{\mathbb{R}^n\setminus B_4}\frac{|f|^2(x)}{|x|^{n+s}}dx \Big)^{\frac 12}\Big]^2,\nonumber
\end{eqnarray}
(where the last inequality follows  by applying Lemma 5.2 to $u-u(x_0)$ with $\frac{1}{2}B$ and $B$ replaced by $B_4$ and $B_6$, respectively.) and
\begin{eqnarray}
\|\nabla^{s+\varrho_0}_{B_4}(1-\eta)[u-u(x_0)]\|^2_{L^2(B_4)}
&\leq&2\|\nabla^{s+\varrho_0}_{B_4}u\|^2_{L^2(B_4)}+c(n,s)\|u-u(x_0)\|^2_{L^2(B_4)}\nonumber\\
&\leq&c([u]^{2}_{W^{s+\varrho_0,2+\tau_0}(B_4)}+\|u-\bar{u}_{B_4}\|^2_{L^{2+\tau_0}(B_4)})+c_1(n,s)\|u-u(x_0)\|^2_{L^2(B_4)}\nonumber\\
&\leq &C\Big(\nabla^{s}_{\mathbb{R}^n}u(x_0)+\mathcal{M}^{\frac{1}{p_0}}_{B_8}(|f|^{p_0})(x_0)+\Big(\int_{\mathbb{R}^n\setminus B_5(x_0)}\frac{|f|^{2}(x)}{|x-x_0|^{n+s}}dx\Big)^{\frac{1}{2}}\Big)^2,\nonumber
\end{eqnarray}
where we have used the fact that $W^{s+\varrho_0,2+\tau_0}(B_4)\hookrightarrow W^{s+\varrho_0,2}(B_4)$ and the fractional Sobolev-Poincar\'{e} inequality. The constant $C$ depends on $n$, $s$, $\lambda$, $\Lambda$, $\beta$ and $p$.

Notice that  $\eta [u-u(x_0)]-v\in\mathcal{H}^s_0(B_4)$ is a weak solution of
\begin{equation}
 L^{2s}_{{K}}(\eta (u-u(x_0))-v)(x)=L^{2s}_{K_0-{K}}v+L^s_bf\chi_{B_4}(x),\quad \hbox{ in }B_4.\nonumber\\
\end{equation}
Take the test function $\phi=\eta (u-u(x_0))-v$,  we have
\begin{eqnarray}\label{5.3.3}
&&\lambda[\eta (u-u(x_0))-v]^2_{H^s(\mathbb{R}^n)}\nonumber\\
&\leq&\int_{B_3}\int_{B_3}[\bar{A}_{B_3\times B_3}-A(x,y)]\frac{[v(x)-v(y)][(\eta (u-u(x_0))-v)(x)-(\eta (u-u(x_0))-v)(y)]}{|x-y|^{n+2s}}dydx\nonumber\\
& &+\int_{B_4}f(x)L^s_b(\eta (u-u(x_0))-v)(x)dx\nonumber\\
&=&I+II.
\end{eqnarray}
By H\"{o}lder's inequality, we get
\begin{eqnarray}
I&=&\int_{B_3}\int_{B_3}[\Bar{A}_{B_3\times B_3}-A(x,y)]\frac{[v(x)-v(y)][(\eta (u-u(x_0))-v)(x)-(\eta (u-u(x_0))-v)(y)]}{|x-y|^{n+2s}}dydx\nonumber\\
&\leq&\Big(\int_{B_3}\int_{B_3}|\Bar{A}_{B_3\times B_3}-A(x,y)|^{\frac{2(2+\gamma)}{\gamma}}dydx\Big)^{\frac{\gamma}{2(2+\gamma)}}
\cdot\Big(\int_{B_3}\int_{B_3}\frac{|v(x)-v(y)|^{2+\gamma}}{|x-y|^{n+(s+\frac{n\gamma}{2(2+\gamma)})(2+\gamma)}}dydx\Big)^{\frac{1}{2+\gamma}}\nonumber\\
& &\cdot\Big(\int_{B_3}\int_{B_3}\frac{[(\eta (u-u(x_0))-v)(x)-(\eta (u-u(x_0))-v)(y)]^2}{|x-y|^{n+2s}}dydx\Big)^{\frac{1}{2}}\nonumber\\
&\leq&\Big(\int_{B_3}\int_{B_3}|\Bar{A}_{B_3\times B_3}-A(x,y)|^{\frac{2(2+\gamma)}{\gamma}}dydx\Big)^{\frac{\gamma}{2(2+\gamma)}}[\eta( u-u(x_0))-v]_{H^s(\mathbb{R}^n)}[v]_{W^{s+\frac{n\gamma}{2(2+\gamma)},2+\gamma}(B_3)}\nonumber\\
&\leq&(2\Lambda)^{\frac{2}{(2+\gamma)}}|B_3\times B_3|^{\frac{\gamma}{2(2+\gamma)}}\Big(\fint_{B_3}\fint_{B_3}|\Bar{A}_{B_3\times B_3}-A(x,y)|^2dydx\Big)^{\frac{\gamma}{2(2+\gamma)}}[v]_{W^{s+\frac{n\gamma}{2(2+\gamma)},2+\gamma}(B_3)}[\eta (u-u(x_0))-v]_{H^s(\mathbb{R}^n)}\nonumber\\
&\leq&\delta^{\frac{\gamma}{2+\gamma}}(2\Lambda)^{\frac{2}{(2+\gamma)}}|B_3\times B_3|^{\frac{\gamma}{2(2+\gamma)}}[v]_{W^{s+\frac{n\gamma}{2(2+\gamma)},2+\gamma}(B_3)}[\eta (u-u(x_0))-v]_{H^s(\mathbb{R}^n)},\nonumber
\end{eqnarray}
where $\gamma$ is a positive constant satisfying $\gamma\leq\tau$,
$\frac{n\gamma}{2(2+\gamma)}\leq\varrho$. Then by the fact that $W^{s+\varrho,2+\tau}(B_3)\hookrightarrow W^{s+\frac{n\gamma}{2(2+\gamma)},2+\gamma}(B_3)$ and fractional Poincar\'{e} inequality,
we have
\begin{eqnarray}\label{5.3.4}
[v]_{W^{s+\frac{n\gamma}{2(2+\gamma)},2+\gamma}(B_3)}&=&[v-\bar{v}_{B_3}]_{W^{s+\frac{n\gamma}{2(2+\gamma)},2+\gamma}(B_3)}\nonumber\\
&\leq& C(n,s,\tau,\varrho,\gamma)\Big([v]_{{W^{s+\varrho,2+\tau}(B_3)}}+\|v-\bar{v}_{B_3}\|_{L^{2+\tau}(B_3)}\Big)\nonumber\\
&\leq& \hat{C}(n,s,\tau,\varrho,\gamma)[v]_{{W}^{s+\varrho,2+\tau}(B_3)}\nonumber\\
&\leq&C_1(n,s,\Lambda,\beta,p)\Big(\nabla^s_{\mathbb{R}^n}u(x_0)+\mathcal{M}^{\frac{1}{p_0}}_{B_8}(|f|^{p_0})(x_0)
+\Big(\int_{\mathbb{R}^n\setminus B_5(x_0)}\frac{|f|^{2}(x)}{|x-x_0|^{n+s}}dx\Big)^{\frac{1}{2}}\Big)\nonumber\\
&\leq&C_1(n,s,\Lambda,\beta,p)(1+\delta).
\end{eqnarray}
Notice that, by H\"{o}lder's inequality, assumptions (\ref{5.3.1}) and (\ref{1.3.1}), we get
\begin{eqnarray}\label{5.3.5}
|II|&=&\Big|\int_{B_4}f(x)L^s_b(\eta (u-u(x_0))-v)(x)dx\Big|\nonumber\\
&\le&\|f\|_{L^2(B_4)}\|L^s_b(\eta (u-u(x_0))-v)\|_{L^2(B_4)}\nonumber\\
&\leq&  C(n,s,\beta)\delta|B_4|^{\frac{1}{2}}[\eta (u-u(x_0))-v]_{H^s(\mathbb{R}^n)}.
\end{eqnarray}
Combining (\ref{5.3.3})-(\ref{5.3.5}), we derive that
\begin{eqnarray}
\|\nabla^s_{\mathbb{R}^n}(\eta (u-u(x_0))-v)\|_{L^2(B_4)}\leq[\eta (u-u(x_0))-v]_{H^s(\mathbb{R}^n)}\leq C[\delta^{\frac{\gamma}{2+\gamma}}(1+\delta)+\delta]<\epsilon,\nonumber
\end{eqnarray}
if we take $\delta>0$ small enough, where $C$ is a constant depending on $n$, $s$, $\lambda$, $\Lambda$, $\beta$ and $p$.

Notice that for any $\phi\in C^\infty_0(B_{\frac{5}{2}})$, by H\"{o}lder's inequality, we can deduce that
\begin{eqnarray}
\int_{\mathbb{R}^n}L^{2s}_K[(1-\eta)(u-u(x_0))](x)\phi(x)dx&=&\int_{\mathbb{R}^n}\int_{\mathbb{R}^n}K(x,y)[(1-\eta)(u-u(x_0))](x)[\phi(x)-\phi(y)]dydx\nonumber\\
&=&-\int_{\mathbb{R}^n\setminus B_3}\int_{B_{\frac{5}{2}}}K(x,y)[(1-\eta)(u-u(x_0))](x)\phi(y)dydx\nonumber\\
&\leq&\Lambda\int_{\mathbb{R}^n\setminus B_3}\int_{B_{\frac{5}{2}}}\frac{|u(x)-u(x_0)|}{|x-y|^{n+2s}}|\phi(y)|dydx\nonumber\\
&\leq&c(n,s,\Lambda)\int_{\mathbb{R}^n\setminus B_3}\frac{|u(x)-u(x_0)|}{|x|^{n+2s}}dx\int_{B_{\frac{5}{2}}}|\phi(y)|dy,\nonumber
\end{eqnarray}
and
\begin{eqnarray}
\int_{\mathbb{R}^n}L^{s}_bf\chi_{\mathbb{R}^n\setminus B_4}(x)\phi(x)dx&=&\int_{\mathbb{R}^n}\int_{\mathbb{R}^n}b(x,y)f\chi_{\mathbb{R}^n\setminus B_4}(x)[\phi(x)-\phi(y)]dydx\nonumber\\
&=&-\int_{\mathbb{R}^n\setminus B_4}\int_{B_{\frac{5}{2}}}b(x,y)f(x)\phi(y)dydx\nonumber\\
&\leq&\beta\int_{\mathbb{R}^n\setminus B_4}\int_{B_{\frac{5}{2}}}\frac{|f(x)|}{|x-y|^{n+s}}|\phi(y)|dydx\nonumber\\
&\leq& C(n,s,\beta)\int_{\mathbb{R}^n\setminus B_4}\frac{|f(x)|}{|x|^{n+s}}dx\int_{B_{\frac{5}{2}}}|\phi(y)|dy,\nonumber
\end{eqnarray}
which imply that
\begin{equation}
\|g\|_{L^\infty(B_{\frac{5}{2}})}\leq C(n,s,\Lambda,\beta)\Big(\int_{\mathbb{R}^n\setminus B_3}\frac{|u(x)-u(x_0)|}{|x|^{n+2s}}dx
+\int_{\mathbb{R}^n\setminus B_4}\frac{|f(x)|}{|x|^{n+s}}dx\Big).\nonumber
\end{equation}
By Theorem 3.1 (with $B_1$  and $B_3$ replaced by $B_2$ and $B_{\frac{5}{2}}$, respectively ), we can get $v\in C^\alpha(B_2)$ for every $\alpha\in(0,\min\{2s,1\})$ and
$$
[v]_{C^\alpha(B_2)}\leq C\Big(\|\nabla^s_{\mathbb{R}^n}v\|_{L^2(B_{\frac{5}{2}})}+\|g\|_{L^\infty(B_{\frac{5}{2}})})\leq C(\|\nabla^s_{\mathbb{R}^n} v\|_{L^2(B_{\frac{5}{2}})}+\int_{\mathbb{R}^n\setminus B_3}\frac{|u(x)-u(x_0)|}{|x|^{n+2s}}dx
+\int_{\mathbb{R}^n\setminus B_4}\frac{|f(x)|}{|x|^{n+s}}dx\Big).
$$
By Lemma 5.1, we have
\begin{eqnarray}
\|\nabla^s_{\mathbb{R}^n} v\|_{L^2(B_{\frac{5}{2}})}&\leq&\|\nabla^s_{\mathbb{R}^n}v\|_{L^2(B_4)}\nonumber\\
&\leq&c\Big[\|\nabla^s_{\mathbb{R}^n}(1-\eta)(u-u(x_0))\|_{L^2(B_4)}+\|f\|_{L^2(B_4)}+\Big(\int_{\mathbb{R}^n\setminus B_4}\frac{|f|^2(x)}{|x|^{n+s}}dx \Big)^{\frac{1}{2}}\Big]\nonumber\\
&\leq& C\Big(\nabla^s_{\mathbb{R}^n}u(x_0)+\mathcal{M}^{\frac{1}{p_0}}_{B_8}(|f|^{p_0})(x_0)+\Big(\int_{\mathbb{R}^n\setminus B_6}\frac{|f|^2(x)}{|x|^{n+s}}dx \Big)^{\frac{1}{2}}\Big).\nonumber
\end{eqnarray}
By a  computation similar to the  above, we get
\begin{eqnarray}
\int_{\mathbb{R}^n\setminus {B_3}}\frac{|u(x)-u(x_0)|}{|x|^{n+2s}}dx
\leq c(n,s)\Big(\int_{\mathbb{R}^n\setminus {B_3}}\frac{|u(x)-u(x_0)|^2}{|x|^{n+2s}}dx\Big)^{\frac{1}{2}}
\leq c_1(n,s)\nabla^s_{\mathbb{R}^n}u(x_0),\nonumber
\end{eqnarray}
and
\begin{eqnarray}
\int_{\mathbb{R}^n\setminus B_4}\frac{|f(x)|}{|x|^{n+s}}dx
&\leq& c(n,s)\Big[\|f\|_{L^{p_0}(B_6)}+\Big(\int_{\mathbb{R}^n\setminus {B_6}}\frac{|f|^{2}(x)}{|x|^{n+2s}}dx\Big)^{\frac{1}{2}}\Big]\nonumber\\
&\leq& C(n,s)\Big[\mathcal{M}^{\frac{1}{p_0}}_{B_8}(|f|^{p_0})(x_0)+\Big(\int_{\mathbb{R}^n\setminus B_5(x_0)}\frac{|f|^{2}(x)}{|x-x_0|^{n+s}}dx\Big)^{\frac{1}{2}}\Big].\nonumber
\end{eqnarray}
Therefore,  we have
$$
[v]_{C^\alpha(B_2)}\leq C\Big[\nabla^s_{\mathbb{R}^n}u(x_0)+\mathcal{M}^{\frac{1}{p_0}}_{B_8}(|f|^{p_0})(x_0)+\Big(\int_{\mathbb{R}^n\setminus B_5(x_0)}\frac{|f|^{2}(x)}{|x-x_0|^{n+s}}dx\Big)^{\frac{1}{2}}\Big],\nonumber
$$
where  the constant $C$ depends on $n$, $s$, $\lambda$, $\Lambda$, $\alpha$, $\beta$ and $p$.

Take some $\vartheta\in(s,\min\{2s,1\})$,  there exists $N_0>0$ such that for any $x\in B_1$,
\begin{eqnarray}
(\nabla^s_{B_1(x)}v(x))^2
&=&\int_{B_1(x)}\frac{[v(x)-v(y)]^2}{|x-y|^{n+2s}}dy\nonumber\\
&\leq&[v]^2_{C^\vartheta(B_2)}\int_{B_1(x)}\frac{1}{|x-y|^{n+2s-2\vartheta}}dy\nonumber\\
&\leq& C(n,s,\lambda,\Lambda,\beta)\Big[\nabla^s_{\mathbb{R}^n}u(x_0)+\mathcal{M}^{\frac{1}{p_0}}_{B_8}(|f|^{p_0})(x_0)+\Big(\int_{\mathbb{R}^n\setminus B_5(x_0)}\frac{|f|^{2}(x)}{|x-x_0|^{n+s}}dx\Big)^{\frac{1}{2}}\Big]\nonumber\\
&\leq& C(n,s,\lambda,\Lambda,\beta) (1+\delta)\nonumber\\
&<&N^2_0.\nonumber
\end{eqnarray}
\end{proof}

\begin{lem}Let $2<p_0<p$. There is a constant $N_1>0$ such that for any $\epsilon>0$, there exists a small  $\delta=\delta(n,s,\lambda, \Lambda,\beta, p,\epsilon)>0$, and if $u\in{H}^s(\mathbb{R}^n)$ is a weak solution of (\ref{1.1}) in $\Omega\supset B_8$ with $A$ $(\delta,8)$-vanishing and $f\in L^2({\mathbb{R}^n})\cap L^p(B_8)$ such that
\begin{equation}
\{x\in B_1|\hbox{ } \nabla^s_{\mathbb{R}^n}u(x)\leq1\}\cap\{x\in B_1|\hbox{ }\mathcal{M}_{B_8}(|f|^{p_0})(x)\leq\delta^{p_0}\}\cap\Big\{x\in B_1\Big|\hbox{ } \int_{\mathbb{R}^n\setminus B_5(x)}\frac{|f|^2(y)}{|y-x|^{n+s}}dy\leq\delta^2\Big\}\neq{\O},\nonumber\\
\end{equation}
then
$$
|\{x\in B_1|\hbox{ }\nabla^s_{\mathbb{R}^n}u(x)>N_1\}|<\epsilon |B_1|.
$$
\end{lem}
\begin{proof}
Assume there exists $x_1\in B_1$ such that
\begin{equation}
\nabla^s_{\mathbb{R}^n}u(x_1)\leq1,\quad \mathcal{M}_{B_8}(|f|^{p_0})(x_1)\leq\delta^{p_0}\hbox{ and }\int_{\mathbb{R}^n\setminus B_5(x_1)}\frac{|f|^2(y)}{|y-x_1|^{n+s}}dy\leq\delta^2.\nonumber\\
\end{equation}
Applying Lemma 5.3, we  find that there exists $N_0>0$ such that for any $\mu>0$, there exist $\delta=\delta(n,s,\Lambda,\beta,\mu)>0$ and a weak solution $v\in \mathcal{H}^s_0(B_4)$ of
\begin{equation}
\left\{
\begin{array}{lll}
 L^{2s}_{K_0}v(x)=g,\quad &\hbox{ in }B_4,\nonumber\\
\quad\quad v=0,\quad &\hbox{ in }\mathbb{R}^n\setminus B_4,\nonumber
\end{array}\right.
\end{equation}
satisfying
\begin{equation}
\|\nabla^s_{\mathbb{R}^n}(\eta (u-u(x_1))-v)\|_{L^2(B_4)}<\mu\hbox{ and }\|\nabla^s_{B_1(x)}v(x)\|_{L^\infty(B_1)}\leq N_0.\nonumber\\
\end{equation}

Now we claim that there exists $N_1=N_1(n,s,N_0)$ such that
$$
\{x\in B_1|\nabla^s_{\mathbb{R}^n}u(x)>N_1\}\subset\{x\in B_1|\nabla^s_{B_1(x)}(u-v)(x)>N_0\}.
$$

Indeed, assume there is $x_2\in B_1$ such that $\nabla^s_{B_1(x_2)}(u-v)(x_2)\leq N_0$, then
\begin{eqnarray}\label{5.4.1}
\Big(\nabla^s_{B_1(x_2)}u\Big)^2(x_2)
&\leq&2\Big[\Big(\nabla^s_{B_1(x_2)}(u-v)\Big)^2(x_2)
+\Big(\nabla^s_{B_1(x_2)}(v)\Big)^2(x_2)\Big]\nonumber\\
&\leq&4N^2_0,
\end{eqnarray}
which implies that
\begin{equation}\label{5.4.2}
\int_{B_1(x_2)\setminus B_{\frac{1}{2}}(x_2)}[u(x_2)-u(y)]^2dy\leq\int_{B_1(x_2)}\frac{[u(x_2)-u(y)]^2}{|x_2-y|^{n+2s}}dy=\Big(\nabla^s_{B_1(x_2)}u\Big)^2(x_2)\leq4N^2_0.
\end{equation}
Since  $\nabla^s_{\mathbb{R}^n}u(x_1)\leq1$, we have
\begin{equation}
\int_{B_1(x_2)\setminus B_{\frac{1}{2}}(x_2)}\frac{[u(x_1)-u(y)]^2}{|x_1-y|^{n+2s}}dy\leq1.\nonumber\\
\end{equation}
Notice the fact that
$$
|x_1-y|\leq|x_2-y|+|x_1-x_2|=|x_2-y|\Big(1+\frac{|x_1-x_2|}{|x_2-y|}\Big)\leq5|x_2-y|, \hbox{ for any }y\in B_1(x_2)\setminus B_{\frac{1}{2}}(x_2),
$$
then we  get
\begin{eqnarray}\label{5.4.3}
\int_{B_1(x_2)\setminus B_{\frac{1}{2}}(x_2)}[u(x_1)-u(y)]^2dy
&\leq&\int_{B_1(x_2)\setminus B_{\frac{1}{2}}(x_2)}\frac{[u(x_1)-u(y)]^2}{|x_2-y|^{n+2s}}dy\nonumber\\
&\leq&5^{n+2s}\int_{B_1(x_2)\setminus B_{\frac{1}{2}}(x_2)}\frac{[u(x_1)-u(y)]^2}{|x_1-y|^{n+2s}}dy\nonumber\\
&\leq&5^{n+2s}.
\end{eqnarray}
Combining (\ref{5.4.2}) and (\ref{5.4.3}), we obtain  $[u(x_1)-u(x_2)]^2\leq C$,
where $C=2(4N^2_0 +5^{n+2s})/|B_1\setminus B_{1/2}|$.

Now for any $y\in\mathbb{R}^n\setminus B_1(x_2)$, noticing that $|x_1-y|\leq|x_2-y|+|x_1-x_2|\leq3|x_2-y|$, we have
\begin{eqnarray}
&&\Big(\nabla^s_{\mathbb{R}^n\setminus B_1(x_2)}u(x_2)\Big)^2\nonumber\\
&=&\int_{\mathbb{R}^n\setminus B_1(x_2)}\frac{[u(x_2)-u(y)]^2}{|x_2-y|^{n+2s}}dy\nonumber\\
&\leq&2\int_{\mathbb{R}^n\setminus B_1(x_2)}\frac{[u(x_2)-u(x_1)]^2}{|x_2-y|^{n+2s}}dy+2\int_{\mathbb{R}^n\setminus B_1(x_2)}\frac{[u(x_1)-u(y)]^2}{|x_2-y|^{n+2s}}dy\nonumber\\
&\leq&2\int_{\mathbb{R}^n\setminus B_1(x_2)}\frac{[u(x_2)-u(x_1)]^2}{|x_2-y|^{n+2s}}dy+2\cdot3^{n+2s}\int_{\mathbb{R}^n\setminus B_1(x_2)}\frac{[u(x_1)-u(y)]^2}{|x_1-y|^{n+2s}}dy\nonumber\\
&\leq&\frac{4(4N^2_0+5^{n+2s})}{|B_1\setminus B_{\frac{1}{2}}|}\int_{\mathbb{R}^n\setminus B_1(x_2)}\frac{1}{|x_2-y|^{n+2s}}dy+2\cdot 3^{n+2s}\nabla^s_{\mathbb{R}^n}u(x_1)\nonumber\\
&\leq&C(n,s)(1+N^2_0),\nonumber
\end{eqnarray}
which along with (\ref{5.4.1}) imply that
\begin{eqnarray}
\Big(\nabla^s_{\mathbb{R}^n}u(x_2)\Big)^2&=&\Big(\nabla^s_{B_1(x_2)}u(x_2)\Big)^2+\Big(\nabla^s_{\mathbb{R}^n\setminus
B_1(x_2)}u(x_2)\Big)^2\nonumber\\
&\leq& 4N^2_0 +C(n,s)(1+N^2_0)=N^2_1.\nonumber
\end{eqnarray}
Thus the claim is proved. Then it follows that
\begin{eqnarray}
|\{x\in B_1|\nabla^s_{\mathbb{R}^n}u(x)>N_1\}|
&\leq&|\{x\in B_1|\nabla^s_{B_1(x)}(u-v)(x)>N_0\}|\nonumber\\
&\leq&\frac{1}{N^2_0}\int_{B_1}\Big(\nabla^s_{B_2}(u-v)\Big)^2(x)dx\nonumber\\
&\leq&\frac{1}{N^2_0}\int_{B_4}\Big(\nabla^s_{\mathbb{R}^n}(\eta (u-u(x_1))-v)\Big)^2(x)dx\nonumber\\
&\leq& \mu^2=\epsilon|B_1|,\nonumber
\end{eqnarray}
by taking $\mu>0$ satisfying the last identity above. This finishes the proof.
\end{proof}

\begin{lem} Let $2<p_0<p$. There exists a constant $N_1>0$ such that for any $\epsilon>0$ there exists a small  $\delta=\delta(n,s,\lambda, \Lambda,\beta, p, \epsilon, R)>0$, and if $u\in {H}^s(\mathbb{R}^n)$ is a weak solution of (\ref{1.1}) in $\Omega\supset 8B_r(y_0)=8B$ with  $A(x,y)$ $(\delta, 8r)$-vanishing  and $f\in L^2(\mathbb{R}^n)\cap L^p(8B)$ such that
\begin{equation}
\{x\in B|\hbox{ } \nabla^s_{\mathbb{R}^n}u(x)\leq1\}\cap\{x\in B|\hbox{ }\mathcal{M}_{8B}(|f|^{p_0})(x)\leq\delta^{p_0}\}\cap\Big\{x\in B\Big|\hbox{ } r^s\int_{\mathbb{R}^n\setminus B_{5r}(x)}\frac{|f|^2(y)}{|y-x|^{n+s}}dy\leq\delta^2\Big\}\neq{\O},\nonumber\\
\end{equation}
then
$$
|\{x\in B|\hbox{ }\nabla^s_{\mathbb{R}^n}u(x)>N_1\}|<\epsilon |B|.
$$
\end{lem}
\begin{proof}
For any $x\in B_1$, let $u_r(x)=r^{-s}u(rx+y_0)$, $f_r(x)=f(rx+y_0)$, $K_r(x,y)=r^{n+2s}K(rx+y_0,ry+y_0)$  and $b_r(x,y)=r^{n+s}b(rx+y_0,ry+y_0)$, then by a direct calculation, we can see that $u_r$ is a weak solution of
$$
L^{2s}_{K_r}u_r(x)=L^s_{b_r}f_r(x)\quad\hbox{ in }\Big\{\frac{z-y_0}{r}|\hbox{ }z\in\Omega\Big\}\supset B_8.
$$
 By Lemma 5.3, there is $N_1>0$ such that for any $\epsilon>0$, there is $\delta>0$, if
\begin{equation}
\{x\in B_1| \hbox{ }\nabla^s_{\mathbb{R}^n}u_r(x)\leq1\}\cap\{x\in B_1|\hbox{ }\mathcal{M}_{B_8}(|f_r|^{p_0})(x)\leq\delta^{p_0}\}\cap\Big\{x\in B_1\Big|\hbox{ } \int_{\mathbb{R}^n\setminus B_{5}(x)}\frac{|f_r|^2(y)}{|y-x|^{n+s}}dy\leq\delta^2\Big\}\neq{\O},\nonumber\\
\end{equation}
then
$$
|\{x\in B_1|\hbox{ }\nabla^s_{\mathbb{R}^n}u_r(x)>N_1\}|<\epsilon |B_1|,
$$
which means, for the above constants $\epsilon$, $\delta$ and $N_1$, if
\begin{equation}
\{x\in B|\hbox{ } \nabla^s_{\mathbb{R}^n}u(x)\leq1\}\cap\{x\in B|\hbox{ }\mathcal{M}_{8B}(|f|^{p_0})(x)\leq\delta^{p_0}\}\cap\Big\{x\in B\Big|\hbox{ } r^s\int_{\mathbb{R}^n\setminus B_{5r}(x)}\frac{|f|^2(y)}{|y-x|^{n+s}}dy\leq\delta^2\Big\}\neq{\O},\nonumber\\
\end{equation}
then
$$
|\{x\in B|\hbox{ }\nabla^s_{\mathbb{R}^n}u(x)>N_1\}|<\epsilon |B|,
$$
where we used the facts that
$
\nabla^s_{\mathbb{R}^n}u_r(x)=\nabla^s_{\mathbb{R}^n}u(rx+y_0)$,
\begin{equation}
\frac{1}{|B_{R}|}\int_{8B\cap B_{R}(z)}|f|^{p_0}(y)dy=\frac{1}{|B_{\frac{R}{r}}|}\int_{B_8\cap B_{\frac{R}{r}}\Big(\frac{z-y_0}{r}\Big)}|f_r|^{p_0}(x)dx,\nonumber
\end{equation}
and
\begin{equation}
r^s\int_{\mathbb{R}^n\setminus B_{5r}(z)}\frac{|f|^{p_0}(y)}{|y-z|^{n+s}}dy=\int_{\mathbb{R}^n\setminus B_{5}\Big(\frac{z-y_0}{r}\Big)}\frac{|f_r|^{p_0}(x)}{\Big|x-\frac{z-y_0}{r}\Big|^{n+s}}dx,\nonumber
\end{equation}
for any $z\in B$.
\end{proof}

 An immediate consequence of the above lemma is the following corollary.
\begin{cor}Let $2<p_0<p$. There is a constant $N_1>0$ such that for any $\epsilon>0$ there exists  $\delta=\delta(n,s,\lambda,\Lambda,\beta,\epsilon)>0$,  if $u\in{H}^s(\mathbb{R}^n)$  is a weak solution of (\ref{1.1}) in $\Omega\supset 8B_r(y_0)=8B$ with $A(x,y)$ $(\delta, 8r)$-vanishing and $f\in L^2(\mathbb{R}^n)\cap L^p(8B)$ such that 
$$
|\{x\in B|\hbox{ }\nabla^s_{\mathbb{R}^n}u(x)>N_1\}|\geq\epsilon |B|,
$$
then
\begin{equation}
B\cap B_1\subset\{x\in B_1|\hbox{ } \nabla^s_{\mathbb{R}^n}u(x)>1\}\cup\{x\in B_1|\hbox{ }\mathcal{M}_{8B}(|f|^{p_0})(x)>\delta^{p_0}\}\cup\Big\{x\in B_1\Big|\hbox{ } \int_{\mathbb{R}^n\setminus B_{5}(x)}\frac{|f|^2(y)}{|y-x|^{n+s}}dy>\delta^2\Big\}.\nonumber\\
\end{equation}
\end{cor}

\begin{lem}Let $2<p_0<p$.
Suppose $u\in{H}^s(\mathbb{R}^n)$ is a weak solution of (\ref{1.1}) in $\Omega$. Assume that $A(x,y)$ is $(\delta, 8)$-vanishing and $f\in L^2(\mathbb{R}^n)\cap L^p(B_8))$ such that
\begin{equation}
|\{x\in B_1|\hbox{ }\nabla^s_{\mathbb{R}^n}u(x)>N_1\}|<\epsilon |B_1|.\nonumber\\
\end{equation}
Let $k\in\mathbb{N}^+$, set $\epsilon_1=20^n\epsilon$, then
\begin{eqnarray}
&&|\{x\in B_1|\hbox{ }\nabla^s_{\mathbb{R}^n}u(x)>N^{k}_1\}|\nonumber\\
&\leq&\epsilon^k_1|\{x\in B_1|\hbox{ }\nabla^s_{\mathbb{R}^n}u(x)>1\}|
+\sum^k_{i=1}\epsilon^i_1|\hbox{ }\{x\in B_1|\hbox{ }\mathcal{M}_{B_8}(|f|^{p_0})(x)>\delta^{p_0} N^{p_0(k-i)}_1\}|\nonumber\\
& &+\sum^k_{i=1}\epsilon^i_1\Big|\Big\{x\in B_1|\hbox{ }\int_{\mathbb{R}^n\setminus B_{5}(x)}\frac{|f|^2(y)}{|y-x|^{n+s}}dy>\delta^{2} N^{2(k-i)}_1\Big\}\Big|.\nonumber
\end{eqnarray}
\end{lem}
\begin{proof}
We intend to prove this lemma by induction on $k$.  The conclusion is obtained for the case $k=1$ if we let
$$
C=\{x\in B_1|\hbox{ }\nabla^s_{\mathbb{R}^n}u(x)>N_1\}
$$
and
$$
D=\{x\in B_1|\hbox{ }\nabla^s_{\mathbb{R}^n}u(x)>1\}\cup\{x\in B_1|\hbox{ }\mathcal{M}_{B_8}(|f|^{p_0})(x)>\delta^{p_0}\}\cup\Big\{x\in B_1\Big|\hbox{ } \int_{\mathbb{R}^n\setminus B_{5}(x)}\frac{|f|^2(y)}{|y-x|^{n+s}}dy>\delta^2\Big\}.
$$
and then apply Lemma 2.3.

Assume that the conclusion is valid for some positive integer $k\geq2$. Set $u_1=u/{N_1}$ and $f_1=f/{N_1}$. Then $u_1$ is a weak solution of
\begin{equation}
L^{2s}_Ku_1(x)=L^s_bf_1(x)\quad\hbox{ in }\Omega\nonumber\\
\end{equation}
and
\begin{equation}
|\{x\in B_1|\hbox{ }\nabla^s_{\mathbb{R}^n}u_1(x)>N_1\}|<\epsilon |B_1|.\nonumber\\
\end{equation}
By the induction assumption and a direct calculation, we deduce that
\begin{eqnarray}
&&|\{x\in B_1|\hbox{ }\nabla^s_{\mathbb{R}^n}u(x)>N^{(k+1)}_1\}|\nonumber\\
&=&|\{x\in B_1|\hbox{ }\nabla^s_{\mathbb{R}^n}u_1(x)>N^{k}_1\}|\nonumber\\
&\leq&\epsilon^k_1|\{x\in B_1|\hbox{ }\nabla^s_{\mathbb{R}^n}u_1(x)>1\}|+\sum^k_{i=1}\epsilon^i_1|\{x\in B_1|\hbox{ }\mathcal{M}_{B_8}(|f_1|^{p_0})(x)>\delta^{p_0} N^{(k-i)p_0}_1\}|\nonumber\\
& &+\sum^k_{i=1}\epsilon^i_1\Big|\Big\{x\in B_1|\hbox{ }\int_{\mathbb{R}^n\setminus B_{5}(x)}\frac{|f_1|^2(y)}{|y-x|^{n+s}}dy>\delta^{2} N^{2(k-i)}_1\Big\}\Big|\nonumber\\
&=&\epsilon^k_1|\{x\in B_1|\hbox{ }\nabla^s_{\mathbb{R}^n}u(x)>N_1\}|+\sum^k_{i=1}\epsilon^i_1|\{x\in B_1|\hbox{ }\mathcal{M}_{B_8}(|f|^{p_0})(x)>\delta^{p_0} N^{(k+1-i)p_0}_1\}|\nonumber\\
& &+\sum^k_{i=1}\epsilon^i_1\Big|\Big\{x\in B_1|\hbox{ }\int_{\mathbb{R}^n\setminus B_{5}(x)}\frac{|f|^2(y)}{|y-x|^{n+s}}dy>\delta^{2} N^{2(k+1-i)}_1\Big\}\Big|\nonumber\\
&\leq&\epsilon^{k+1}_1|\{x\in B_1|\hbox{ }\nabla^s_{\mathbb{R}^n}u(x)>1\}|+\sum^{k+1}_{i=1}\epsilon^i_1|\{x\in B_1|\hbox{ }\mathcal{M}_{B_8}(|f|^{p_0})(x)>\delta^{p_0} N^{(k+1-i)p_0}_1\}|\nonumber\\
& &+\sum^{k+1}_{i=1}\epsilon^i_1\Big|\Big\{x\in B_1|\hbox{ }\int_{\mathbb{R}^n\setminus B_{5}(x)}\frac{|f|^2(y)}{|y-x|^{n+s}}dy>\delta^{2} N^{2(k+1-i)}_1\Big\}\Big|.\nonumber
\end{eqnarray}
Thus the conclusion is valid for the case $k+1$. This completes the proof.
\end{proof}

\begin{proof}[ Proof of Theorem 1.1.]

For $p>2$ fixed, take $\epsilon\in(0,1)$ such that $N^p_120^n\epsilon<\frac{1}{2}$. Let $N_1$ and $\delta$ be as in Lemma 5.4. Assume $\nabla^s_{\mathbb{R}^n}u\not\equiv0$ in $B_8$.  Let $\tilde{u}(x)=\frac{\theta u}{M}$ and $\tilde{f}(x)=\frac{\theta f}{M}$,  where $$M=\|u\|_{L^2(B_8)}+\Big(\int_{\mathbb{R}^n\setminus B_8}\frac{u^2(x)}{|x|^{n+2s}}dx\Big)^{\frac{1}{2}}+\|f\|_{L^{2}(B_8)}+\Big(\int_{\mathbb{R}^n\setminus B_8}\frac{|f|^{2}(x)}{|x|^{n+s}}dx\Big)^{\frac{1}{2}},$$
and $\theta>0$ is a constant. Then $\tilde{u}$ is a weak solution of
$$
L^{2s}_K\tilde{u}(x)=L^s_b\tilde{f}(x)\quad\text{ in }\Omega.
$$
Moreover, by Lemma 5.2, we have
\begin{eqnarray}
    \|\nabla^s_{\mathbb{R}^n}\tilde{u}\|_{L^2(B_4)}=\frac{\theta}{M}\|\nabla^s_{\mathbb{R}^n}{u}\|_{L^2(B_4)}\leq C\theta.\nonumber
\end{eqnarray}
We can take $\theta>0$ small such that
$$
|\{x\in B_1|\nabla^s_{\mathbb{R}^n}\tilde{u}(x)>N_1\}|<\epsilon|B_1|.
$$
which is due to the fact that
$$
|\{x\in B_1|\nabla^s_{\mathbb{R}^n}\tilde{u}(x)>N_1\}|\leq\frac{1}{N^2_1}\int_{B_4}(\nabla^s_{\mathbb{R}^n}\tilde{u})^2(x)dx\leq\frac{C^2\theta^2}{N^2_1}<\epsilon|B_1|,
$$
if $\theta>0$ small enough. Then by Lemma 5.6, Lemma 2.1 and Theorem 2.1, we have
\begin{eqnarray}
&&\sum^\infty_{k=1}N^{kp}_1\Big|\left\{x\in B_1|\nabla^s_{\mathbb{R}^n}\tilde{u}(x)>N^{k}_1\right\}\Big|\nonumber\\
&\leq&\sum^\infty_{k=1}N^{kp}_1\Big[\epsilon^k_1\Big|\left\{x\in B_1|\nabla^s_{\mathbb{R}^n}\tilde{u}(x)>1\right\}\Big|+\sum^k_{i=1}\epsilon^i_1\Big|\left\{x\in B_1|\mathcal{M}_{B_8}(|\tilde{f}|^{p_0})(x)>\delta^{p_0} N^{p_0(k-i)}_1\right\}\Big|\nonumber\\
& &+\sum^k_{i=1}\epsilon^i_1\Big|\Big\{x\in B_1|\hbox{ }\int_{\mathbb{R}^n\setminus B_{5}(x)}\frac{|\tilde{f}|^2(y)}{|y-x|^{n+s}}dy>\delta^{2} N^{2(k-i)}_1\Big\}\Big|\Big]\nonumber\\
&=&\sum^\infty_{k=1}(N^p_1\epsilon_1)^k\Big|\{x\in B_1|\nabla^s_{\mathbb{R}^n}\tilde{u}(x)>1\}\Big|\nonumber\\
& &+\sum^\infty_{k=0}N^{kp}_1\Big|\{x\in B_1|\mathcal{M}_{B_8}(|\tilde{f}|^{p_0})(x)>\delta^{p_0} N^{kp_0}_1\}\Big|\sum^\infty_{i=1}(N^p_1\epsilon_1)^i\nonumber\\
& &+\sum^\infty_{k=0}N^{kp}_1\Big|\Big\{x\in B_1|\hbox{ }\int_{\mathbb{R}^n\setminus B_{5}(x)}\frac{|\tilde{f}|^2(y)}{|y-x|^{n+s}}dy>\delta^{2} N^{2k}_1\Big\}\Big|\sum^\infty_{i=1}(N^p_1\epsilon_1)^i\nonumber\\
&\leq&\sum^\infty_{k=1}N^{kp}_1\Big|\{x\in B_1|\mathcal{M}_{B_8}(|\tilde{f}|^{p_0})(x)>\delta^{p_0} N^{kp_0}_1\}\Big|\nonumber\\
& &+\sum^\infty_{k=1}N^{kp}_1\Big|\Big\{x\in B_1|\hbox{ }\int_{\mathbb{R}^n\setminus B_{5}(x)}\frac{|\tilde{f}|^2(y)}{|y-x|^{n+s}}dy>\delta^{2} N^{2k}_1\Big\}\Big|+3|B_1|\nonumber\\
&\leq&C\Big[\int_{B_1}\mathcal{M}^{\frac{p}{p_0}}_{B_8}(|\tilde{f}|^{p_0})(x)dx+\int_{B_1}\Big(\int_{\mathbb{R}^n\setminus B_{5}(x)}\frac{|\tilde{f}|^2(y)}{|y-x|^{n+s}}dy\Big)^{\frac{p}{2}}dx\Big]+3|B_1|\nonumber\\
&\leq&\tilde{C}\Big[\int_{B_8}|\tilde{f}|^p(x)dx+\Big(\int_{\mathbb{R}^n\setminus B_8}\frac{|\tilde{f}|^2(y)}{|y|^{n+s}}dy\Big)^{\frac{p}{2}}+1\Big].\nonumber
\end{eqnarray}
Hence, by Lemma 2.1, we have
$$
\|\nabla^s_{\mathbb{R}^n}\tilde{u}\|_{L^p(B_1)}\leq C\Big[\|\tilde{f}\|_{L^p(B_8)}+\Big(\int_{\mathbb{R}^n\setminus B_8}\frac{|\tilde{f}|^2(y)}{|y|^{n+s}}dy\Big)^{\frac{1}{2}}+1\Big],
$$
which implies
$$
\|\nabla^s_{\mathbb{R}^n}u\|_{L^p(B_1)}\leq C\Big[\|u\|_{L^2(B_8)}+\Big(\int_{\mathbb{R}^n\setminus B_8}\frac{u^2(x)}{|x|^{n+2s}}dx\Big)^{\frac{1}{2}}+\|f\|_{L^p(B_8)}+\Big(\int_{\mathbb{R}^n\setminus B_8}\frac{|f|^2(x)}{|x|^{n+s}}dx\Big)^{\frac{1}{2}}\Big],
$$
where the constant $C$ depends on $n$, $s$, $\lambda$, $\Lambda$, $\beta$ and $p$.   Then by fractional Poincar\'{e} inequality, we also have
\begin{eqnarray}
\|u\|_{L^p(B_1)}&\leq&\|u-\bar{u}_{B_1}\|_{L^p(B_1)}+\|\bar{u}_{B_1}\|_{L^p(B_1)}\nonumber\\
&\leq&c\|\nabla^s_{\mathbb{R}^n}u\|_{L^p(B_1)}+\bar{u}_{B_1}|B_1|^{\frac{1}{p}}\nonumber\\
&\leq&C\Big[\|u\|_{L^2(B_8)}+\Big(\int_{\mathbb{R}^n\setminus B_8}\frac{u^2(x)}{|x|^{n+2s}}dx\Big)^{\frac{1}{2}}+\|f\|_{L^p(B_8)}+\Big(\int_{\mathbb{R}^n\setminus B_8}\frac{|f|^2(x)}{|x|^{n+s}}dx\Big)^{\frac{1}{2}}\Big].\nonumber
\end{eqnarray}
Thus we obtain that $u\in H^{s,p}(B_1)$.
\end{proof}

\section*{Acknowledgments} The third author was partially supported by NSFC Grant 12571223.

\section*{Data Availability}
Data sharing is not applicable to this article as no new data were created or analyzed in this study. The research is entirely theoretical and mathematical.

\section*{Conflict of Interest}
The authors declare no competing interests.

\bibliographystyle{plain}

\end{document}